\documentclass[10.99pt]{article}
\usepackage{epic,latexsym,amssymb}
\usepackage{color}
\usepackage{tikz}
\usepackage{amsmath}
\usepackage{comment}

\newtheorem{theorem}{Theorem}
\newtheorem{lemma}{Lemma}
\newtheorem{proposition}{Proposition}
\newtheorem{problem}{Problem}
\newtheorem{observation}{Observation}
\newtheorem{claim}{Claim}
\newtheorem{subclaim}{Claim}[claim]
\newtheorem{subsubclaim}{Claim}[subclaim]

\newcommand{\QED}{$\Box$}
\newcommand{\smallqed}{{\tiny ($\Box$)}}

\newcommand{\cB}{{\cal B}}

\newcommand{\cS}{{\cal S}}
\newcommand{\cL}{{\cal L}}
\newcommand{\rdom}{{\rm rdom}}
\newcommand{\ndom}{{\rm ndom}}
\newcommand{\dom}{{\rm dom}}
\newcommand{\w}{{\rm w}}
\newcommand{\modo}{{\rm mod \,}}

\def\cp{\,\square\,}

\newcommand{\barS}{{\overline{S}}}

\newcommand{\proof}{\noindent\textbf{Proof. }}
\newcommand{\2}{ \vspace{0.2cm} }
\newcommand{\1}{ \vspace{0.1cm} }

\let\oldenumerate\enumerate
\renewcommand{\enumerate}{
  \oldenumerate
  \setlength{\itemsep}{0pt}
  \setlength{\parskip}{0pt}
  \setlength{\parsep}{0pt}
}

\begin{document}

\title{Best possible upper bounds on the \\
restrained domination number of cubic graphs}
\author{$^{1,2}$Bo\v{s}tjan Bre\v{s}ar \, and $^3$Michael A. Henning \\
\\
$^1$Faculty of Natural Sciences and Mathematics \\ University of Maribor, Slovenia \1 \\
$^2$Institute of Mathematics, Physics and Mechanics, Slovenia \\
{\tt bostjan.bresar@um.si} \\
\\
$^3$Department of Mathematics and Applied Mathematics \\
University of Johannesburg \\
Auckland Park, 2006 South Africa\\
{\tt mahenning@uj.ac.za} }

\date{}
\maketitle

\begin{abstract}
A dominating set in a graph $G$ is a set $S$ of vertices such that every vertex in $V(G) \setminus S$ is adjacent to a vertex in $S$. A restrained dominating set of $G$ is a dominating set $S$ with the additional restraint that the graph $G - S$ obtained by removing all vertices in $S$ is isolate-free.
The domination number $\gamma(G)$ and the restrained domination number $\gamma_{r}(G)$ are the minimum cardinalities of a dominating set and restrained dominating set, respectively, of $G$. Let $G$ be a cubic graph of order~$n$.
A classical result of Reed [Combin. Probab. Comput. 5 (1996), 277--295]
states that $\gamma(G) \le \frac{3}{8}n$, and this bound is best possible. 
To determine a best possible upper bound on the restrained domination number of $G$ is more challenging, and
we prove that $\gamma_{r}(G) \le \frac{2}{5}n$.
\end{abstract}


\section{Introduction}
\label{Intro}

A \emph{dominating set} of a graph $G$ is a set $S$ of vertices of $G$ such that every vertex not in $S$ has a neighbor in $S$, where two vertices are neighbors in $G$ if they are adjacent. The \emph{domination number} of $G$, denoted by $\gamma(G)$, is the minimum cardinality of a dominating set of $G$. A set $S$ dominates a vertex~$v$ is $v \in S$ or if $v$ has a neighbor in $S$. A \emph{restrained dominating set}, abbreviated RD-set, of $G$ is a dominating set $S$ of $G$ with the additional property that every vertex not in $S$ has a neighbor not in $S$, that is, the subgraph of $G$ induced by the set $V(G) \setminus S$ is isolate-free. The \emph{restrained domination number} of $G$, denoted by $\gamma_{r}(G)$, is the minimum cardinality of a RD-set of $G$. A $\gamma_{r}$-\emph{set of $G$} is a RD-set of $G$ of minimum cardinality~$\gamma_{r}(G)$. Restrained domination in graphs is well studied in the literature  with over 100 publications according to MathSciNet. We refer the reader to the excellent book chapter by Hattingh and Joubert in 2020 on restrained domination in graphs that gives the state of the art on the topic. For recent books on domination in graphs, we refer the reader to~\cite{HaHeHe-20,HaHeHe-21,HaHeHe-22,HeYebook-13}.

A \emph{cubic graph}, also called a $3$-\emph{regular graph}, is a graph in which every vertex has degree~$3$. A \emph{subcubic graph} is a graph with maximum degree at most~$3$. Domination in cubic and subcubic graphs is very well studied in the literature (see, for example,~\cite{AnKu-22,AbHe-22,BrHe-19,Cyman-18,DaHe-19,Dorbec-15,GoHe-09,HaJo-11,HaJo-20,HeKa-18,KoSt-09,Kral-12,RaRe-08,SoHe-13,St-08,OWe-16,Re-96}).
We define a \emph{special subcubic graph} as a subcubic graph $G$ with minimum degree at least~$2$. In this paper, we continue the study of restrained domination in cubic graphs. We consider the following problem.

\begin{problem}
\label{problem-1}
Determine the best possible constant $c_{\rdom}$ such that $\gamma_r(G) \le c_{\rdom} \cdot n(G)$ for all cubic graphs~$G$.
\end{problem}

\begin{figure}[htb]
\begin{center}
\begin{tikzpicture}[scale=.75,style=thick,x=0.75cm,y=0.75cm]
\def\vr{2.5pt}
\path (0.00,1.49) coordinate (x1);
\path (1.27,2.41) coordinate (x2);
\path (2.54,1.49) coordinate (x3);
\path (2.05,0.00) coordinate (x4);
\path (0.48,0.00) coordinate (x5);
\path (0.63,1.28) coordinate (y1);
\path (1.27,1.75) coordinate (y2);
\path (1.90,1.28) coordinate (y3);
\path (1.66,0.54) coordinate (y4);
\path (0.88,0.54) coordinate (y5);
\draw (x1)--(x2)--(x3)--(x4)--(x5)--(x1);
\draw (y1)--(y3)--(y5)--(y2)--(y4)--(y1);
\draw (x1)--(y1);
\draw (x2)--(y2);
\draw (x3)--(y3);
\draw (x4)--(y4);
\draw (x5)--(y5);
\draw (x1) [fill=black] circle (\vr);
\draw (x2) [fill=black] circle (\vr);
\draw (x3) [fill=black] circle (\vr);
\draw (x4) [fill=white] circle (\vr);
\draw (x5) [fill=white] circle (\vr);
\draw (y1) [fill=white] circle (\vr);
\draw (y2) [fill=black] circle (\vr);
\draw (y3) [fill=white] circle (\vr);
\draw (y4) [fill=white] circle (\vr);
\draw (y5) [fill=white] circle (\vr);
\end{tikzpicture}
\end{center}
\vskip -0.45 cm
\caption{The Petersen graph $G$}
\label{Petersen}
\end{figure}
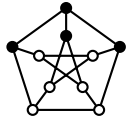

The best known upper bound to date, prior to this paper, on $c_{\rdom}$ is due to Hattingh and Joubert~\cite{HaJo-11} who proved that $c_{\rdom} \le \frac{5}{11}$. Their proof is nontrivial and uses intricate and ingenious counting arguments. We observe that the Petersen graph $G$, illustrated in Figure~\ref{Petersen}, has order~$n(G) = 10$ and $\gamma_r(G) = 4 = \frac{2}{5}n(G)$, where the set consisting of the four shaded vertices is an example of a $\gamma_r$-set of~$G$. This yields the trivial lower bound $c_{\rdom} \ge \frac{2}{5}$.

\begin{theorem}{\rm (\cite{HaJo-11})}
\label{rdom:main-known-HJ}
$\frac{2}{5} \le c_{\rdom} \le \frac{5}{11}$.
\end{theorem}

It is conjectured in~\cite{He-22} that the lower bound in Theorem~\ref{rdom:main-known-HJ} is the correct value of $c_{\rdom}$. In this paper we prove that this is indeed the case.

\begin{theorem}
\label{rdom:main-1}
$c_{\rdom} = \frac{2}{5}$.
\end{theorem}

To prove Theorem~\ref{rdom:main-1}, it suffices to show that if $G$ is a cubic graph of order~$n$, then $\gamma_r(G) \le \frac{2}{5}n$. However in order to prove this result, we relax the $3$-regularity condition to allow vertices of degree~$2$ in the mix to make the inductive hypothesis easier to handle. If $n_2(G)$ and $n_3(G)$ denote the number of vertices of degree~$2$ and~$3$, respectively, in such a graph $G$, then we would like to prove that $10\gamma_r(G) \le 5n_2(G) + 4n_3(G)$ since if $G$ is $3$-regular this yields $\gamma_r(G) \le \frac{2}{5}n$. However, relaxing the $3$-regularity condition results in a family $\cB_{\rdom}$ of ``troublesome graphs'' for which the desired inequality $10\gamma_r(G) \le 5n_2(G) + 4n_3(G)$ does not hold. Therefore we add a function $\Omega(G)$ such that the statement becomes true even for these troublesome graphs. However we try to keep $\Omega(G)$ as small as possible in order to establish a bound on $\gamma_r(G)$ that remains as strong as possible. The resulting bound will be our key result that will enable us to prove Theorem~\ref{rdom:main-1}.

We proceed as follows. In Section~\ref{S:key} we formally state our key result, namely Theorem~\ref{thm:main-1}. In Section~\ref{S:notation} we present the necessary graph theory notation. In Section~\ref{S:nearRDset} we introduce the concept of a near restrained dominating sets, which we will need when proving our key result. Known results are discussed in Section~\ref{S:known}. In Section~\ref{S:cB} we discuss properties of troublesome graphs that belong to the family $\cB_{\rdom}$. A preliminary result is proven in Section~\ref{S:prelim}. A proof of our key result is given in Section~\ref{S:proof-key}, and thereafter in Section~\ref{S:proof-main} we deduce our main result.

\section{Key result}
\label{S:key}

In order to prove our main result, namely Theorem~\ref{rdom:main-1}, we identify a family $\cB_{\rdom} = \{R_1, R_2, \ldots, R_{11}\}$ of eleven troublesome graphs $G$ shown in Figure~\ref{rdom:fig-2} that satisfy $10\gamma_r(G) > 5n_2(G) + 4n_3(G)$. We note that $R_6$ is obtained from $R_4$ by selecting two arbitrary vertices of degree~$2$ at distance~$2$, adding a new vertex and joining it to the two selected vertices of degree~$2$. Let {$\cB_{\rdom,1} = \{R_6,R_7,R_8,R_9,R_{11}\}$, $\cB_{\rdom,2} = \{R_2,R_3\}$, $\cB_{\rdom,3} = \{R_{10}\}$, $\cB_{\rdom,4} = \{R_4,R_5\}$ and $\cB_{\rdom,5} = \{R_1\}$. Let $f_i(G)$ denote the number of components of a special subcubic graph $G$ that belong to $\cB_{\rdom,i}$ for $i \in [5]$. We define
\[
\Omega(G) = \sum_{i=1}^5 i f_i(G).
\]
We note that if $G$ is a connected graph and $G \notin \cB_{\rdom}$, then $\Omega(G) = 0$, while if $G \in \cB_{\rdom}$, then $G \in \cB_{\rdom,i}$ for some $i \in [5]$ in which case $\Omega(G) = i \le 5$. We define a weight function $\w(G)$ associated with $G$ by
\[
\w(G) = 5n_2(G) + 4n_3(G) + \Omega(G).
\]

We define the \emph{weight} $\w_G(v)$ of a vertex $v$ in $G$ as its contribution to the weight $5n_2(G) + 4n_3(G)$. Thus, if $\deg_G(v) = 2$, then $\w_G(v) = 5$, and if $\deg_G(v) = 3$, then $\w_G(v) = 4$. We define the \emph{weight} $\w_G(S)$ of a set $S$ of vertices in $G$ as the sum of the weights of vertices in $S$, that is, $\w_G(S) = \sum_{v \in S} \w_G(v)$. We are now in a position to state our key result, a proof of which will be given in Section~\ref{S:proof-key}.

\begin{theorem}
\label{thm:main-1}
If $G$ is a special subcubic graph, then $10\gamma_r(G) \le \w(G)$.
\end{theorem}

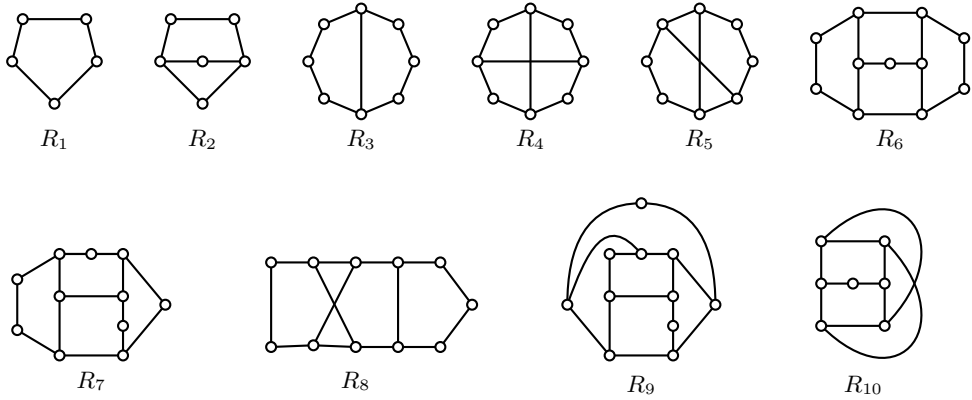
\begin{figure}[htb]
\begin{center}
\begin{tikzpicture}[scale=.75,style=thick,x=0.75cm,y=0.75cm]
\def\vr{2.5pt} 
\path (0,1.25) coordinate (v1);
\path (1,0.25) coordinate (v2);
\path (2,1.25) coordinate (v3);
\path (1.75,2.25) coordinate (v4);
\path (0.25,2.25) coordinate (v5);
\draw (v1)--(v2)--(v3)--(v4)--(v5)--(v1);
\draw (v1) [fill=white] circle (\vr);
\draw (v2) [fill=white] circle (\vr);
\draw (v3) [fill=white] circle (\vr);
\draw (v4) [fill=white] circle (\vr);
\draw (v5) [fill=white] circle (\vr);
\draw (1,-0.6) node {{\small $R_1$}};
\path (3.5,1.25) coordinate (w1);
\path (4.5,0.25) coordinate (w2);
\path (5.5,1.25) coordinate (w3);
\path (5.25,2.25) coordinate (w4);
\path (3.75,2.25) coordinate (w5);
\path (4.5,1.25) coordinate (w6);
\draw (w1)--(w2)--(w3)--(w4)--(w5)--(w1);
\draw (w1)--(w6)--(w3);
\draw (w1) [fill=white] circle (\vr);
\draw (w2) [fill=white] circle (\vr);
\draw (w3) [fill=white] circle (\vr);
\draw (w4) [fill=white] circle (\vr);
\draw (w5) [fill=white] circle (\vr);
\draw (w6) [fill=white] circle (\vr);
\draw (4.5,-0.6) node {{\small $R_2$}};
\path (8.25,0.00) coordinate (v1);
\path (9.13,0.37) coordinate (v2);
\path (9.50,1.25) coordinate (v3);
\path (9.13,2.13) coordinate (v4);
\path (8.25,2.50) coordinate (v5);
\path (7.37,2.13) coordinate (v6);
\path (7.00,1.25) coordinate (v7);
\path (7.37,0.37) coordinate (v8);
\draw (v1)--(v2)--(v3)--(v4)--(v5)--(v6)--(v7)--(v8)--(v1);
\draw (v1)--(v5);
\draw (v1) [fill=white] circle (\vr);
\draw (v2) [fill=white] circle (\vr);
\draw (v3) [fill=white] circle (\vr);
\draw (v4) [fill=white] circle (\vr);
\draw (v5) [fill=white] circle (\vr);
\draw (v6) [fill=white] circle (\vr);
\draw (v7) [fill=white] circle (\vr);
\draw (v8) [fill=white] circle (\vr);
\draw (8.25,-0.6) node {{\small $R_3$}};
\path (12.25,0.00) coordinate (w1);
\path (13.13,0.37) coordinate (w2);
\path (13.50,1.25) coordinate (w3);
\path (13.13,2.13) coordinate (w4);
\path (12.25,2.50) coordinate (w5);
\path (11.37,2.13) coordinate (w6);
\path (11.00,1.25) coordinate (w7);
\path (11.37,0.37) coordinate (w8);
\draw (w1)--(w2)--(w3)--(w4)--(w5)--(w6)--(w7)--(w8)--(w1);
\draw (w1)--(w5);
\draw (w3)--(w7);
\draw (w1) [fill=white] circle (\vr);
\draw (w2) [fill=white] circle (\vr);
\draw (w3) [fill=white] circle (\vr);
\draw (w4) [fill=white] circle (\vr);
\draw (w5) [fill=white] circle (\vr);
\draw (w6) [fill=white] circle (\vr);
\draw (w7) [fill=white] circle (\vr);
\draw (w8) [fill=white] circle (\vr);
\draw (12.25,-0.6) node {{\small $R_4$}};
\path (16.25,0.00) coordinate (x1);
\path (17.13,0.37) coordinate (x2);
\path (17.50,1.25) coordinate (x3);
\path (17.13,2.13) coordinate (x4);
\path (16.25,2.50) coordinate (x5);
\path (15.37,2.13) coordinate (x6);
\path (15.00,1.25) coordinate (x7);
\path (15.37,0.37) coordinate (x8);
\draw (x1)--(x2)--(x3)--(x4)--(x5)--(x6)--(x7)--(x8)--(x1);
\draw (x1)--(x5);
\draw (x2)--(x6);
\draw (x1) [fill=white] circle (\vr);
\draw (x2) [fill=white] circle (\vr);
\draw (x3) [fill=white] circle (\vr);
\draw (x4) [fill=white] circle (\vr);
\draw (x5) [fill=white] circle (\vr);
\draw (x6) [fill=white] circle (\vr);
\draw (x7) [fill=white] circle (\vr);
\draw (x8) [fill=white] circle (\vr);
\draw (16.25,-0.6) node {{\small $R_5$}};

\path (19,0) coordinate (u1);
\path (19,2.4) coordinate (u2);
\path (19.75,0.5) coordinate (u3);
\path (20.5,1) coordinate (u4);
\path (20.5,1.7) coordinate (u5);
\path (20.5,2.4) coordinate (u6);
\path (21.25,0.5) coordinate (u7);
\path (22,0) coordinate (u8);
\path (22,2.4) coordinate (u9);
\draw (u1)--(u2)--(u6)--(u9)--(u8)--(u1);
\draw (u2)--(u5)--(u9);
\draw (u1)--(u3)--(u4)--(u7)--(u8);
\draw (u4)--(u5);
\draw (u1) [fill=white] circle (\vr);
\draw (u2) [fill=white] circle (\vr);
\draw (u3) [fill=white] circle (\vr);
\draw (u4) [fill=white] circle (\vr);
\draw (u5) [fill=white] circle (\vr);
\draw (u6) [fill=white] circle (\vr);
\draw (u7) [fill=white] circle (\vr);
\draw (u8) [fill=white] circle (\vr);
\draw (u9) [fill=white] circle (\vr);
\draw (20.5,-0.6) node {{\small $R_6$}};
\draw[anchor = south] (u6) node {{\small $v$}};

\end{tikzpicture}


\begin{tikzpicture}[scale=.75,style=thick,x=0.75cm,y=0.75cm]
\def\vr{2.5pt} 
\path (-5.5,0.6) coordinate (u1);
\path (-5.5,1.8) coordinate (u2);
\path (-4.5,0) coordinate (u3);
\path (-4.5,1.2) coordinate (u4);
\path (-4.5,2.4) coordinate (u5);
\path (-3.75,1.2) coordinate (u6);
\path (-3,0) coordinate (u7);
\path (-3,1.2) coordinate (u8);
\path (-3,2.4) coordinate (u9);
\path (-2,0.6) coordinate (u10);
\path (-2,1.8) coordinate (u11);
\draw (u1)--(u2);
\draw (u3)--(u4)--(u5);
\draw (u7)--(u8)--(u9);
\draw (u10)--(u11);
\draw (u1)--(u3)--(u7)--(u10);
\draw (u2)--(u5)--(u9)--(u11);
\draw (u4)--(u6)--(u8);
\draw (u1) [fill=white] circle (\vr);
\draw (u2) [fill=white] circle (\vr);
\draw (u3) [fill=white] circle (\vr);
\draw (u4) [fill=white] circle (\vr);
\draw (u5) [fill=white] circle (\vr);
\draw (u6) [fill=white] circle (\vr);
\draw (u7) [fill=white] circle (\vr);
\draw (u8) [fill=white] circle (\vr);
\draw (u9) [fill=white] circle (\vr);
\draw (u10) [fill=white] circle (\vr);
\draw (u11) [fill=white] circle (\vr);
\draw (-3.75,-0.6) node {{\small $R_7$}};
\path (0,0.6) coordinate (x1);
\path (0,1.8) coordinate (x2);
\path (1,0) coordinate (x3);
\path (1,1.4) coordinate (x4);
\path (1,2.4) coordinate (x5);
\path (1.75,2.4) coordinate (x6);
\path (2.5,0) coordinate (x7);
\path (2.5,1.4) coordinate (x8);
\path (2.5,2.4) coordinate (x9);
\path (3.5,1.2) coordinate (x10);
\path (2.5,0.7) coordinate (x11);
\draw (x1)--(x2);
\draw (x3)--(x4)--(x5);
\draw (x7)--(x11)--(x8)--(x9);
\draw (x1)--(x3)--(x7)--(x10);
\draw (x2)--(x5)--(x9)--(x10);
\draw (x5)--(x6)--(x9);
\draw (x4)--(x8);
\draw (x1) [fill=white] circle (\vr);
\draw (x2) [fill=white] circle (\vr);
\draw (x3) [fill=white] circle (\vr);
\draw (x4) [fill=white] circle (\vr);
\draw (x5) [fill=white] circle (\vr);
\draw (x6) [fill=white] circle (\vr);
\draw (x7) [fill=white] circle (\vr);
\draw (x8) [fill=white] circle (\vr);
\draw (x9) [fill=white] circle (\vr);
\draw (x10) [fill=white] circle (\vr);
\draw (x11) [fill=white] circle (\vr);
\draw (1.75,-0.6) node {{\small $R_8$}};
\path (5,0.2) coordinate (v1);
\path (5,2.2) coordinate (v2);
\path (6,0.24) coordinate (v3);
\path (6,2.2) coordinate (v4);
\path (7,0.2) coordinate (v5);
\path (7,2.2) coordinate (v6);
\path (8,0.2) coordinate (v7);
\path (8,2.2) coordinate (v8);
\path (9,0.2) coordinate (v9);
\path (9,2.2) coordinate (v10);
\path (9,1.2) coordinate (v11);
\draw (v1)--(v2);
\draw (v3)--(v6);
\draw (v4)--(v5);
\draw (v7)--(v8);
\draw (v1)--(v3)--(v5)--(v7)--(v9)--(v11);
\draw (v2)--(v4)--(v6)--(v8)--(v10)--(v11);
\draw (v1) [fill=white] circle (\vr);
\draw (v2) [fill=white] circle (\vr);
\draw (v3) [fill=white] circle (\vr);
\draw (v4) [fill=white] circle (\vr);
\draw (v5) [fill=white] circle (\vr);
\draw (v6) [fill=white] circle (\vr);
\draw (v7) [fill=white] circle (\vr);
\draw (v8) [fill=white] circle (\vr);
\draw (v9) [fill=white] circle (\vr);
\draw (v10) [fill=white] circle (\vr);
\draw (v11) [fill=white] circle (\vr);
\draw (7,-0.6) node {{\small $R_9$}};
\path (11,1.2) coordinate (x1);
\path (12.75,3.6) coordinate (x2);
\path (12,0) coordinate (x3);
\path (12,1.4) coordinate (x4);
\path (12,2.4) coordinate (x5);
\path (12.75,2.4) coordinate (x6);
\path (13.5,0) coordinate (x7);
\path (13.5,1.4) coordinate (x8);
\path (13.5,2.4) coordinate (x9);
\path (14.5,1.2) coordinate (x10);
\path (13.5,0.7) coordinate (x11);
\draw (x3)--(x4)--(x5);
\draw (x7)--(x11)--(x8)--(x9);
\draw (x1)--(x3)--(x7)--(x10);
\draw (x5)--(x9)--(x10);
\draw (x5)--(x6)--(x9);
\draw (x4)--(x8);
\draw (x1) to[out=90,in=180, distance=1cm] (x2);
\draw (x1) to[out=70,in=130, distance=1cm] (x6);
\draw (x2) to[out=0,in=90, distance=1cm] (x10);
\draw (x1) [fill=white] circle (\vr);
\draw (x2) [fill=white] circle (\vr);
\draw (x3) [fill=white] circle (\vr);
\draw (x4) [fill=white] circle (\vr);
\draw (x5) [fill=white] circle (\vr);
\draw (x6) [fill=white] circle (\vr);
\draw (x7) [fill=white] circle (\vr);
\draw (x8) [fill=white] circle (\vr);
\draw (x9) [fill=white] circle (\vr);
\draw (x10) [fill=white] circle (\vr);
\draw (x11) [fill=white] circle (\vr);
\draw (12.75,-0.7) node {{\small $R_{10}$}};
\path (16.5,0.7) coordinate (f3);
\path (16.5,1.7) coordinate (f4);
\path (16.5,2.7) coordinate (f5);
\path (17.25,1.7) coordinate (f6);
\path (18,0.7) coordinate (f7);
\path (18,1.7) coordinate (f8);
\path (18,2.7) coordinate (f9);
\draw (f3)--(f4)--(f5);
\draw (f7)--(f8)--(f9);
\draw (f3)--(f7);
\draw (f5)--(f9);
\draw (f4)--(f6)--(f8);
\draw (f5) to[out=45,in=45, distance=2cm] (f7);
\draw (f3) to[out=-45,in=-45, distance=2cm] (f9);
\draw (f3) [fill=white] circle (\vr);
\draw (f4) [fill=white] circle (\vr);
\draw (f5) [fill=white] circle (\vr);
\draw (f6) [fill=white] circle (\vr);
\draw (f7) [fill=white] circle (\vr);
\draw (f8) [fill=white] circle (\vr);
\draw (f9) [fill=white] circle (\vr);
\draw (17.5,-0.7) node {{\small $R_{11}$}};
\end{tikzpicture}
\end{center}
\begin{center}
\vskip -0.70 cm
\caption{The family $\cB_{\rdom}$}
\label{rdom:fig-2}
\end{center}
\end{figure}

\subsection{Notation}
\label{S:notation}

For notation and graph theory terminology, we in general follow~\cite{HeYebook-13}. Specifically, let $G$ be a graph with vertex set $V(G)$ and edge set $E(G)$, and of order~$n(G) = |V(G)|$ and size $m(G) = |E(G)|$. For a set of vertices $S\subseteq V(G)$, the subgraph induced by $S$ is denoted by $G[S]$. Two vertices in $G$ are \emph{neighbors} if they are adjacent. The \emph{open neighborhood} $N_G(v)$ of a vertex $v$ in $G$ is the set of neighbors of $v$, while the \emph{closed neighborhood} of $v$ is the set $N_G[v] = \{v\} \cup N(v)$. Two vertices are \emph{open twins} if they have the same open neighborhood. We denote the \emph{degree} of $v$ in $G$ by $\deg_G(v) = |N_G(v)|$. The minimum and maximum degree in $G$ is denoted by $\delta(G)$ and $\Delta(G)$, respectively. An \emph{isolated vertex} is a vertex of degree~$0$. A graph is \emph{isolate}-\emph{free} if it contains no isolated vertex.

We denote a \emph{path}, a \emph{cycle}, and a \emph{complete} graph on $n$ vertices by $P_n$, $C_n$, and $K_n$, respectively. A \emph{diamond} is the graph $K_4 - e$ where $e$ is an arbitrary edge of the $K_4$. A \emph{domino} is a graph that can be obtained from a $6$-cycle by adding an edge between two antipodal vertices of the $6$-cycle. An $F$-\emph{component} of a graph $G$ is a component of $G$ that is isomorphic to $F$. An edge-cut of a connected graph is a set of edges whose removal disconnects the graph. A $k$-\emph{edge}-\emph{cut} is an edge-cut of cardinality~$k$. The \emph{girth} of $G$ is the length of a shortest cycle in $G$.

If $G$ is a special subcubic graph, then we denote by $n_2(G)$ and $n_3(G)$ the number of vertices of degree~$2$ and~$3$, respectively, in $G$. For a special subcubic graph $G$, let $\cS$ and $\cL$ be the set of all vertices of degree~$2$ and~$3$ in $G$, respectively, that is, $\cL = \{v \in V(G) \, \colon \deg_G(v) = 3\}$ and $\cS = \{v \in V(G) \, \colon \deg_G(v) = 2\}$. We call a vertex in $\cL$ a \emph{large vertex}, and a vertex in $\cS$ a \emph{small vertex}. For $k \ge 3$, we define a \emph{$k$-handle} to be a $k$-cycle that contains exactly one large vertex. For $k \ge 1$, a \emph{$k$-linkage} is a path on $k+2$ vertices that starts and ends at distinct large vertices and with $k$ internal vertices of degree~$2$ in $G$. A \emph{handle} is a $k$-handle for some $k \ge 3$, and a \emph{linkage} is a $k$-linkage for some $k \ge 1$. We use the standard notation $[k] = \{1,\ldots,k\}$.

\subsection{Near restrained dominating sets}
\label{S:nearRDset}

In order to prove our main result, we introduce the concept of a near restrained dominating set. Given a graph $G$ and a set $S$ of vertices in $G$, we let $\barS$ denote the complement of $S$, that is, $\barS = V(G) \setminus S$. We define a \emph{near restrained dominating set}, abbreviated NeRD-set, of $G$ with respect to a subset $X$ of vertices of $G$ as a relaxed variant of a RD-set $S$ of $G$ such that either the vertices in $X$ need not be dominated by $S$ but every vertex in $\barS$ is still required to have a neighbor in $\barS$ or the vertices in $X$ are dominated by $S$ but need not have a neighbor in~$\barS$. 
Formally, a NeRD-set of $G$ with respect to a specified subset $X$ is a set $S \subseteq V(G)$ such that exactly one of the following two conditions hold:
\\[-20pt]
\begin{enumerate}
\item[{\rm (C1)}] The set $S$ dominates the set $V(G) \setminus X$ and every vertex in $\barS$ has a neighbor in $\barS$.
\item[{\rm (C2)}] The set $S$ dominates the set $V(G)$, the set $X \subseteq \barS$, and every vertex in $\barS \setminus X$ has a neighbor in~$\barS$.
\end{enumerate}

If condition (C1) holds, then we refer to the NeRD-set as a type-$1$ such set, while if condition (C2) holds, then we refer to the NeRD-set as a type-$2$ such set. We denote by $\gamma_{r,\ndom}(G;X)$ the minimum cardinality of a type-$1$ NeRD-set with respect to the set $X$ (where ``$\ndom$'' stands for ``not dominated'' since the vertices in $X$ are not required to be dominated), and we denote by $\gamma_{r,\dom}(G;X)$ the minimum cardinality of a type-$2$ NeRD-set with respect to the set $X$ (where ``$\dom$'' stands for ``dominated'' since the vertices in $X$ are dominated but not required to have a neighbor that is not dominated). If $X = \{v\}$, we simply write $\gamma_{r,\ndom}(G;v)$ and $\gamma_{r,\dom}(G;v)$ rather than $\gamma_{r,\ndom}(G;\{v\})$ and $\gamma_{r,\dom}(G;\{v\})$, respectively. Since every RD-set is also an NeRD-set, we note that $\gamma_{r,\ndom}(G;X) \le \gamma_{r}(G)$ and $\gamma_{r,\dom}(G;X) \le \gamma_{r}(G)$.

\subsection{Known bounds on restrained domination}
\label{S:known}

Closed formulas for the restrained domination number of paths and cycles are given in~\cite{Domke-99} where it is shown that for $n \ge 1$, $\gamma_r(P_n) = n - 2\lfloor \frac{n-1}{3} \rfloor$ and for $n \ge 3$, $\gamma_r(C_n) = n - 2\lfloor \frac{n}{3} \rfloor$.
The following theorem summarizes classical results on bounds on the restrained domination number of a graph.

\begin{theorem}
\label{rdom:known}
If $G$ is a connected graph of order~$n$, then the following hold. \\[-22pt]
\begin{enumerate}
\item[{\rm (a)}] {\rm (\cite{Domke-99})} If $\delta(G) \ge 1$, then $\gamma_r(G) \le n - 2$, unless $G$ is a star $K_{1,n-1}$, in which case $\gamma_r(G) = n$. \1
\item[{\rm (b)}] {\rm (\cite{Domke-20})} If $\delta(G) \ge 2$ and $G \ne C_5$, then $\gamma_r(G) \le \frac{1}{2}n$.  \1
\item[{\rm (c)}] {\rm (\cite{HaJo-11})} If $G$ is a cubic graph, then $\gamma_r(G) \le \frac{5}{11}n$.
\end{enumerate}
\end{theorem}

\section{Properties of graph in the family $\cB_{\rdom}$}
\label{S:cB}

In this section, we present properties of graphs that belong to the family $\cB_{\rdom} = \{R_1, \ldots, R_{11}\}$. We note that there are no open twins in the graphs in the family $\cB_{\rdom}$ with the exception of $R_2$ which contains two vertices of degree~$2$ that have two common neighbors (of degree~$3$). We shall need the following properties of graphs in the family $\cB_{\rdom}$. These properties are straightforward to check (or can be checked by computer).

\begin{observation}
\label{obser-1}
If $G \in \cB_{\rdom}$ and $v$ is a vertex of degree~$2$ in $G$, then the following properties hold. \\[-22pt]
\begin{enumerate}
\item[{\rm (a)}] $\gamma_r(R_i) = 3$ for $i \in \{1,2,11\}$, $\gamma_r(R_i) = 4$ for $i \in \{3,4,5,6\}$, and $\gamma_r(R_i) = 5$ for $i \in \{7,8,9,10\}$.
\item[{\rm (b)}] There exists a $\gamma_r$-set of $G$ that contains~$v$.
\item[{\rm (c)}] There exists a $\gamma_r$-set of $G$ that does not contains~$v$.
\item[{\rm (d)}] $\gamma_{r,\ndom}(G;v) \le \gamma_{r}(G) - 1$.
\item[{\rm (e)}] $\gamma_{r,\dom}(G;v) \le \gamma_{r}(G) - 1$, unless $v$ is an open twin of $R_2$ or $v$ belongs to the $4$-cycle in $R_6$. 
\item[{\rm (f)}] If $X$ consists of two vertices of degree~$2$, then $\gamma_{r,\dom}(G;X) \le \gamma_{r}(G) - 1$.
\end{enumerate}
\end{observation}

\begin{observation}
\label{obser-2}
Let $G \in \cB_{\rdom}$ and let $e = xy$ be an arbitrary edge of $G$. If $G^*$ is obtained from $G$ by subdividing the edge $e$ resulting in a new vertex $v^*$ of degree~$2$ (with neighbors $x$ and $y$), then $\gamma_r(G^*) \le \gamma_r(G)$. Furthermore, there exists a $\gamma_r$-set of $G^*$ that contains~$v^*$ and contains neither $x$ nor~$y$.
\end{observation}

\begin{observation}
\label{obser-3}
Let $G \in \cB_{\rdom}$ and let $e = xy$ be an arbitrary edge of $G$. If $G^*$ is obtained from $G$ by subdividing the edge $e$ twice resulting in a path $x x_1 y_1 y$, then $\gamma_r(G^*) \le \gamma_r(G)$. Furthermore, there exists a $\gamma_r$-set of $G^*$ that contains~$x_1$ but not $y_1$.
\end{observation}

\begin{observation}
\label{obser-4}
Let $G \in \cB_{\rdom}$ and let $e = xy$ be an arbitrary edge of $G$. If $G^*$ is obtained from $G$ by subdividing the edge $e$ three times resulting in a path $x v_1 v_2 v_3 y$, then the following properties hold.
\\[-22pt]
\begin{enumerate}
\item[{\rm (a)}] $\gamma_{r,\ndom}(G^*;v_1) \le \gamma_{r}(G)$.
\item[{\rm (b)}] $\gamma_{r,\dom}(G^*;v_1) \le \gamma_{r}(G)$, except if $e$ is an edge that belongs to the $4$-cycle in $R_6$.
\item[{\rm (c)}] If $G \in \{R_4,R_5,R_{10}\}$, then $\gamma_{r,\dom}(G^*;v_2) \le \gamma_{r}(G)$.
\end{enumerate}
\end{observation}

\begin{observation}
\label{obser-5}
Let $G \in \cB_{\rdom}$ and let $e = xy$ be an arbitrary edge of $G$. If $G^*$ is obtained from $G$ by subdividing the edge $e$ four times resulting in a path $x v_1 v_2 v_3 v_4 y$, then there exists a RD-set $S^*$ of $G^*$ such that $S^* \cap \{v_1,v_2,v_3,v_4\} = \{v_1,v_4\}$ and the following properties hold.
\\[-22pt]
\begin{enumerate}
\item[{\rm (a)}] If $G \notin \{R_4,R_5\}$, then $|S^*| \le \gamma_{r}(G) + 1$.
\item[{\rm (b)}] If $G \in \{R_4,R_5\}$, then $|S^*| \le \gamma_{r}(G)$.
\end{enumerate}
\end{observation}

\begin{observation}
\label{obser-6}
Let $G \in \cB_{\rdom}$ and let $e = xy$ be an arbitrary edge of $G$. If $G^*$ is obtained from $G$ by subdividing the edge $e$ four times resulting in a path $x v_1 v_2 v_3 v_4 y$, then there exists a RD-set $S^*$ of $G^*$ such that $v_2 \in S^*$ and the following properties hold.
\\[-22pt]
\begin{enumerate}
\item[{\rm (a)}] If $G \ne R_2$ or if $G = R_2$ and neither $x$ nor $y$ is an open twin in $G$, then $|S^*| \le \gamma_{r}(G)$.
\item[{\rm (b)}] If $G = R_2$ and $x$ or $y$ is an open twin in $G$, then $|S^*| \le \gamma_{r}(G) + 1$.
\end{enumerate}
\end{observation}

\section{Preliminary result}
\label{S:prelim}

In this section we present a preliminary result that we will need when proving our main result.

\begin{lemma}
\label{lem:1}
If $G$ is a bipartite special subcubic graph with partite sets $\cS$ and $\cL$, then $\gamma_r(G) \le |\cL|$.
\end{lemma}
\proof Let $G$ be a bipartite subcubic graph with partite sets $\cS$ and $\cL$. Thus $\cS$ and $\cL$ are independent sets, and every vertex in $\cS$ has degree~$2$ with two neighbors in $\cL$ and every vertex in $\cL$ has degree~$3$ with three neighbors in $\cS$. Let $s = |\cS|$ and $\ell = |\cL|$.

Let $F$ be the graph with $V(F) = \cL$, where two vertices are adjacent in $F$ if and only if they have a common neighbor (that belongs to $\cS$) in the graph $G$. Let $\cL_1$ be a maximal independent set in $F$, and let $\cL_2 = \cL \setminus \cL_1$. Let $\ell_1 = |\cL_1|$ and let $\ell_2 = |\cL_2|$. Let $\cS_1$ be the set of vertices dominated by $\cL_1$ in the graph $G$, and let $\cS_2 = \cS \setminus \cS_1$. Possibly, $\cS_ 2 = \emptyset$.

If a vertex in $\cS_1$ has both its neighbors in $\cL_1$, then the set $\cL_1$ would contain two adjacent vertices in $F$, contradicting the fact that $\cL_1$ is an independent set in $F$. Hence every vertex in $\cS_1$ is adjacent to exactly one vertex of $\cL_1$ and to exactly one vertex in $\cL_2$. In particular, this implies that the subgraph $G[\cL_1 \cup \cS_1]$ of $G$ induced by the set $\cL_1 \cup \cS_1$ consists of $\ell_1$ vertex disjoint copies of $K_{1,3}$ where the central vertex of each star belongs to $\cL_1$.

By the maximality of the independent set $\cL_1$, the set $\cL_1$ is a dominating set in $F$, implying that every vertex in $\cL_2$ must have at least one neighbor in $G$ that belongs to the set $\cS_1$, that is, the set $\cS_1$ dominates the set $\cL_2$ in $G$. Let $\cL_{2.i}$ be the set of vertices in $\cL_2$ that have exactly~$i$ neighbors in $\cS_1$ for $i \in [3]$. Further, let $\ell_{2.i} = |\cL_{2.i}|$ for $i \in [3]$, and so $\ell_2 = \ell_{2.1} + \ell_{2.2} + \ell_{2.3}$.

Since each vertex in $\cS_1$ has exactly one neighbor in $\cL_2$, no two vertices in $\cL_2$ have a common neighbor in $\cS_1$. For each vertex $v$ in $\cL_{2.3}$, we select an arbitrary neighbor $v'$ in $\cS_1$ and let $\cS_{1.1}$ be the resulting subset of vertices in $\cS_1$, that is,
\[
\cS_{1.1} = \bigcup_{v \in \cL_{2.3}} \{v'\}.
\]

By our earlier observations, $|\cS_{1.1}| = \ell_{2,3}$. Let $\cS_{1.2} = \cS_1 \setminus \cS_{1.1}$. Each vertex in $\cL_{2.3}$ has one neighbor in $\cS_{1.1}$ and two neighbors in $\cS_{1.2}$, while each vertex in $\cL_{2.i}$ has $i$ neighbors in $\cS_{1.2}$ and $3-i$ neighbors in $\cS_2$ for $i \in \{1,2\}$. Each vertex in $\cL_2$ therefore has at least one neighbor in $\cS_{1.2}$, and each vertex in $\cS_{1.2}$ has exactly one neighbor in $\cL_2$. Therefore, the subgraph of $G$ induced by the set $\cS_{1.2} \cup \cL_2$ is isolate-free.

We now consider the set $D = \cL_1 \cup \cS_{1.1} \cup \cS_2$. By construction, $V(G) \setminus D = \cS_{1.2} \cup \cL_2$. As observed earlier, the subgraph of $G$ induced by the set $\cS_{1.2} \cup \cL_2$ is isolate-free. Moreover, every vertex in $\cS_{1.2}$ is dominated by the set $\cL_1 \subseteq D$ and every vertex of $\cL_2$ is dominated by the set $\cS_{1.1} \cup S_2 \subseteq D$. Hence, $D$ is indeed a RD-set. It remains for us to show that $|D| \le \ell$.
Each vertex in $\cS_2$ has no neighbor in $\cL_1 \cup \cL_{2.3}$, and therefore has both its neighbors in $\cL_{2.1} \cup \cL_{2.2}$. Counting edges between the set $\cS_2$ and the sets $\cL_{2.1} \cup \cL_{2.2}$, we therefore have $2|\cS_2| = 2\ell_{2.1} + \ell_{2.2} \le 2\ell_{2.1} + 2\ell_{2.2}$, and so $|\cS_2| \le \ell_{2.1} + \ell_{2.2}$. Recall that $|\cS_{1.1}| = \ell_{2,3}$. Hence,
$|D| = |\cL_1| + |\cS_2| + |\cS_{1.1}| \le \ell_1 + (\ell_{2.1} + \ell_{2.2}) + \ell_{2.3}  = \ell_1 + \ell_2  = \ell$, as required. Therefore, $\gamma_r(G) \le |D| \le \ell$.~\QED

\section{Proof of key result}
\label{S:proof-key}

In this section we present a proof of our key result, namely Theorem~\ref{thm:main-1}. Recall its statement.

\medskip
\noindent \textbf{Theorem~\ref{thm:main-1}}. \emph{If $G$ is a special subcubic graph, then $10\gamma_r(G) \le \w(G)$.
}

\noindent
\proof Suppose, to the contrary, that there exists a counterexample to the theorem. Among all counterexamples, let $G$ be chosen to have minimum order. Thus if $G'$ is a special subcubic graph of order less than~$n(G)$, then $G'$ is not a counterexample, that is, $10\gamma_r(G) > \w(G)$ and $10\gamma_r(G') \le \w(G')$ for all special subcubic graphs $G'$ with $n(G') < n(G)$. The restrained domination number of a graph is the sum of the restrained domination numbers of its components. Hence by the minimality of $G$, the counterexample $G$ is connected. For notational simplicity, we adopt the following notation throughout the proof. Let $n = n(G)$, $n_2 = n_2(G)$, and $n_3 = n_3(G)$. If $G'$ is a special subcubic graph, then we let $n' = n(G')$, $n_2' = n_2(G')$, and $n_3' = n_3(G')$. Further, let $k'$ be the number of components of $G'$ that belong $\cB_{\rdom}$, and let $r'$ be the remaining components of $G'$. If $G'$ is a connected graph, then we note that $k' + r' = 1$.
Since $\delta(G) \ge 2$, we note that $n \ge 3$. If $G \in \cB_{\rdom}$, then $10\gamma_r(G) = \w(G)$, contradicting the fact that $G$ is a counterexample. Hence, $G \notin \cB_{\rdom}$. If $n \in \{3,4,5\}$, then it is straightforward to check that $10\gamma_r(G) \le \w(G)$, a contradiction. Hence, $n \ge 6$.
In what follows we present a series of claims describing some structural properties of $G$ which culminate in the implication of its non-existence.

\begin{claim}
\label{claim.1}
$\Delta(G) = 3$.
\end{claim}
\proof Suppose, to the contrary, that $\Delta(G) = 2$, and so $G$ is a cycle $C_n$ (and $n \ge 6$). In this case, $\w(G) = 5n$ and $\gamma_r(C_n) = n - 2\lfloor \frac{n}{3} \rfloor$. Thus if $n \equiv 0 \, (\modo \, 3)$, then $10\gamma_r(C_n) = 10n/3$. If $n \equiv 1 \, (\modo \, 3)$, then $n \ge 7$ and $10\gamma_r(C_n) = 10(n+2)/3$. If $n \equiv 2 \, (\modo \, 3)$, then $n \ge 8$ and $10\gamma_r(C_n) = 10(n+4)/3$. In all cases, $10\gamma_r(G) \le \w(G)$, a contradiction.~\smallqed


\begin{claim}
\label{claim.longpath}
The graph $G$ does not contain a path on five vertices with the internal vertices all of degree~$2$ in $G$ and such that either the two ends of the path are not adjacent or the two ends are adjacent and both have degree~$3$ in $G$.
\end{claim}
\proof Suppose, to the contrary, that $P \colon uv_1v_2v_3w$ is a path in $G$, where $\deg_G(v_i) = 2$ for $i \in [3]$ and if $uw$ is an edge, then $\deg_G(u) = \deg_G(w) = 3$. Since $\delta(G) = 2$ and $\Delta(G) = 3$, we can choose the path $P$ so that $\deg_G(u) = 3$.
Let $G'$ be the graph of order~$n' = n - 3$ obtained from $G$ by deleting the set of vertices $\{v_1,v_2,v_3\}$. Further, if $u$ and $w$ are not adjacent, then we add the edge $uw$ to $G'$. Let $S'$ be a $\gamma_r$-set of $G'$. If $\{u,w\} \subseteq S'$, let $S = S' \cup \{v_1\}$. If $u \in S'$ and $w \notin S'$, let $S = S' \cup \{v_3\}$. If $u \notin S'$ and $w \in S'$, let $S = S' \cup \{v_1\}$. If $u \notin S'$ and $w \notin S'$, let $S = S' \cup \{v_2\}$. In all cases, $S$ is a RD-set of $G$, and so $\gamma_r(G) \le \gamma_r(G') + 1$.

Suppose that $u$ and $w$ are not adjacent in $G$. In this case the edge $uw$ was added to $G'$, implying that the degree of the vertices $u$ and $w$ remain unchanged. In particular, $\deg_{G'}(u) = 3$. The graph $G'$ is a connected special subcubic and is not a counterexample, and so $10\gamma_r(G') \le \w(G')$. Suppose that $G' \notin \cB_{\rdom}$. In this case, $\w(G) = \w(G') + 15$, and so $10\gamma_r(G) \le 10(\gamma_r(G') + 1) \le \w(G') + 10 < \w(G)$, a contradiction. Hence, $G' \in \cB_{\rdom}$. Thus, $G$ is obtained from one of the graphs in $\cB_{\rdom}$ by subdividing the (added) edge $uw$ in $G'$ three times, where as observed earlier $\deg_{G'}(u) = 3$ (and $\deg_{G'}(w) \in \{2,3\}$). Since $R_1$ has no vertex of degree~$3$, we note that $G \ne R_1$. If $G' = R_2$, then $\gamma_r(G) \le 4$ and $\w(G) = 43$. If $G' \in \{R_3,R_4,R_5\}$, then $\gamma_r(G) \le 5$ and $\w(G) \ge 51$. If $G' = R_6$, then $\gamma_r(G) \le 5$ and $\w(G) \ge 54$. If $G' \in \{R_7,R_8,R_9,R_{10}\}$, then $\gamma_r(G) \le 6$ and $\w(G) \ge 64$. If $G' = R_{11}$, then $\gamma_r(G) \le 4$ and $\w(G) = 44$. In all cases, $10\gamma_r(G) \le \w(G)$, a contradiction.

Hence, $u$ and $w$ are adjacent in $G$. As before the graph $G'$ is a connected special subcubic graph and $10\gamma_r(G') \le \w(G')$. By supposition, both $u$ and $w$ have degree~$3$ in $G$, and therefore have degree~$2$ in $G'$. Hence the weight of each of $u$ and $w$ decreases by~$1$ from weight~$5$ in $G'$ to weight~$4$ in~$G$. If $G' \notin \cB_{\rdom}$, then $\w(G) = \w(G') + 15 - 2 = \w(G') + 13$, and so $10\gamma_r(G) \le 10(\gamma_r(G') + 1) \le \w(G') + 10 < \w(G)$, a contradiction. Hence, $G' \in \cB_{\rdom}$. Thus, $G$ is obtained from one of the graphs in $\cB_{\rdom}$ by adding an extra edge between two vertices of degree~$2$ in $G'$, and then subdividing this added edge three times. Since none of $R_4$, $R_6$, $R_{10}$ and $R_{11}$ has two adjacent vertices of degree~$2$, we note that $G' \ne \{R_4,R_6,R_{10},R_{11}\}$. If $G' = R_1$, then $G = R_3$, while if $G' = R_5$, then $G = R_9$. In both cases, $G \in \cB_{\rdom}$, a contradiction. If $G' = R_2$, then $\gamma_r(G) = 4$ and $\w(G) = 41$. If $G' = R_3$, then $\gamma_r(G) = 5$ and $\w(G) = 51$. If $G' \in \{R_7,R_8,R_9\}$, then $\gamma_r(G) = 6$ and $\w(G) = 62$. In all cases, $10\gamma_r(G) \le \w(G)$, a contradiction.~\smallqed

\smallskip
As a consequence of Claim~\ref{claim.longpath}, we have the following structure of handles and linkages.

\begin{claim}
\label{claim.3}
The following properties hold in the graph $G$. \\[-22pt]
\begin{enumerate}
\item[{\rm (a)}]  If $G$ contains a $k$-handle, then $k \in \{3,4,5\}$.
\item[{\rm (b)}]  If $G$ contains a $k$-linkage, then $k \in \{1,2\}$.
\end{enumerate}
\end{claim}

\begin{claim}
\label{claim.key1}
Let $G$ be obtained from the disjoint union of a special subcubic graph $G'$ of order less than~$n$ and a graph $H$ by adding at least one edge between $H$ and $G'$. If $\gamma_r(G) \le \gamma_r(G') + p$ for some integer $p \ge 0$, then $\w(G) < \w(G') + 10p$.
\end{claim}
\proof Suppose that $\gamma_r(G) \le \gamma_r(G') + p$ for some integer $p \ge 0$. Since $G'$ is not a counterexample, no component of $G'$ is a counterexample, implying by linearity that $10\gamma_r(G') \le \w(G')$. If $\w(G) \ge \w(G') + 10p$, then $10\gamma_r(G) \le 10(\gamma_r(G') + p) \le \w(G') + 10p \le \w(G)$, a contradiction.~\smallqed

\begin{claim}
\label{claim.key2}
Let $G$ be obtained from the disjoint union of a special subcubic graph $G'$ of order less than~$n$ and a graph $H$ by adding at least one edge between $H$ and $G'$. If there exists a $\gamma_r$-set $S_H$ of $H$ such that every component of $G'$ in $\cB_{\rdom}$ has at least one neighbor that belongs to $S_H$ in the graph $G$, then $\w(G) < \w(G') + 10p$ where $p = \gamma_r(H) - k'$.
\end{claim}
\proof If $k' \ge 1$, let $G_1, \ldots, G_{k'}$ denote the component of $G'$ that belong to $\cB_{\rdom}$. By supposition, there exists a $\gamma_r$-set $S_H$ of $H$ such that the component $G_i$ contains a vertex $v_i$ that is adjacent to a vertex in $S_H$ for all $i \in [k']$. By Observation~\ref{obser-1}(d), $\gamma_{r,\ndom}(G_i;v_i) \le \gamma_r(G_i) - 1$ for all $i \in [k']$. If $G'$ has $r' \ge 1$ components that do not belong to $\cB_{\rdom}$, let $G_{k'+1}, \ldots, G_{k'+r'}$ denote these components of $G'$. Hence,
\[
\begin{array}{lcl}
\gamma_r(G) & \le & \displaystyle{ |S_H| + \left( \sum_{i=1}^{k'} \gamma_{r,\ndom}(G_i;v_i) \right) + \left( \sum_{i=k'+1}^{k'+r'} \gamma_r(G_i) \right) } \2 \\
& \le & \displaystyle{ \gamma_r(H) + \left( \sum_{i=1}^{k'+r'} \gamma_r(G_i) \right) - k'  } \2 \\
& = & \gamma_r(H) + \gamma_r(G') - k' \2 \\
& = & \gamma_r(G') + p, \1
\end{array}
\]
where $p = \gamma_r(H)  - k'$. By Claim~\ref{claim.key1}, $\w(G) < \w(G') + 10p$.~\smallqed

\begin{claim}
\label{claim.6}
There is no $3$-handle in $G$.
\end{claim}
\proof Suppose that $C \colon v v_1 v_2 v$ is a $3$-handle, where $\deg_G(v) = 3$. Let $v_3$ be the third neighbor of $v$. Suppose that $\deg_G(v_3) = 3$. Let $G' = G - V(C)$. We note that $G'$ is a connected special subcubic graph and $k' + r' = 1$. Applying Claim~\ref{claim.key2} with $H = C$ and $S_H = \{v\}$, we have $\w(G) < \w(G') + 10p$ where $p = \gamma_r(H) - k' = 1 - k'$. The weights of the vertices in $G'$ remain unchanged in $G$, except for $v_3$ whose weight increases by~$1$ from weight~$4$ in $G$ to weight $5$ in $G'$. Moreover if $k' = 1$ (that is, if $G' \in \cB_{\rdom}$), then there is an additional weight increase of at most~$5$ for creating the component $G'$ that belongs to~$\cB_{\rdom}$. Hence, $\w(G) \ge \w_G(V(C)) + ( \w(G') - 6k' - r' ) = 14 + ( \w(G') - 6k' + (k'-1) ) = 14 + ( \w(G') - 5k' - 1 ) = 14 + ( \w(G') - 5(1-p) - 1 ) = \w(G') + 5p + 8$.
Therefore, $\w(G') + 5p + 8 \le \w(G) < \w(G') + 10p$, and so $8 < 5p$, implying that $p \ge 2$. However, $p = 1 - k' \le 1$, a contradiction.

Hence, $\deg_G(v_3) = 2$. Let $v_4$ be the neighbor of $v_3$ different from~$v$. Suppose that $\deg_G(v_4) = 3$. Let $G' = G - \{v,v_1,v_2,v_3\}$. Suppose that $G' \notin \cB_{\rdom}$, implying that $\w(G) = 19 + (\w(G') - 1) = \w(G') + 18$. Let $S'$ be a $\gamma_r$-set of $G'$. If $v_4 \in S'$, let $S = S' \cup \{v_1\}$, and if $v_4 \notin S'$, let $S = S' \cup \{v\}$. In both cases, $S$ is a RD-set of $G$, and so $\gamma_r(G) \le \gamma_r(G') + 1$. Hence, $10\gamma_r(G) \le 10(\gamma_r(G') + 1) \le \w(G') + 10 = (\w(G) - 18) + 10 < \w(G)$, a contradiction. Hence, $G' \in \cB_{\rdom}$, and so the graph $G$ is determined. If $G' = R_2$ and $v_4$ is an open twin of $G'$, then $\gamma_r(G) = 4$ and $\w(G) = 46$, and so $10\gamma_r(G) < \w(G)$, a contradiction. If $G' = R_6$ and $v_4$ belongs to the $4$-cycle in $R_6$, then $\gamma_r(G) = 5$ and $\w(G) = 57$, and so $10\gamma_r(G) < \w(G)$, a contradiction. Hence, if $G' = R_2$, then $v_4$ is not an open twin of $G'$, and if $G' = R_6$, then $v_4$ does not belong to the $4$-cycle in $R_6$. By Observation~\ref{obser-1}(e), $\gamma_r(G) \le |\{v\}| + \gamma_{r,\dom}(G';v_4) \le 1 + (\gamma_r(G') - 1) = \gamma_r(G')$, and so $10\gamma_r(G) \le 10\gamma_r(G') \le \w(G')$. However, $\w(G) \ge 19 + ( \w(G') - 6 ) = \w(G') + 13$, a contradiction.

Hence, $\deg_G(v_4) = 2$. Let $v_5$ be the neighbor of $v_4$ different from~$v_3$. By Claim~\ref{claim.3}, we have $\deg_G(v_5) = 3$. Let $Q = \{v,v_1,v_2,v_3,v_4\}$ and let $G' = G - Q$. We note that $G'$ is a connected special subcubic graph and $k' + r' = 1$. Applying Claim~\ref{claim.key2} with $H = G[Q]$ and $S_H = \{v_1,v_4\}$, we have $\w(G) < \w(G') + 10p$ where $p = \gamma_r(H) - k' = 2 - k'$. Hence, $\w(G) \ge \w_G(Q) + ( \w(G') - 6k' - r' ) = 24 + ( \w(G') - 6k' + (k'-1) ) = 24 + ( \w(G') - 5k' - 1 ) = 24 + ( \w(G') - 5(2-p) - 1 ) = \w(G') + 5p + 13$.
Therefore, $\w(G') + 5p + 13 \le \w(G) < \w(G') + 10p$, and so $13 < 5p$, implying that $p \ge 3$. However, $p = 2 - k' \le 2$, a contradiction.~\smallqed

\medskip
By Claim~\ref{claim.6}, there is no $3$-handle.

\begin{claim}
\label{claim.7}
There is no $4$-handle in $G$.
\end{claim}
\proof Suppose that $C \colon v v_1 v_2 v_3 v$ is a $4$-handle, where $\deg_G(v) = 3$. Let $v_4$ be the neighbor of $v$ not on $C$.

\begin{subclaim}
\label{claim.7.1}
$\deg_G(v_4) = 2$.
\end{subclaim}
\proof Suppose, to contrary, that $\deg_G(v_4) = 3$. Let $x$ and $y$ be the two neighbors of $v_4$ different from~$v_3$. Suppose that $x$ and $y$ are both large vertices. Let $Q = \{v,v_1,v_2,v_3,v_4\}$ and let $G' = G - Q$. We note that $G'$ has at most two components, and so $k' + r' \le 2$. Applying Claim~\ref{claim.key2} with $H = G[Q]$ and $S_H = \{v_2,v_4\}$, we have $\w(G) < \w(G') + 10p$ where $p = \gamma_r(H) - k' = 2 - k'$. On the other hand, $\w(G) \ge 23 + ( \w(G') - 6k' - r' ) \ge 23 + ( \w(G') - 6k' +(k'-2) ) = 23 + ( \w(G') - 5k' - 2 ) = 23 + ( \w(G') - 5(2-p) - 2 ) = \w(G') + 5p + 11$.
Therefore, $\w(G') + 5p + 11 \le \w(G) < \w(G') + 10p$, and so $11 < 5p$, implying that $p \ge 3$. However, $p = 2 - k' \le 2$, a contradiction.

Hence, at least one of $x$ and $y$ is a small vertex. Renaming vertices if necessary, we may assume that $\deg_G(x) = 2$. Note that $\deg_G(y) \in \{2,3\}$. Suppose $xy \in E(G)$. Since there is no $3$-handle, the vertex $y$ is large. Let $z$ be the neighbor of $y$ different from $x$ and $v_4$. Let $G'$ be obtained from $G$ by deleting $x,y$ and $v_4$, and adding the edge $vz$. The graph $G'$ is a connected special subcubic graph of order less than~$n$. Since no graph in $\cB_{\rdom}$ contains a $4$-handle, we note that $G' \notin \cB_{\rdom}$, implying that $\w(G) = \w(G') + 13$. Let $S'$ be a $\gamma_r$-set of $G'$. If $v \in S'$, let $S = S' \cup \{y\}$. If $v \notin S'$ and $z \in S'$, let $S = S' \cup \{v_4\}$. If $v \notin S'$ and $z \notin S'$, let $S = S' \cup \{x\}$. In all cases, $S$ is a RD-set of $G$, and so $\gamma_r(G) \le |S| = |S'| + 1 = \gamma_r(G') + 1$. Hence, $10\gamma_r(G) \le 10(\gamma_r(G') + 1) \le w(G') + 10 = (\w(G) - 13) + 10 < \w(G)$, a contradiction.

Hence, $xy \notin E(G)$. Let $Q = \{v,v_1,v_2,v_3,v_4\}$ and let $G'$ be obtained from $G - Q$ by adding the edge $xy$. The resulting graph $G'$ is a connected special subcubic graph of order less than~$n$. Suppose $G' \notin \cB_{\rdom}$. In this case, $\w(G) = 23 + \w(G')$. Let $S'$ be a $\gamma_r$-set of $G'$. If $x \in S'$ or $y \in S'$, let $S = S' \cup \{v_1,v_4\}$. If $x \notin S'$ and $y \notin S'$, let $S = S' \cup \{v,v_1\}$. In both cases $S$ is a RD-set of $G$, and so $\gamma_r(G) \le |S| = |S'| + 2 = \gamma_r(G') + 2$. Therefore, $10\gamma_r(G) \le 10(\gamma_r(G') + 2) \le \w(G') + 20 < \w(G)$, a contradiction.

Hence, $G' \in \cB_{\rdom}$. Let $G^* = G - V(C)$. We note that in this case, $G^*$ is obtained from $G'$ by subdividing the edge $xy$ of $G'$ where $v_4$ is the resulting vertex of degree~$2$ in $G^*$. By Observation~\ref{obser-2}, $\gamma_r(G^*) \le \gamma_r(G')$ and there exists a $\gamma_r$-set $S^*$ of $G^*$ that contains the vertex $v_4$. The set $S^* \cup \{v_2\}$ is a RD-set of $G$, and so $\gamma_r(G) \le 1 + |S^*| = 1 + \gamma_r(G^*) \le 1 + \gamma_r(G')$. Hence, $\w(G) < 10\gamma_r(G) \le 10(\gamma_r(G') + 1) \le \w(G') + 10$. Moreover noting that the degrees of the vertices in $G'$ are the same as their degrees in $G$, we have $\w(G) \ge 23 + ( \w(G') - 5) = \w(G') + 18$, a contradiction.~\smallqed

\medskip
By Claim~\ref{claim.7.1}, we have $\deg_G(v_4) = 2$. Let $v_5$ be the neighbor of $v_4$ different from $v$. Suppose that $\deg_G(v_5) = 3$. Let $Q = \{v,v_1,v_2,v_3,v_4\}$ and let $G' = G - Q$. The resulting graph $G'$ is a connected special subcubic graph of order less than~$n$. We note that $k' + r' = 1$. Applying Claim~\ref{claim.key2} with $H = G[Q]$ and $S_H = \{v_2,v_4\}$, we have $\w(G) < \w(G') + 10p$ where $p = \gamma_r(H) - k = 2 - k$, implying $p \le 2$. On the other hand using the same calculations as in the earlier proofs, we have $\w(G) \ge 24 + ( \w(G') - 6k' - r' ) \ge \w(G') + 5p + 13$. Therefore, $\w(G') + 5p + 13 \le \w(G) < \w(G') + 10p$, and so $13 < 5p$, that is, $p \ge 3$. However, $p = 2 - k' \le 2$, a contradiction.

Hence, $\deg_G(v_5) = 2$. Let $v_6$ be the neighbor of $v_5$ different from $v_4$. By Claim~\ref{claim.3}, $\deg_G(v_6) = 3$. Let $Q = \{v,v_1,v_2,v_3,v_4,v_5\}$ and let $G' = G - Q$. The resulting graph $G'$ is a connected special subcubic graph of order less than~$n$. Let $S'$ be a $\gamma_r$-set of $G'$. If $v_6 \in S'$, let $S = S' \cup \{v,v_1\}$. If $v_6 \notin S'$, let $S = S' \cup \{v_2,v_4\}$. In both cases $S$ is a RD-set of $G$, and so $\gamma_r(G) \le |S| = |S'| + 2 = \gamma_r(G') + 2$.  We note that $k' + r' = 1$. Applying Claim~\ref{claim.key1} with $H = G[Q]$ and $p = 2$, we have $\w(G) < \w(G') + 10p = \w(G') + 20$. However, $\w(G) \ge 29 + ( \w(G') - 6k' - r' ) \ge 29 + \w(G') - 6(k'+r') = 29 + \w(G') - 6 = \w(G') + 23$, a contradiction. This completes the proof of Claim~\ref{claim.7}.~\smallqed

\medskip
By Claim~\ref{claim.7}, there is no $4$-handle.

\begin{claim}
\label{claim.8}
There is no handle in $G$.
\end{claim}
\proof Suppose, to the contrary, that $G$ contains a handle. By our earlier observations, it must be a $5$-handle. Let $C \colon v v_1 v_2 v_3 v_4 v$ is a $5$-handle, where $\deg_G(v) = 3$. Let $v_5$ be the third neighbor of $v$ not on $C$.

\begin{subclaim}
\label{claim.8.1}
$\deg_G(v_5) = 2$.
\end{subclaim}
\proof Suppose, to contrary, that $\deg_G(v_5) = 3$. Let $G' = G - V(C)$. The resulting graph $G'$ is a connected special subcubic graph of order less than~$n$. Suppose that $G' \in \{R_1,R_4,R_5,R_{10}\}$. If $G' = R_1$, then $\gamma_r(G) = 4$ and $\w(G) = 48$. If $G' \in \{R_4,R_5\}$, then $\gamma_r(G) = 5$ and $\w(G) = 59$. If $G' = R_{10}$, then $\gamma_r(G) = 4$ and $\w(G) = 52$. In all cases, $10 \gamma_r(G) \le \w(G)$, a contradiction. Hence, $G' \not\in \{R_1,R_4,R_5,R_{10}\}$. Let $S'$ be a $\gamma_r$-set of $G'$. If $v_5 \in S'$, let $S = S' \cup \{v_2,v_3\}$. If $v_5 \notin S'$, let $S = S' \cup \{v_1,v_4\}$. In both cases $S$ is a RD-set of $G$, and so $\gamma_r(G) \le |S| = |S'| + 2 = \gamma_r(G') + 2$.  We note that $k' + r' = 1$. Applying Claim~\ref{claim.key1} with $H = C$ and $p = 2$, we have $\w(G) < \w(G') + 10p = \w(G') + 20$. Since $G' \notin \{R_1,R_4,R_5,R_{10}\}$, when reconstructing the graph $G$ the contribution of the weight of $G'$ to the weight of $G$ decreases by at most~$3k' + r'$. Thus, $\w(G) \ge 24 + ( \w(G') - 3k' - r' ) \ge 24 + \w(G') - 3(k' + r') = 24 + \w(G') - 3 = \w(G') + 21$, a contradiction.~\smallqed

\medskip
By Claim~\ref{claim.8.1}, we have $\deg_G(v_5) = 2$. Let $v_6$ be the neighbor of $v_5$ different from $v$.

\begin{subclaim}
\label{claim.8.2}
$\deg_G(v_6) = 2$.
\end{subclaim}
\proof Suppose, to contrary, that $\deg_G(v_6) = 3$. Let $x$ and $y$ be the two neighbors of $v_6$ different from~$v_5$. Suppose that $x$ and $y$ are both large vertices. Let $Q = \{v,v_1,v_2,v_3,v_4,v_5,v_6\}$ and let $G' = G - Q$. We note that $G'$ has at most two components, and so $k' + r' \le 2$. Applying Claim~\ref{claim.key2} with $H = G[Q]$ and $S_H = \{v_1,v_4,v_6\}$, we have $\w(G) < \w(G') + 10p$ where $p = \gamma_r(H) - k' = 3 - k'$. On the other hand, $\w(G) \ge 33 + ( \w(G') - 6k' - r' ) \ge 33 + ( \w(G') - 6k' +(k'-2) ) = 33 + ( \w(G') - 5k' - 2 ) = 33 + ( \w(G') - 5(3-p) - 2 ) = \w(G') + 5p + 16$.
Therefore, $\w(G') + 5p + 16 \le \w(G) < \w(G') + 10p$, and so $16 < 5p$, that is, $p \ge 4$. However, $p = 3 - k' \le 3$, a contradiction.

Hence at least one of $x$ and $y$ is a small vertex. Renaming vertices if necessary, we may assume that $\deg_G(x) = 2$. Note that $\deg_G(y) \in \{2,3\}$. Suppose $xy \in E(G)$. Since there is no $3$-handle, the vertex $y$ is large. Let $z$ be the neighbor of $y$ different from $x$ and $v_6$. Let $G'$ be obtained from $G$ by deleting $x,y$ and $v_6$, and adding the edge $v_5z$. The graph $G'$ is a connected special subcubic graph of order less than~$n$. Since no graph in $\cB_{\rdom}$ contains a $5$-handle, we note that $G' \notin \cB_{\rdom}$, implying that $\w(G) = \w(G') + 13$. Let $S'$ be a $\gamma_r$-set of $G'$. If $v_5 \in S'$, let $S = S' \cup \{y\}$. If $v_5 \notin S'$ and $z \in S'$, let $S = S' \cup \{v_6\}$. If $v_5 \notin S'$ and $z \notin S'$, let $S = S' \cup \{x\}$. In all cases, $S$ is a RD-set of $G$, and so $\gamma_r(G) \le |S| = |S'| + 1 = \gamma_r(G') + 1$. Hence, $10\gamma_r(G) \le 10(\gamma_r(G') + 1) \le w(G') + 10 = (\w(G) - 13) + 10 < \w(G)$, a contradiction.

Hence, $xy \notin E(G)$. Let $Q = \{v,v_1,v_2,v_3,v_4,v_5,v_6\}$ and let $G'$ be obtained from $G - Q$ by adding the edge $xy$. The resulting graph $G'$ is a connected special subcubic graph of order less than~$n$. Suppose $G' \notin \cB_{\rdom}$. In this case, $\w(G) = 33 + \w(G')$. Let $S'$ be a $\gamma_r$-set of $G'$. If $x \in S'$ or $y \in S'$, let $S = S' \cup \{v_1,v_4,v_6\}$. If $x \notin S'$ and $y \notin S'$, let $S = S' \cup \{v_2,v_3,v_5\}$. In both cases $S$ is a RD-set of $G$, and so $\gamma_r(G) \le |S| = |S'| + 3 = \gamma_r(G') + 3$. Therefore, $10\gamma_r(G) \le 10(\gamma_r(G') + 3) \le \w(G') + 30 < \w(G)$, a contradiction.
Hence, $G' \in \cB_{\rdom}$. Let $G^* = G - \{v,v_1,v_2,v_3,v_4,v_5\}$. We note that in this case, $G^*$ is obtained from $G'$ by subdividing the edge $xy$ of $G'$ where $v_6$ is the resulting vertex of degree~$2$ in $G^*$. By Observation~\ref{obser-2}, $\gamma_r(G^*) \le \gamma_r(G')$ and there exists a $\gamma_r$-set $S^*$ of $G^*$ that contains the vertex $v_6$. The set $S^* \cup \{v_1,v_4\}$ is a RD-set of $G$, and so $\gamma_r(G) \le 2 + |S^*| = 2 + \gamma_r(G^*) \le \gamma_r(G') + 2$. Hence, $\w(G) < 10\gamma_r(G) \le 10(\gamma_r(G') + 2) \le \w(G') + 20$. However noting that the degrees of the vertices in $G'$ are the same as their degrees in $G$, we have $\w(G) \ge 33 + ( \w(G') - 5) = \w(G') + 28$, a contradiction.~\smallqed

\medskip
By Claim~\ref{claim.8.2}, we have $\deg_G(v_6) = 2$. Let $v_7$ be the neighbor of $v_6$ different from $v_5$. By Claim~\ref{claim.3}, $\deg_G(v_7) = 3$. Let $Q = \{v,v_1,v_2,v_3,v_4,v_5,v_6\}$ and let $G' = G - Q$. The resulting graph $G'$ is a connected special subcubic graph of order less than~$n$. We note that $k'+r'=1$. Applying Claim~\ref{claim.key2} with $H = G[Q]$ and $S_H = \{v_1,v_4,v_6\}$, we have $\w(G) < \w(G') + 10p$ where $p = \gamma_r(H) - k' = 3 - k'$. On the other hand, $\w(G) \ge 33 + ( \w(G') - 6k' - r' ) = 33 + ( \w(G') - 6k' +(k'-1) ) = 33 + ( \w(G') - 5k' - 1 ) =  33 + ( \w(G') - 5(3-p) - 1 ) = \w(G') + 15p + 17$.
Therefore, $\w(G') + 15p + 17 \le \w(G) < \w(G') + 10p$, and so $17 < 5p$, that is, $p \ge 4$. However, $k' \ge 0$ and $p = 3 - k' \le 3$, a contradiction. This completes the proof of Claim~\ref{claim.8}.~\smallqed

\medskip
By Claim~\ref{claim.8}, there is no handle in $G$. In particular, the removal of a bridge cannot create a $C_5$-component. Recall that there is no $k$-linkage for any $k \ge 3$. Hence if $\delta(G) = 2$, then every vertex of degree~$2$ in $G$ belongs to a $k$-linkage for some $k \in \{1,2\}$.

\begin{claim}
\label{claim.9}
If $G$ contains a $2$-linkage, then the two large vertices on the linkage are not adjacent.
\end{claim}
\proof Suppose, to the contrary, that $G$ contains a $2$-linkage $P \colon v v_1 v_2 u$ where $u$ and $v$ are adjacent. We note that $u, v \in \cL$ and $v_1,v_2 \in \cS$.

\begin{subclaim}
\label{claim.9.1}
The vertices $u$ and $v$ have no common neighbor.
\end{subclaim}
\proof Suppose that $u$ and $v$ have a common neighbor $v_3$. Since $n \ge 6$, the vertex $v_3$ is large. Let $v_4$ be the neighbor of $v_3$ not on $P$. Suppose that $\deg_G(v_4) = 3$. Let $Q = \{v,v_1,v_2,v_3,u\}$ and let $G' = G - Q$. The graph $G'$ is a connected special subcubic graph of order less than~$n$. We note that $k' + r' = 1$. Applying Claim~\ref{claim.key2} with $H = G[Q]$ and $S_H = \{v_1,v_3\}$, we have $\w(G) < \w(G') + 10p$ where $p = \gamma_r(H) - k' = 2 - k'$. Since $G$ has no handle, we note that $G' \ne R_1$, implying that $\w(G) \ge 22 + ( \w(G') - 5k' - r' ) = 22 + ( \w(G') - 5k' + k'-1 ) = 22 + ( \w(G') - 4k' - 1 ) = 22 + ( \w(G') - 4(2-p) - 1 ) \ge \w(G') + 4p + 13$.
Therefore, $\w(G') + 4p + 13 \le \w(G) < \w(G') + 10p$, and so $13 < 6p$, that is, $p \ge 3$. However, $p = 2 - k' \le 2$, a contradiction.

Hence, $\deg_G(v_4) = 2$. Let $v_5$ be the neighbor of $v_4$ different from $v_3$. Suppose that $\deg_G(v_5) = 3$. Let $G' = G - \{v,v_1,v_2,v_3,v_4,u\}$. Let $S'$ be a $\gamma_r$-set of $G'$. If $v_5 \in S'$, let $S = S' \cup \{u,v\}$. If $v_5 \notin S'$, let $S = S' \cup \{v_1,v_3\}$. In both cases, $S$ is a RD-set of $G$, and so $\gamma_r(G) \le |S| = |S'| + 2 = \gamma_r(G') + 2$. Applying Claim~\ref{claim.key1} with $p = 2$, we have $\w(G) < \w(G') + 10p = \w(G') + 20$. Recall that $G$ has no handle, and so $G' \ne R_1$. Therefore, $\w(G) \ge 27 + (\w(G') - 1 - 4) = \w(G') + 22$, a contradiction.

Hence, $\deg_G(v_5) = 2$. Let $v_6$ be the neighbor of $v_5$ different from $v_4$. By Claim~\ref{claim.3}, $\deg_G(v_6) = 3$. Let $Q = \{v,v_1,v_2,v_3,v_4,v_5,u\}$ and let $G' = G - Q$. The resulting graph $G'$ is a connected special subcubic graph of order less than~$n$. We note that $k' + r' = 1$. Applying Claim~\ref{claim.key2} with $H = G[Q]$ and $S_H = \{u,v,v_5\}$, we have $\w(G) < \w(G') + 10p$ where $p = \gamma_r(H) - k' = 3 - k'$. On the other hand, noting that $G' \ne R_1$, we have $\w(G) \ge 32 + ( \w(G') - 5k' - r' ) \ge 32 + ( \w(G') - 5k' +k'-1 ) = 32 + ( \w(G') - 4k' - 1 ) = 32 + ( \w(G') - 4(3-p) - 1 ) = \w(G') + 4p + 19$.
Therefore, $\w(G') + 4p + 19 \le \w(G) < \w(G') + 10p$, and so $19 < 6p$, that is, $p \ge 4$. However, $k' \ge 0$ and $p = 3 - k' \le 3$, a contradiction.~\smallqed

\medskip
By Claim~\ref{claim.9.1}, the vertices $u$ and $v$ have no common neighbor. Let $v_3$ be the third neighbor of $v$ not on $P$. Since $u$ and $v$ have no common neighbor, $u$ and $v_3$ are not adjacent. Let $G'$ be obtained from $G - \{v,v_1,v_2\}$ by adding the edge $uv_3$. The resulting graph $G'$ is a connected special subcubic graph of order less than~$n$. Suppose that $G' \notin \cB_{\rdom}$. In this case, $\w(G) = 14 + (\w(G') - 1) = \w(G') + 13$. Let $S'$ be a $\gamma_r$-set of $G'$. If $u \in S'$, let $S = S' \cup \{v\}$. If $u \notin S'$ and $v_3 \in S'$, let $S = S' \cup \{v_2\}$. If $u \notin S'$ and $v_3 \notin S'$, let $S = S' \cup \{v_1\}$. In all cases $S$ is a RD-set of $G$, and so $\gamma_r(G) \le |S| = |S'| + 1 = \gamma_r(G') + 1$. Therefore, $\w(G) < 10\gamma_r(G) \le 10(\gamma_r(G') + 1) \le \w(G') + 10 < \w(G)$, a contradiction.

Hence, $G' \in \cB_{\rdom}$. If $G' = R_1$, then $G$ would contain a $4$-linkage, a contradiction. If $G \in \{R_4,R_5\}$, then $\gamma_r(G) = 4$ and $\w(G) = 49$, and so $10\gamma_r(G) \le \w(G)$, a contradiction. Hence,  $G' \notin \{R_1,R_4,R_5\}$. Let $G^* = G - \{v_1,v_2\}$. Thus, $G^*$ is obtained from $G'$ by subdividing the edge $uv_3$ of $G'$ where $v$ is the resulting vertex of degree~$2$ in $G^*$. By Observation~\ref{obser-2}, $\gamma_r(G^*) \le \gamma_r(G')$ and there exists a $\gamma_r$-set $S^*$ of $G^*$ that contains the vertex $v$ and does not contain $u$ or $v_3$. The set $S^* \cup \{v_1\}$ is a RD-set of $G$, and so $\gamma_r(G) \le 1 + |S^*| = 1 + \gamma_r(G^*) \le 1 + \gamma_r(G')$. Hence, $\w(G) < 10\gamma_r(G) \le 10(\gamma_r(G') + 1) \le \w(G') + 10$. We note that the degrees of the vertices in $G'$ are the same as their degrees in $G$, except for the vertex~$u$ which has degree~$3$ in $G$ and degree~$2$ in $G'$. As observed earlier, $G' \notin \{R_1,R_4,R_5\}$, implying that $\w(G) \ge 14 + (\w(G') - 1 - 3) = \w(G') + 10$, a contradiction. This completes the proof of Claim~\ref{claim.9}.~\smallqed

\begin{claim}
\label{claim.10}
If $G$ contains a $1$-linkage, then the two large vertices on the linkage are not adjacent.
\end{claim}
\proof Suppose, to the contrary, that $G$ contains a $1$-linkage $P \colon v v_1 u$ where $u$ and $v$ are adjacent. We note that $u, v \in \cL$ and $v_1 \in \cS$.

\begin{subclaim}
\label{claim.10.1}
The vertices $u$ and $v$ have no common neighbor.
\end{subclaim}
\proof Suppose that $u$ and $v$ have a common neighbor $v_2$, and so $G[\{v,v_1,v_2,u\}]$ is a diamond. Since $n \ge 6$, the vertex $v_2$ is large. Let $v_3$ be the third neighbor of $v_2$ not on $P$. Suppose that $\deg_G(v_3) = 3$. Let $Q = \{v,v_1,v_2,u\}$ and let $G' = G - Q$. The graph $G'$ is a connected special subcubic graph of order less than~$n$. We note that $k' + r' = 1$. Every $\gamma_r$-set of $G'$ can be extended to a RD-set of $G$ by adding to it the vertex~$v$, and so $\gamma_r(G) \le \gamma_r(G') + 1$. Thus, $\w(G) < 10\gamma_r(G) \le 10(\gamma_r(G') + 1) \le \w(G') + 10$. Since there is no handle in $G$, we note that $G' \ne R_1$, implying that $\w(G) \ge 17 + (\w(G') - 1 - 4) = \w(G') + 12$, a contradiction.

Hence, $\deg_G(v_3) = 2$. Let $v_4$ be the neighbor of $v_3$ different from $v_2$. Suppose that $\deg_G(v_4) = 3$. Let $Q = \{v,v_1,v_2,v_3,u\}$ and let $G' = G - Q$. The graph $G'$ is a connected special subcubic graph of order less than~$n$ different from $R_1$. We note that $k' + r' = 1$. Applying Claim~\ref{claim.key2} with $H = G[Q]$ and $S_H = \{v_1,v_3\}$, we have $\w(G) < \w(G') + 10p$ where $p = \gamma_r(H) - k' = 2 - k'$. On the other hand,
$\w(G) \ge 22 + ( \w(G') - 5k' - r' ) = \w(G') + 4p + 13$. Therefore, $\w(G') + 5p + 13 \le \w(G) < \w(G') + 10p$, and so $13 < 5p$, that is, $p \ge 3$. However, $p = 2 - k' \le 2$, a contradiction.
Hence, $\deg_G(v_4) = 2$. Let $v_5$ be the neighbor of $v_4$ different from $v_3$. By Claim~\ref{claim.3}, $\deg_G(v_5) = 3$. Let $Q = \{v,v_1,v_2,v_3,v_4,u\}$ and let $G' = G - Q$. The resulting graph $G'$ is a connected special subcubic graph of order less than~$n$. We note that $k' + r' = 1$. Applying Claim~\ref{claim.key2} with $H = G[Q]$ and $S_H = \{v,v_4\}$, we have $\w(G) < \w(G') + 10p$ where $p = \gamma_r(H) - k' = 2 - k'$. On the other hand noting that $G' \ne R_1$, we have $\w(G) \ge 27 + ( \w(G') - 5k' - r' ) \ge \w(G') + 4p + 18$.
Therefore, $\w(G') + 4p + 18 \le \w(G) < \w(G') + 10p$, and so $18 < 6p$, that is, $p \ge 4$. However, $k' \ge 0$ and $p = 2 - k' \le 2$, a contradiction.~\smallqed 

\medskip
By Claim~\ref{claim.10.1}, the vertices $u$ and $v$ have no common neighbor. Let $v_2$ and $u_2$ be the third neighbors of $v$ and $u$, respectively, not on $P$. Since $u$ and $v$ have no common neighbor, $u_1 \ne v_2$.

\begin{subclaim}
\label{claim.10.2}
The vertices $u_2$ and $v_2$ are not adjacent.
\end{subclaim}
\proof Suppose that $u_2$ and $v_2$ are adjacent. Since $n \ge 6$, at least one of $u_2$ and $v_2$ is large. Renaming vertices if necessary, assume that $u_2 \in \cL$. Suppose that $v_2 \in \cS$ and $N(v_2) = \{v,u_2\}$. Let $u_3$ be the neighbor of $u_2$ different from $u$ and $v_2$. Suppose that $u_3 \in \cL$. Let $Q = \{v,v_1,v_2,u,u_2\}$ and let $G' = G - Q$. The graph $G'$ is a connected special subcubic graph of order less than~$n$. We note that $k' + r' = 1$. Applying Claim~\ref{claim.key2} with $H = G[Q]$ and $S_H = \{v_1,u_2\}$, we have $\w(G) < \w(G') + 10p$ where $p = \gamma_r(H) - k' = 2 - k'$. On the other hand, $\w(G) \ge 22 + ( \w(G') - 5k' - r' ) = \w(G') + 4p + 13$. Therefore, $\w(G') + 4p + 13 \le \w(G) < \w(G') + 10p$, and so $13 < 6p$, that is, $p \ge 3$. However, $p = 2 - k' \le 2$, a contradiction. Hence, $u_3 \in \cS$. Let $u_4$ be the neighbor of $u_3$ different from $u_2$.

Suppose that $u_4 \in \cL$. Let $Q = \{v,v_1,v_2,u,u_2,u_3\}$ and let $G' = G - Q$. The graph $G'$ is a connected special subcubic graph of order less than~$n$. We note that $k' + r' = 1$. Applying Claim~\ref{claim.key2} with $H = G[Q]$ and $S_H = \{v,u_3\}$, we have $\w(G) < \w(G') + 10p$ where $p = \gamma_r(H) - k' = 2 - k'$. On the other hand, $\w(G) \ge 27 + ( \w(G') - 5k' - r' ) = \w(G') + 4p + 18$. Therefore, $\w(G') + 4p + 18 \le \w(G) < \w(G') + 10p$, and so $18 < 6p$, that is, $p \ge 4$. However, $p = 2 - k' \le 2$, a contradiction.

Hence, $u_4 \in \cS$. Let $u_5$ be the neighbor of $u_4$ different from $u_3$. By Claim~\ref{claim.3}, $\deg_G(u_3) = 3$. Let $Q = \{v,v_1,v_2,u,u_2,u_3,u_4\}$ and let $G' = G - Q$. The resulting graph $G'$ is a connected special subcubic graph of order less than~$n$. We note that $k' + r' = 1$. Let $S'$ be a $\gamma_r$-set of $G'$. If $u_5 \in S'$, let $S = S' \cup \{u_2,v_1\}$. If $u_5 \notin S'$, let $S = S' \cup \{v,u_3\}$. In both cases, $S$ is a RD-set of $G$, and so $\gamma_r(G) \le |S| = |S'| + 2 = \gamma_r(G') + 2$. Applying Claim~\ref{claim.key1} with $p = 2$, we have $\w(G) < \w(G') + 10p = \w(G') + 20$. However, $\w(G) \ge 32 + (\w(G') - 1 - 4) = \w(G') + 27$, a contradiction.

Hence, $v_2 \in \cL$. Recall that $u_2 \in \cL$. Let $Q = \{v,v_1,u\}$ and let $G' = G - Q$. We note that $k' + r' = 1$. Applying Claim~\ref{claim.key2} with $H = G[Q]$ and $S_H = \{v\}$, we have $\w(G) < \w(G') + 10p$ where $p = \gamma_r(H) - k' = 1 - k'$. Since there is no $3$-linkage in $G$, we note that $G' \ne R_1$, implying that $\w(G) \ge 13 + ( \w(G') - 5k' - r' ) = \w(G') + 4p + 7$.
Therefore, $\w(G') + 4p + 7 \le \w(G) < \w(G') + 10p$, and so $7 < 6p$, that is, $p \ge 2$. However, $p = 1 - k' \le 1$, a contradiction.~\smallqed

\medskip
By Claim~\ref{claim.10.2}, the vertices $u_2$ and $v_2$ are not adjacent. Let $G'$ be obtained from $G - \{v,v_1,u\}$ by adding the edge $u_2v_2$. The resulting graph $G'$ is a connected special subcubic graph of order less than~$n$. Suppose that $G' \notin \cB_{\rdom}$, implying that $\w(G) = 13 + \w(G')$. Let $S'$ be a $\gamma_r$-set of $G'$. If $v_2 \in S'$, let $S = S' \cup \{u\}$. If $v_2 \notin S'$ and $u_2 \in S'$, let $S = S' \cup \{v\}$. If $u_2 \notin S'$ and $v_2 \notin S'$, let $S = S' \cup \{v_1\}$. In all cases $S$ is a RD-set of $G$, and so $\gamma_r(G) \le |S| = |S'| + 1 = \gamma_r(G') + 1$. Therefore, $\w(G) < 10\gamma_r(G) \le 10(\gamma_r(G') + 1) \le \w(G') + 10 = \w(G) - 3 < \w(G)$, a contradiction.
Hence, $G' \in \cB_{\rdom}$. Let $G^* = G - v_1$. We note that in this case, $G^*$ is obtained from $G'$ by subdividing the edge $u_2v_3$ of $G'$ twice where $v_2 v u u_2$ is the resulting path in $G^*$. By Observation~\ref{obser-3}, $\gamma_r(G^*) \le \gamma_r(G')$ and there exists a $\gamma_r$-set $S^*$ of $G^*$ that contains the vertex $v$ and does not contain $u$. The set $S^*$ is a RD-set of $G$, and so $\gamma_r(G) \le |S^*| = \gamma_r(G^*) \le \gamma_r(G')$. Hence, $\w(G) < 10\gamma_r(G) \le 10\gamma_r(G') \le \w(G')$.
However, $\w(G) \ge 13 + (\w(G') - 5) = \w(G') + 9$, a contradiction. This completes the proof of Claim~\ref{claim.10}.~\smallqed

\medskip
Recall that $G$ has no handle. By Claim~\ref{claim.10}, no small vertex belongs to a triangle.  We state this formally.

\begin{claim}
\label{claim.11a}
No small vertex belongs to a triangle.
\end{claim}

\begin{claim}
\label{claim.11}
Two large vertices cannot be the ends of two common $2$-linkages.
\end{claim}
\proof Suppose, to the contrary, that there are two large vertices $u$ and $v$ that belong to two common $2$-linkages $u v_1 v_2 v$ and $v v_3 v_4 u$ in $G$. Thus, $C \colon u v_1 v_2 v v_3 v_4 u$ is a $6$-cycle in $G$, where $u,v \in \cL$ and $v_1,v_2,v_3,v_4 \in \cS$.

\begin{subclaim}
\label{claim.11.1}
The vertices $u$ and $v$ have no common neighbor.
\end{subclaim}
\proof Suppose that $u$ and $v$ have a common neighbor $v_5$.
If $v_5 \in \cS$, then the graph $G$ is determined and $\gamma_r(G) = 3$ and $\w(G) = 33$, a contradiction. Hence, $v_5 \in \cL$. Let $v_6$ be the neighbor of $v_5$ different from $u$ and $v$. Suppose that $v_6 \in \cL$. Let $Q = \{u,v,v_1,v_2,v_3,v_4,v_5\}$ and let $G' = G - Q$. The graph $G'$ is a connected special subcubic graph of order less than~$n$. We note that $k' + r' = 1$. Applying Claim~\ref{claim.key2} with $H = G[Q]$ and $S_H = \{v_1,v_3,v_5\}$, we have $\w(G) < \w(G') + 10p$ where $p = \gamma_r(H) - k' = 3 - k'$. On the other hand, $\w(G) \ge 32 + ( \w(G') - 5k' - r' ) = \w(G') + 4p + 19$.
Therefore, $\w(G') + 4p + 19 \le \w(G) < \w(G') + 10p$, and so $19 < 6p$, that is, $p \ge 4$. However,
$p = 3 - k' \le 3$, a contradiction.

Hence, $v_6 \in \cS$. Let $v_7$ be the neighbor of $v_6$ different from $v_5$. Suppose $v_7 \in \cL$. Let $Q = \{u,v,v_1,v_2,v_3,v_4,v_5,v_6\}$ and let $G' = G - Q$. The graph $G'$ is a connected special subcubic graph of order less than~$n$. We note that $k' + r' = 1$. Applying Claim~\ref{claim.key2} with $H = G[Q]$ and $S_H = \{v_1,v_3,v_6\}$, we have $\w(G) < \w(G') + 10p$ where $p = \gamma_r(H) - k' = 3 - k'$. On the other hand, $\w(G) \ge 37 + ( \w(G') - 5k' - r' ) = \w(G') + 4p + 24$.
Therefore, $\w(G') + 4p + 24 \le \w(G) < \w(G') + 10p$, and so $24 < 6p$, that is, $p \ge 5$. However,
$p = 3 - k' \le 3$, a contradiction.

Hence, $v_7 \in \cS$. Let $v_8$ be the neighbor of $v_7$ different from $v_6$. By Claim~\ref{claim.3}, $\deg_G(v_8) = 3$. Let $Q = \{u,v,v_1,v_2,v_3,v_4,v_5,v_6,v_7\}$ and let $G' = G - Q$. The graph $G'$ is a connected special subcubic graph of order less than~$n$. Let $S'$ be a $\gamma_r$-set of $G'$. If $v_8 \in S'$, let $S = S' \cup \{v_1,v_3,v_5\}$. If $v_5 \notin S'$, let $S = S' \cup \{v_1,v_3,v_6\}$. In both cases, $S$ is a RD-set of $G$, and so $\gamma_r(G) \le |S| = |S'| + 3 = \gamma_r(G') + 3$. Applying Claim~\ref{claim.key1} with $H = G[Q]$ and $p = 3$, we have $\w(G) < \w(G') + 10p = \w(G') + 30$. However, $\w(G) \ge 42 + (\w(G') - 1 - 4) = \w(G') + 37$, a contradiction.~\smallqed


\medskip
By Claim~\ref{claim.11.1}, the vertices $u$ and $v$ have no common neighbor. Let $x$ be the neighbor of $u$ different from $v_1$ and $v_4$, and let $y$ be the neighbor of $v$ different from $v_2$ and $v_3$.
Suppose that $x$ and $y$ are adjacent. If both $x$ and $y$ have degree~$2$, then the graph $G$ is determined and $\gamma_r(G) = 2$ and $\w(G) = 38$, a contradiction. Hence at least one of $x$ and $y$ are large. Renaming vertices if necessary, assume that $y$ is large. An analogous proof as before shows that $x \in \cL$. Let $Q = \{u,v,v_1,v_2,v_3,v_4\}$ and consider the graph $G' = G - \{u,v,v_1,v_2,v_3,v_4\}$. Every $\gamma_r$-set of $G'$ can be extended to a RD-set of $G$ by adding to it the set $\{u,v\}$, and so $\gamma_r(G) \le \gamma_r(G') + 2$. Applying Claim~\ref{claim.key1} with $H = G[Q]$ and $p = 2$, we have $\w(G) < \w(G') + 10p = \w(G') + 20$. However, $\w(G) \ge 28 + (\w(G') - 1 - 4) = \w(G') + 23$, a contradiction.

Hence, the vertices $x$ and $y$ are not adjacent. Let $Q = \{u,v,v_1,v_2,v_3,v_4\}$, and let $G'$ be obtained from $G' = G - Q$ by adding the edge $xy$. The resulting graph $G'$ is a connected special subcubic graph of order less than~$n$. Suppose that $G' \notin \cB_{\rdom}$, implying that $\w(G) = 28 + \w(G')$. Every $\gamma_r$-set of $G'$ can be extended to a RD-set of $G$ by adding to it $u$ and $v$ or $v_1$ and $v_3$, implying that $\gamma_r(G) \le \gamma_r(G') + 2$. Therefore, $\w(G) < 10\gamma_r(G) \le 10(\gamma_r(G') + 2) \le \w(G') + 20 = \w(G) - 8 < \w(G)$, a contradiction.
Hence, $G' \in \cB_{\rdom}$. Let $G^* = G - \{v_3,v_4\}$. Thus, $G^*$ is obtained from $G'$ by subdividing the edge $xy$ of $G'$ four times where $x u v_1 v_2 v y$ is the resulting path in $G^*$. By Observation~\ref{obser-5}, there exists a RD-set $S^*$ of $G^*$ such that $|S^*| \le \gamma_{r}(G') + 1$ and $S^* \cap \{u,v_1,v_2,v\} = \{u,v\}$. The set $S^*$ is a RD-set of $G$, and so $\gamma_r(G) \le |S^*| = \gamma_{r}(G') + 1$. Hence, $\w(G) < 10\gamma_r(G) \le 10(\gamma_r(G') + 1) \le \w(G') + 10$. Since $G$ has no $3$-linkage, we note that $G' \ne R_1$, implying that $\w(G) \ge 28 + (\w(G') - 4) = \w(G') + 24$, a contradiction. This completes the proof of Claim~\ref{claim.11}.~\smallqed

\begin{claim}
\label{claim.12}
Two large vertices cannot be the ends of a common $1$-linkage and a common $2$-linkage.
\end{claim}
\proof Suppose, to the contrary, that there are two large vertices $u$ and $v$ such that $u v_1 v_2 v$ is a $2$-linkage and $u v_3 v$ is a $1$-linkage in $G$. Thus, $C \colon u v_1 v_2 v v_3 u$ is a $5$-cycle in $G$, where $u,v \in \cL$ and $v_1,v_2,v_3 \in \cS$.

\begin{subclaim}
\label{claim.12.1}
The vertex $v_3$ is the only common neighbor of $u$ and $v$.
\end{subclaim}
\proof Suppose that $u$ and $v$ have two common neighbors. Let $v_4$ be the common neighbor of~$u$ and~$v$ different from $v_3$.
If $v_4 \in \cS$, then $G = R_2$, a contradiction.
Hence, $v_4 \in \cL$. Let $v_5$ be the neighbor of $v_4$ different from $u$ and~$v$.

Suppose that $v_5 \in \cS$. Let $v_6$ be the neighbor of $v_5$ different from $v_4$. If $v_6 \in \cL$, then let $Q = \{u,v,v_1,v_2,v_3,v_4,v_5\}$ and let $G' = G - Q$. Applying Claim~\ref{claim.key2} with $H = G[Q]$ and $S_H = \{v_1,v_3,v_5\}$ we obtain a contradiction. Hence, $v_6 \in \cS$. Let $v_7$ be the neighbor of $v_6$ different from $v_5$. By Claim~\ref{claim.3}, we have $\deg_G(v_7) = 3$. Let $Q = \{u,v,v_1,v_2,v_3,v_4,v_5,v_6\}$ and let $G' = G - Q$. In this case, $\gamma_r(G) \le \gamma_r(G') + 3$, and applying Claim~\ref{claim.key1} with $H = G[Q]$ and $p = 3$ we obtain a contradiction.

Hence, $v_5 \in \cL$. Let $x$ and $y$ be the two neighbors of $v_5$ different from $v_4$. Suppose that $x$ and $y$ are not adjacent. Let $Q = \{u,v,v_1,v_2,v_3,v_4,v_5\}$ and let $G'$ be obtained from $G - Q$ by adding the edge $xy$. If $G' \notin \cB_{\rdom}$, then $\gamma_r(G) \le \gamma_r(G') + 3$, and applying Claim~\ref{claim.key1} with $H = G[Q]$ and $p = 3$ we obtain a contradiction. Hence, $G' \in \cB_{\rdom}$. In this case, we let $Q^* = \{u,v,v_1,v_2,v_3,v_4\}$ and let $G^* = G - Q^*$. Thus, $G^*$ is obtained from $G'$ by subdividing the edge $xy$ of $G'$ where $x v_5 y$ is the resulting path in $G^*$. Applying Observation~\ref{obser-2}, $\gamma_r(G) \le \gamma_{r}(G') + 2$, implying that $\w(G) < \w(G') + 20$. However since $G$ contains no $6$-handle, $G' \ne R_1$, and so $\w(G) \ge 31 + (\w(G') - 4) = \w(G') + 27$, a contradiction.

Hence, $xy \in E(G)$. Since there is no $3$-handle in $G$, at least one of $x$ and $y$ is a large vertex. Hence by Claim~\ref{claim.11a}, $x \in \cL$ and $y \in \cL$. Let $w = v_5$, and so $G[\{w,x,y\}]$ is a triangle. Let $x_1$ and $y_1$ be the neighbors of $x$ and $y$, respectively, different from~$w$.

We show next that $x_1 \ne y_1$. Suppose that $x_1 = y_1$. Since no vertex of degree~$2$ belongs to a triangle, $x_1 \in \cL$. Let $x_2$ be the neighbor of $x_1$ different from $x$ and $y$. If $x_2 \in \cL$, then we let $Q = \{u,v,v_1,v_2,v_3,v_4,w,x,y,x_1\}$ and $G' = G - Q$, and applying Claim~\ref{claim.key2} with $H = G[Q]$ and $S_H = \{v,v_2,v_4,x_1\}$ we obtain a contradiction. Hence, $x_2 \in \cS$. Let $x_3$ be the neighbor of $x_2$ different from $x_1$. If $x_3 \in \cL$, then we let $Q = \{u,v,v_1,v_2,v_3,v_4,w,x,y,x_1,x_2\}$ and $G' = G - Q$. In this case, $\gamma_r(G) \le \gamma_r(G') + 4$, and applying Claim~\ref{claim.key1} with $H = G[Q]$ and $p = 4$, we obtain a contradiction. Hence, $x_3 \in \cS$. Let $x_4$ be the neighbor of $x_3$ different from $x_2$. By Claim~\ref{claim.3}, $\deg_G(x_4) = 3$. Thus, $G$ contains the subgraph illustrated in Figure~\ref{rdom:fig-6d}. We now let $Q = \{u,v,v_1,v_2,v_3,v_4,w,x,y,x_1,x_2,x_3\}$ and $G' = G - Q$. In this case, $\gamma_r(G) \le \gamma_r(G') + 4$, and applying Claim~\ref{claim.key1} with $H = G[Q]$ and $p = 4$, we obtain a contradiction.

\begin{figure}[htb]
\begin{center}
\begin{tikzpicture}[scale=.8,style=thick,x=0.8cm,y=0.8cm]
\def\vr{2.5pt} 
\path (0,0.5) coordinate (u1);
\path (0,1.5) coordinate (u2);
\path (1,1) coordinate (u3);
\path (2,0) coordinate (u4);
\path (2,2) coordinate (u5);
\path (3,1) coordinate (u6);
\path (4,1) coordinate (u7);
\path (5,0) coordinate (u8);
\path (5,2) coordinate (u9);
\path (6,1) coordinate (u10);
\path (6.1,1) coordinate (u10p);
\path (7,1) coordinate (u11);
\path (8,1) coordinate (u12);
\path (9,1) coordinate (u13);
\path (8.9,1) coordinate (u13p);
\path (10,0.25) coordinate (u14);
\path (10,1.75) coordinate (u15);
\draw (u4)--(u1)--(u2)--(u5)--(u3)--(u4)--(u1);
\draw (u4)--(u6)--(u5);
\draw (u6)--(u7)--(u8)--(u9)--(u7);
\draw (u8)--(u10)--(u9);
\draw (u9)--(u10)--(u11)--(u12)--(u13);
\draw (u14)--(u13)--(u15);
\draw (u1) [fill=white] circle (\vr);
\draw (u2) [fill=white] circle (\vr);
\draw (u3) [fill=white] circle (\vr);
\draw (u4) [fill=white] circle (\vr);
\draw (u5) [fill=white] circle (\vr);
\draw (u6) [fill=white] circle (\vr);
\draw (u7) [fill=white] circle (\vr);
\draw (u8) [fill=white] circle (\vr);
\draw (u9) [fill=white] circle (\vr);
\draw (u10) [fill=white] circle (\vr);
\draw (u11) [fill=white] circle (\vr);
\draw (u12) [fill=white] circle (\vr);
\draw (u13) [fill=white] circle (\vr);
\draw (u14) [fill=white] circle (\vr);
\draw (u15) [fill=white] circle (\vr);
\draw[anchor = east] (u1) node {{\small $v_1$}};
\draw[anchor = east] (u2) node {{\small $v_2$}};
\draw[anchor = east] (u3) node {{\small $v_3$}};
\draw[anchor = north] (u4) node {{\small $u$}};
\draw[anchor = south] (u5) node {{\small $v$}};
\draw[anchor = south] (u6) node {{\small $v_4$}};
\draw[anchor = south] (u7) node {{\small $w$}};
\draw[anchor = north] (u8) node {{\small $x$}};
\draw[anchor = south] (u9) node {{\small $y$}};
\draw[anchor = south] (u10p) node {{\small $x_1$}};
\draw[anchor = south] (u11) node {{\small $x_2$}};
\draw[anchor = south] (u12) node {{\small $x_3$}};
\draw[anchor = south] (u13p) node {{\small $x_4$}};
\end{tikzpicture}
\end{center}
\begin{center}
\vskip -0.75 cm
\caption{A subgraph in the proof of Claim~\ref{claim.12.1}}
\label{rdom:fig-6d}
\end{center}
\vskip -0.75 cm
\end{figure}

Hence, $x_1 \ne y_1$, and so $G$ contains the subgraph illustrated in Figure~\ref{rdom:fig-6e}.
Suppose that $x_1 \in \cL$ and $y_1 \in \cL$. Let $Q = \{u,v,v_1,v_2,v_3,v_4,w,x,y\}$ and let $G' = G - Q$. Let $G_x$ and $G_y$ be the components of $G'$. Possibly, $G_x = G_y$, in which case $G'$ is connected. By our earlier observations, neither $G_x$ nor $G_y$ is an $R_1$-component. Applying Claim~\ref{claim.key1} with $H = G[Q]$ and $p = 3$ we have $\w(G) < \w(G') + 30$. If at most one component of $G'$ belongs to $\cB_{\rdom}$, then $\w(G) \ge \w(G') + 33$, a contradiction. Hence, $G_x \ne G_y$ and both $G_x$ and $G_y$ belong to $\cB_{\rdom}$. If $G_x \notin \{R_4,R_5\}$, then $\w(G) \ge \w(G') + 30$, a contradiction. Hence, $G_x \in \{R_4,R_5\}$. Analogously, $G_y \in \{R_4,R_5\}$. Let $G_w$ be the component of $G - v_4w$ that contains~$v_4$, and so $G_w = R_2$. We now take a NeRD-set of type-$2$ in $G_x$, and a NeRD-set of type-$1$ in each of $G_y$ and $G_w$, and extend these sets to a RD-set of $G$ by adding to them the vertices $w$ and~$y$. By Observation~\ref{obser-1}, $\gamma_r(G) \le 2 + \gamma_{r,\ndom}(G_w;w) + \gamma_{r,\dom}(G_x;x) + \gamma_{r,\ndom}(G_y;y) \le 2 + (\gamma_r(G_w) - 1) + (\gamma_r(G_x) - 1) + (\gamma_r(G_y) - 1) = \gamma_r(G_w) + \gamma_r(G_x) + \gamma_r(G_y) - 1$. Thus, $\w(G) < 10\gamma_r(G) \le \w(G_w) + \w(G_x) + \w(G_y) - 10$. However, $\w(G) \ge 12 + (\w(G_w) - 3) + (\w(G_x) - 5) + (\w(G_y) - 5) = \w(G_w) + \w(G_x) + \w(G_y) - 1$, a contradiction.

\begin{figure}[htb]
\begin{center}
\begin{tikzpicture}[scale=.8,style=thick,x=0.8cm,y=0.8cm]
\def\vr{2.5pt} 
\path (0,0.5) coordinate (u1);
\path (0,1.5) coordinate (u2);
\path (1,1) coordinate (u3);
\path (2,0) coordinate (u4);
\path (2,2) coordinate (u5);
\path (3,1) coordinate (u6);
\path (4.5,1) coordinate (u7);
\path (5.5,0) coordinate (u8);
\path (5.5,2) coordinate (u9);
\path (7,0) coordinate (u10);
\path (7,2) coordinate (u11);
\draw (u4)--(u1)--(u2)--(u5)--(u3)--(u4)--(u1);
\draw (u4)--(u6)--(u5);
\draw (u6)--(u7)--(u8)--(u9)--(u7);
\draw (u8)--(u10);
\draw (u9)--(u11);
\draw (u1) [fill=white] circle (\vr);
\draw (u2) [fill=white] circle (\vr);
\draw (u3) [fill=white] circle (\vr);
\draw (u4) [fill=white] circle (\vr);
\draw (u5) [fill=white] circle (\vr);
\draw (u6) [fill=white] circle (\vr);
\draw (u7) [fill=white] circle (\vr);
\draw (u8) [fill=white] circle (\vr);
\draw (u9) [fill=white] circle (\vr);
\draw (u10) [fill=white] circle (\vr);
\draw (u11) [fill=white] circle (\vr);
\draw[anchor = east] (u1) node {{\small $v_1$}};
\draw[anchor = east] (u2) node {{\small $v_2$}};
\draw[anchor = east] (u3) node {{\small $v_3$}};
\draw[anchor = north] (u4) node {{\small $u$}};
\draw[anchor = south] (u5) node {{\small $v$}};
\draw[anchor = south] (u6) node {{\small $v_4$}};
\draw[anchor = south] (u7) node {{\small $w$}};
\draw[anchor = north] (u8) node {{\small $x$}};
\draw[anchor = south] (u9) node {{\small $y$}};
\draw[anchor = north] (u10) node {{\small $x_1$}};
\draw[anchor = south] (u11) node {{\small $y_1$}};
\end{tikzpicture}
\end{center}
\begin{center}
\vskip -0.75 cm
\caption{A subgraph in the proof of Claim~\ref{claim.12.1}}
\label{rdom:fig-6e}
\end{center}
\vskip -0.75 cm
\end{figure}

Hence, $x_1 \in \cS$ or $y_1 \in \cS$. Renaming vertices if necessary, we may assume that $x_1 \in \cS$.
Suppose that $x_1y_1 \in E(G)$. By Claim~\ref{claim.9}, $y_1 \in \cL$. Let $y_2$ be the neighbor of $y_1$ different from $x_1$ and $y$. Thus, $G$ contains the subgraph illustrated in Figure~\ref{rdom:fig-7}. If $y_2 \in \cL$, then we let $Q = \{u,v,v_1,v_2,v_3,v_4,w,x,x_1,y,y_1\}$ and $G' = G - Q$, and applying Claim~\ref{claim.key2} with $H = G[Q]$ and $S_H = \{v_1,v_3,w,y_1\}$ we obtain a contradiction. Hence, $y_2 \in \cS$. Let $y_3$ be the neighbor of $y_2$ different from $y_1$. If $y_3 \in \cL$, then we let $Q = \{u,v,v_1,v_2,v_3,v_4,w,x,x_1,y,y_1,y_2\}$ and $G' = G - Q$, and applying Claim~\ref{claim.key2} with $H = G[Q]$ and $S_H = \{v_1,v_3,w,x_1,y_2\}$ we obtain a contradiction. Hence, $y_3 \in \cS$. Let $y_4$ be the neighbor of $y_3$ different from $y_2$. By Claim~\ref{claim.3}, $y_4 \in \cL$. We now let $Q = \{u,v,v_1,v_2,v_3,v_4,w,x,x_1,y,y_1,y_2,y_3\}$ and let $G' = G - Q$. Applying Claim~\ref{claim.key2} with $H = G[Q]$ and $S_H = \{v_1,v_3,w,x_1,y_3\}$, we obtain a contradiction.

\begin{figure}[htb]
\begin{center}
\begin{tikzpicture}[scale=.8,style=thick,x=0.8cm,y=0.8cm]
\def\vr{2.5pt} 
\path (0,0.65) coordinate (u1);
\path (0,1.35) coordinate (u2);
\path (1,1) coordinate (u3);
\path (1.5,0) coordinate (u4);
\path (1.5,2) coordinate (u5);
\path (2,1) coordinate (u6);
\path (2.2,1) coordinate (u6p);
\path (3.5,1) coordinate (u7);
\path (3.4,1) coordinate (u7p);
\path (4.5,0) coordinate (u8);
\path (4.5,2) coordinate (u9);
\path (5.5,0) coordinate (u10);
\path (5.5,2) coordinate (u11);
\path (6.5,2) coordinate (u12);
\draw (u4)--(u1)--(u2)--(u5)--(u3)--(u4)--(u6)--(u5);
\draw (u6)--(u7)--(u8)--(u9)--(u11);
\draw (u7)--(u9);
\draw (u8)--(u10)--(u11)--(u12);
\draw (u1) [fill=white] circle (\vr);
\draw (u2) [fill=white] circle (\vr);
\draw (u3) [fill=white] circle (\vr);
\draw (u4) [fill=white] circle (\vr);
\draw (u5) [fill=white] circle (\vr);
\draw (u6) [fill=white] circle (\vr);
\draw (u7) [fill=white] circle (\vr);
\draw (u8) [fill=white] circle (\vr);
\draw (u9) [fill=white] circle (\vr);
\draw (u10) [fill=white] circle (\vr);
\draw (u11) [fill=white] circle (\vr);
\draw (u12) [fill=white] circle (\vr);
\draw[anchor = east] (u1) node {{\small $v_1$}};
\draw[anchor = east] (u2) node {{\small $v_2$}};
\draw[anchor = east] (u3) node {{\small $v_3$}};
\draw[anchor = north] (u4) node {{\small $u$}};
\draw[anchor = south] (u5) node {{\small $v$}};
\draw[anchor = south] (u6p) node {{\small $v_4$}};
\draw[anchor = south] (u7p) node {{\small $w$}};
\draw[anchor = north] (u8) node {{\small $x$}};
\draw[anchor = south] (u9) node {{\small $y$}};
\draw[anchor = north] (u10) node {{\small $x_1$}};
\draw[anchor = south] (u11) node {{\small $y_1$}};
\draw[anchor = south] (u12) node {{\small $y_2$}};
\end{tikzpicture}
\end{center}
\begin{center}
\vskip -0.75 cm
\caption{A subgraph in the proof of Claim~\ref{claim.12.1}}
\label{rdom:fig-7}
\end{center}
\vskip -0.5 cm
\end{figure}

Hence, $x_1y_1 \notin E(G)$. Let $z$ be the neighbor of $x_1$ different from~$x$. Suppose that $y_1z \notin E(G)$. In this case, let $Q = \{u,v,v_1,v_2,v_3,v_4,w,x,x_1,y\}$ and let $G'$ be obtained from $G - Q$ by adding the edge $y_1z$. If $G' \notin \cB_{\rdom}$, then $\w(G) = 44 + \w(G')$. However, $\gamma_r(G) \le \gamma_r(G') + 4$, and so $\w(G) < \w(G') + 40$, a contradiction. Hence, $G' \in \cB_{\rdom}$. Let $G^* = G - \{u,v,v_1,v_2,v_3,v_4,w\}$, and so $G^*$ is obtained from $G'$ by subdividing the edge $y_1z$ of $G'$ three times where $z x_1 x y y_1$ is the resulting path in $G^*$. A NeRD-set of type-$1$ in $G^*$ with respect to the vertex~$y$ can be extended to a RD-set by adding to it the set $\{v_1,v_3,w\}$, implying by Observation~\ref{obser-4}(a) that $\gamma_r(G) \le \gamma_{r,\ndom}(G^*;y) + 3 \le \gamma_r(G') + 3$, and so $\w(G) < \w(G') + 30$. Since $G' \ne R_1$, we have $\w(G) \ge 44 + (\w(G') - 4) = \w(G) + 40$, a contradiction.

Hence, $y_1z \in E(G)$. Thus $G$ contains the subgraph illustrated in Figure~\ref{rdom:fig-8}, where $x_1 \in \cS$. Since there is no $3$-linkage, $y_1 \in \cL$ or $z \in \cL$. If $y_1 \in \cL$ and $z \in \cL$, then we let $Q = \{u,v,v_1,v_2,v_3,v_4,w,x,x_1,y\}$ and $G' = G - Q$, and applying Claim~\ref{claim.key2} with $H = G[Q]$ and $S_H = \{v_1,v_3,w,x_1\}$ we obtain a contradiction. Hence, either $y_1 \in \cS$ and $z \in \cL$ or $y_1 \in \cL$ and $z \in \cS$.

\begin{figure}[htb]
\begin{center}
\begin{tikzpicture}[scale=.8,style=thick,x=0.8cm,y=0.8cm]
\def\vr{2.5pt} 
\path (0,0.65) coordinate (u1);
\path (0,1.35) coordinate (u2);
\path (1,1) coordinate (u3);
\path (1.5,0) coordinate (u4);
\path (1.5,2) coordinate (u5);
\path (2,1) coordinate (u6);
\path (2.2,1) coordinate (u6p);
\path (3.5,1) coordinate (u7);
\path (3.4,1) coordinate (u7p);
\path (4.5,0) coordinate (u8);
\path (4.5,2) coordinate (u9);
\path (5.5,0) coordinate (u10);
\path (5.5,2) coordinate (u11);
\path (6.0,1) coordinate (u12);
\draw (u4)--(u1)--(u2)--(u5)--(u3)--(u4)--(u6)--(u5);
\draw (u6)--(u7)--(u8)--(u9)--(u11);
\draw (u7)--(u9);
\draw (u8)--(u10);
\draw (u11)--(u12)--(u10);
\draw (u1) [fill=white] circle (\vr);
\draw (u2) [fill=white] circle (\vr);
\draw (u3) [fill=white] circle (\vr);
\draw (u4) [fill=white] circle (\vr);
\draw (u5) [fill=white] circle (\vr);
\draw (u6) [fill=white] circle (\vr);
\draw (u7) [fill=white] circle (\vr);
\draw (u8) [fill=white] circle (\vr);
\draw (u9) [fill=white] circle (\vr);
\draw (u10) [fill=white] circle (\vr);
\draw (u11) [fill=white] circle (\vr);
\draw (u12) [fill=white] circle (\vr);
\draw[anchor = east] (u1) node {{\small $v_1$}};
\draw[anchor = east] (u2) node {{\small $v_2$}};
\draw[anchor = east] (u3) node {{\small $v_3$}};
\draw[anchor = north] (u4) node {{\small $u$}};
\draw[anchor = south] (u5) node {{\small $v$}};
\draw[anchor = south] (u6p) node {{\small $v_4$}};
\draw[anchor = south] (u7p) node {{\small $w$}};
\draw[anchor = north] (u8) node {{\small $x$}};
\draw[anchor = south] (u9) node {{\small $y$}};
\draw[anchor = north] (u10) node {{\small $x_1$}};
\draw[anchor = south] (u11) node {{\small $y_1$}};
\draw[anchor = west] (u12) node {{\small $z$}};
\end{tikzpicture}
\end{center}
\begin{center}
\vskip -0.75 cm
\caption{A subgraph in the proof of Claim~\ref{claim.12.1}}
\label{rdom:fig-8}
\end{center}
\vskip -0.75 cm
\end{figure}

Suppose that $y_1 \in \cS$ and $z \in \cL$. Let $z_1$ be the neighbor of $z$ different from $x_1$ and $y_1$. If $G'$ is obtained from $G - \{w,x,x_1,y,y_1,z\}$ by adding the edge $v_4z_1$, then $\gamma_r(G) \le \gamma_r(G') + 2$, and so $\w(G) < \w(G') + 20$. Since the graph $G'$ contains a bridge, we note that $G' \notin \cB_{\rdom}$, implying that $\w(G) = \w(G') + 26$, a contradiction. Hence, $y_1 \in \cL$ and $z \in \cS$. Let $y_2$ be the neighbor of $y_1$ different from $y$ and $z$. If $G'$ is obtained from $G - \{w,x,x_1,y,y_1,z\}$ by adding the edge $v_4y_2$, then $\gamma_r(G) \le \gamma_r(G') + 2$, and so $\w(G) < \w(G') + 20$. Since $G' \notin \cB_{\rdom}$, we have $\w(G) = \w(G') + 26$, a contradiction. This completes the proof of Claim~\ref{claim.12.1}.~\smallqed

\medskip
By Claim~\ref{claim.12.1}, the vertex $v_3$ is the only common neighbor of $u$ and $v$. Let $x$ be the neighbor of $u$ different from $v_1$ and $v_3$, and let $y$ be the neighbor of $v$ different from $v_2$ and $v_3$.

\begin{subclaim}
\label{claim.12.2}
The vertices $x$ and $y$ are not adjacent.
\end{subclaim}
\proof Suppose that $x$ and $y$ are adjacent. By Claim~\ref{claim.11}, at least one of $x$ and $y$ is large. Renaming vertices if necessary, we assume that $y \in \cL$. Let $y_1$ be neighbor of $y$ different from $x$ and $v$.

Suppose that $x \in \cS$. If $y_1 \in \cL$, then we let $Q = \{u,v,v_1,v_2,v_3,x,y\}$ and $G' = G - Q$. The graph $G'$ is a connected subcubic graph. We note that $k'+r'=1$. Applying Claim~\ref{claim.key2} with $H = G[Q]$ and $S_H = \{v_1,v_3,y\}$, we have $\w(G) < \w(G') + 10p$ where $p = \gamma_r(H) - k' = 3 - k'$. On the other hand, $\w(G) \ge 32 + ( \w(G') - 5k' - r' ) = \w(G') + 4p + 19$. Therefore, $\w(G') + 4p + 19 \le \w(G) < \w(G') + 10p$, and so $19 < 6p$, that is, $p \ge 4$. However, $p = 3 - k' \le 3$, a contradiction. Hence, $y_1 \in \cS$. Let $y_2$ be the neighbor of $y_1$ different from $y$. If $y_2 \in \cL$, then let $G' = G - \{u,v,v_1,v_2,v_3,x,y,y_1\}$. In this case, $\gamma_r(G) \le \gamma_r(G') + 3$, implying that $\w(G) < 10\gamma_r(G) \le \w(G') + 30$. However, $\w(G) \ge 37 - (\w(G') - 1 - 4) \ge \w(G') + 32$, a contradiction. Hence, $y_2 \in \cS$. Let $y_3$ be the neighbor of $y_2$ different from $y_1$. By Claim~\ref{claim.3}, $y_3 \in \cL$. In this case, let $Q = \{u,v,v_1,v_2,v_3,x,y,y_1,y_2\}$ and let $G' = G - Q$. Applying Claim~\ref{claim.key2} with $H = G[Q]$ and $S_H = \{v_2,v_3,x,y_2\}$, we obtain a contradiction.

Hence, $x \in \cL$. We now consider the graph $G' = G - \{u,v,v_1,v_2,v_3\}$. Let $S'$ be a $\gamma_r$-set of $G'$. In this case, $\gamma_r(G) \le \gamma_r(G') + 2$, implying that $\w(G) < 10\gamma_r(G) \le \w(G') + 20$. If $G \notin \cB_{\rdom}$, then $\w(G) \ge 23 + (\w(G') - 2) = \w(G') + 21$, a contradiction. Hence, $G \in \cB_{\rdom}$. We note that $x$ and $y$ are adjacent vertices of degree~$2$ in $G'$. Applying Observation~\ref{obser-1}(f) to the graph $G'$ with $X = \{x,y\}$, we have $\gamma_{r,\dom}(G';X) \le \gamma_{r}(G') - 1$. Let $S''$ be a minimum type-$2$ NeRD-set of $G'$ with respect to the set $X$. The set $S'' \cup \{v_1,v_3\}$ is a RD-set of $G$, implying that $\gamma_r(G) \le |S''| + 2 \le \gamma_{r}(G') + 1$ and $\w(G) < \w(G') + 10$. However, $\w(G) \ge 23 + (\w(G') - 2 - 4) = \w(G') + 15$, a contradiction.~\smallqed

\medskip
By Claim~\ref{claim.12.2}, the vertices $x$ and $y$ are not adjacent.  Let $G'$ be obtained from $G - \{u,v,v_1,v_2,v_3\}$ by adding the edge $xy$. Suppose that $G' \notin \cB_{\rdom}$, implying that $\w(G) \ge 23 + \w(G')$. Let $S'$ be a $\gamma_r$-set of $G'$. If $x \in S'$, let $S = \{v,v_2\}$. If $x \notin S'$ and $y \in S'$, let $S = \{u,v_1\}$. If $x \notin S'$ and $y \notin S'$, let $S = \{v_1,v_3\}$. In all cases, $S$ is a RD-set of $G$, and so $\gamma_r(G) \le |S| + 2 = \gamma_r(G') + 2$, implying that $\w(G) < 10\gamma_r(G) \le \w(G') + 20$, a contradiction. Hence, $G' \in \cB_{\rdom}$. In this case, the graph $G$ is determined and it is straightforward to check (or can be checked by computer) that $10\gamma_r(G) \le \w(G)$, a contradiction.~\smallqed 

\begin{claim}
\label{claim.bridge}
The removal of a bridge joining two large vertices cannot create a component that belongs to $\cB_{\rdom}$.
\end{claim}
\proof Let $e = xy$ be a bridge in $G$ joining two adjacent large vertices $x$ and $y$. Let $G_x$ and $G_y$ be the components of $G - e$ containing $x$ and $y$, respectively. We note that both $G_x$ and $G_y$ are connected special subcubic graphs. Suppose, to the contrary, that at least one of $G_x$ and $G_y$ belongs to $\cB_{\rdom}$. Renaming components if necessary, we may assume that $G_y \in \cB_{\rdom}$.

Suppose that $G_x \in \cB_{\rdom}$. Since there is no handle in $G$, we note that $G_x \ne R_1$ and $G_y \ne R_1$. Therefore, $\w(G) \ge (\w(G_x) - 1 - 4) + (\w(G_y) - 1 - 4) = \w(G_x) + \w(G_y) - 10$. By Observation~\ref{obser-1}(b) there exists a $\gamma_r$-set $S_x$ of $G_x$ that contains~$x$.
A type-$1$ NeRD-set of $G_y$ with respect to the vertex~$y$ can be extended to a RD-set of $G$ by adding to it the set $S_x$. Hence by Observation~\ref{obser-1}(d), $\gamma_r(G) \le \gamma_r(G_x) + \gamma_{r,\ndom}(G;y) \le \gamma_{r}(G_x) + \gamma_{r}(G_y) - 1$.
Hence, $10\gamma_r(G) \le 10(\gamma_{r}(G_x) + \gamma_{r}(G_y) - 1) \le \w(G_x) + \w(G_y) - 10 \le \w(G)$, a contradiction.

Hence, $G_x \notin \cB_{\rdom}$. By Claim~\ref{claim.12} if $G_y = R_2$, then the vertex $y$ cannot be one of the two open twins in $R_2$.
Let $S_x$ be a $\gamma_r$-set of $G_x$. If $x \in S_x$, then let $S_y$ be a minimum type-$1$ NeRD-set of $G_y$ with respect to the vertex $y$. In this case, the set $S_x \cup S_y$ is a RD-set of $G$, implying by Observation~\ref{obser-1}(d) that $\gamma_r(G) \le |S_x| + |S_y| \le \gamma_r(G_x) + \gamma_{r,\ndom}(G;y) \le \gamma_r(G_x) + \gamma_r(G_y) - 1$. If $x \notin S_x$, then let $S_y$ is a minimum type-$2$ NeRD-set of $G_y$ with respect to the vertex $y$. In this case, the set $S_x \cup S_y$ is a RD-set of $G$, implying by Observation~\ref{obser-1}(e) that $\gamma_r(G) \le |S_x| + |S_y| \le \gamma_r(G_x) + \gamma_{r,\dom}(G;y) \le \gamma_r(G_x) + \gamma_r(G_y) - 1$. In both cases, $\gamma_r(G) \le \gamma_r(G_x) + \gamma_r(G_y) - 1$, implying that $\w(G) <  \w(G_x) + \w(G_y) - 10$. However, $\w(G) \ge (\w(G_x) - 1) + (\w(G_y) - 1 - 4) = \w(G_x) + \w(G_y) - 6$, a contradiction.~\smallqed

\begin{claim}
\label{claim.2linkage-bad}
The removal of the two small vertices on a $2$-linkage cannot create a component that belongs to $\cB_{\rdom}$.
\end{claim}
\proof Let $P \colon v v_1 v_2 u$ be a $2$-linkage, and so $u,v \in \cL$ and $v_1,v_2 \in \cS$. By Claim~\ref{claim.9}, $uv \notin E(G)$. Suppose, to the contrary, that $G' = G - \{v_1,v_2\}$ creates a component that belongs to~$\cB_{\rdom}$. Let $G_u$ and $G_v$ be the components of $G - e$ containing $u$ and~$v$, respectively, where we may assume renaming vertices if necessary, that $G_v \in \cB_{\rdom}$. 
Suppose that $G_u = G_v$, and so the graph $G'$ is connected. In this case, let $S_v$ be a minimum type-$1$ NeRD-set of $G_v$ with respect to the vertex $v$. The set $S_v \cup \{v_1\}$ is a RD-set of $G$, implying by Observation~\ref{obser-1} that $\gamma_r(G) \le 1 + \gamma_{r,\ndom}(G;v) \le \gamma_{r}(G_v)$. Hence, $\w(G) < 10\gamma_r(G) \le \w(G_v)$. However, $\w(G) = 10 + (\w(G_v) - 2 - 4) = \w(G_v) + 4$, a contradiction. Hence, $G_u \ne G_v$, and so $G'$ is disconnected with two components $G_u$ and $G_v$.

Let $S_u$ be a $\gamma_r$-set of $G_v$. Suppose that $u \in S_u$. In this case, the set $S_u$ can be extended to a RD-set of $G$ by adding to it a $\gamma_r$-set of $G_v$ that contains~$v$, which exists by Observation~\ref{obser-1}(d), implying that $\gamma_r(G) \le \gamma_r(G_u) + \gamma_r(G_v)$. Suppose that $u \notin S_u$. By Observation~\ref{obser-1}(d), $\gamma_{r,\ndom}(G_v;v) \le \gamma_{r}(G_v) - 1$. In this case, the set $S_u$ can be extended to a RD-set of $G$ by adding to it the vertex $v_2$ and a minimum type-$1$ NeRD-set of $G_v$ with respect to the vertex $v$, implying that $\gamma_r(G) \le |S_u| + 1 + \gamma_{r,\ndom}(G_v;v) \le \gamma_r(G_u) + \gamma_{r}(G_v)$. Thus in both cases, $\gamma_r(G) \le \gamma_r(G_u) + \gamma_r(G_v)$, implying that $\w(G) < \w(G_u) + \w(G_v)$. However, $\w(G) \ge 10 + (\w(G_u) - 1 - 4) + (\w(G_v) - 1 - 4) = \w(G_u) + \w(G_v)$, a contradiction.~\smallqed

\begin{claim}
\label{claim.1linkage-bad}
The removal of the small vertex on a $1$-linkage cannot create a component that belongs to $\cB_{\rdom}$.
\end{claim}
\proof Let $P \colon v v_1 u$ be a $1$-linkage, and so $u,v \in \cL$ and $v_1 \in \cS$. By Claim~\ref{claim.9}, $uv \notin E(G)$. Suppose, to the contrary, that $G' = G - v_1$ creates a component that belongs to~$\cB_{\rdom}$. Let $G_u$ and $G_v$ be the components of $G - v_1$ containing $u$ and~$v$, respectively, where we may assume renaming vertices if necessary, that $G_v \in \cB_{\rdom}$. Let $u_1$ and $u_2$ be the two neighbors of $u$ different from $v_1$. Let $S_v^1$ be a minimum type-$1$ NeRD-set of $G_v$ with respect to the vertex~$v$. By Observation~\ref{obser-1}(d), $|S_v^1| = \gamma_{r,\ndom}(G_v;v) \le \gamma_{r}(G_v) - 1$. Let $S_v^2$ be a minimum type-$2$ NeRD-set of $G_v$ with respect to the vertex~$v$. By Claim~\ref{claim.12}, if $G_v = R_2$, then the vertex $v$ is not one of the open twins in $G_v$, implying by Observation~\ref{obser-1}(e) that $|S_v^2| = \gamma_{r,\dom}(G_v;v) \le \gamma_{r}(G_v) - 1$.

Suppose that $u_1u_2 \in E(G)$. By Claim~\ref{claim.11a}, $u_1, u_2 \in \cL$. Let $Q = V(G_v) \cup \{u,v_1\}$ and let $G' = G - Q$. Suppose that $G' \in \cB_{\rdom}$. In this case, let $S'$ be a minimum type-$1$ NeRD-set of $G'$ with respect to the vertex $u_1$. By Observation~\ref{obser-1}(d), $|S'| \le \gamma_{r,\ndom}(G';u_1) \le \gamma_r(G') - 1$. The set $S' \cup \{u\} \cup S_v^2$ is a RD-set of $G$, and so $\gamma_r(G) \le |S'| + 1 + |S_v^2| \le (\gamma_r(G') - 1) + 1 + (\gamma_{r}(G_v) - 1) = \gamma_r(G') + \gamma_{r}(G_v) - 1$, implying that $\w(G) < \w(G') + \w(G_v) - 10$. However, $\w(G) \ge 9 + (\w(G') - 2 - 4) + (\w(G_v) - 1 - 4) = \w(G') + \w(G_v) - 2$, a contradiction. Hence, $G' \notin \cB_{\rdom}$. Thus, $\w(G) \ge 9 + (\w(G') - 2) + (\w(G_v) - 1 - 4) = \w(G') + \w(G_v) + 2$. Every $\gamma_r$-set of $G'$ can be extended to a RD-set of $G$ by adding to it the set $S_v^2 \cup \{u\}$, implying that $\gamma_r(G) \le \gamma_r(G') + |S_v^2| + 1 \le \gamma_r(G') + (\gamma_{r}(G_v) - 1) + 1 = \gamma_r(G') + \gamma_{r}(G_v)$. Thus, $\w(G) < \w(G') + \w(G_v) \le \w(G) - 2 < \w(G)$, a contradiction.

Hence, $u_1u_2 \notin E(G)$. Let $Q = V(G_v) \cup \{u,v_1\}$ and let $G'$ be obtained from $G - Q$ by adding the edge $u_1u_2$. The resulting graph $G'$ is a connected subcubic graph.
Suppose that $G' \in \cB_{\rdom}$. In this case, let $Q^* = Q \setminus \{u\}$, and let $G^* = G - Q^*$, and so $G^*$ is obtained from $G'$ by subdividing the added edge $u_1u_2$ where $u$ is the resulting new vertex of degree~$2$ in $G^*$. By Observation~\ref{obser-2}, $\gamma_r(G^*) \le \gamma_r(G')$ and there exists a $\gamma_r$-set $S^*$ of $G^*$ that contains~$u$. The set $S^* \cup S_v^2$ is a RD-set of $G$, and so $\gamma_r(G) \le |S^*| + |S_v^2| \le \gamma_r(G') + \gamma_{r}(G_v) - 1$, implying that $\w(G) < \w(G') + \w(G_v) - 10$. However, $\w(G) \ge 9 + (\w(G') - 4) + (\w(G_v) - 1 - 4) = \w(G') + \w(G_v)$, a contradiction.
Hence, $G' \notin \cB_{\rdom}$. Thus, $\w(G) \ge 9 + \w(G') + (\w(G_v) - 1 - 4) = \w(G') + \w(G_v) + 4$. Let $S'$ be a $\gamma_r$-set of $G'$. If at least one of $u_1$ and $u_2$ belongs to $S'$, let $S = S' \cup \{u\} \cup S_v^2$. If $u_1 \notin S'$ and $u_2 \notin S'$, let $S = S' \cup \{v_1\} \cup S_v^1$. In both cases, $S$ is a RD-set of $G$ and $|S| \le |S'| + 1 + \gamma_{r}(G_v) - 1 = \gamma_r(G') + \gamma_{r}(G_v)$. Thus, $\w(G) < \w(G') + \w(G_v) \le \w(G) - 4 < \w(G)$, a contradiction.~\smallqed

\medskip
By our earlier observations, every edge of $G$ either joins two large vertices or belongs to a $2$-linkage or belongs to a $1$-linkage. Hence as an immediate consequence of Claims~\ref{claim.bridge},~\ref{claim.2linkage-bad}, and~\ref{claim.1linkage-bad}, we have the following property of the graph $G$.

\begin{claim}
\label{claim.bridge2}
The removal of a bridge cannot create a component that belongs to $\cB_{\rdom}$.
\end{claim}

As a consequence of Claim~\ref{claim.bridge2}, we have the following claim.

\begin{claim}
\label{claim.no.R10}
The graph $G$ does not contain $R_{11}$ as a subgraph.
\end{claim}
\proof Suppose, to the contrary, that $R'$ is a subgraph of $G$, where $R' = R_{11}$. Let $v$ be small vertex (of degree~$2$) in $R'$. Since $G \notin \cB_{\rdom}$, we note that $R' \ne G$, implying that $v$ is a large vertex in $G$. Let $v'$ be the vertex adjacent to $v$ that does not belong to~$R'$. The edge $vv'$ is a bridge of $G$ whose removal creates a $R_{11}$-component, contradicting Claim~\ref{claim.bridge2}.~\smallqed

\medskip
We are now in a position to prove that there is no $2$-linkage in $G$.

\begin{claim}
\label{claim.16}
There is no $2$-linkage in $G$.
\end{claim}
\proof Suppose, to the contrary, that $G$ contains a $2$-linkage. Let $P \colon u v_1 v_2 v$ be a $2$-linkage, where $x$ and $y$ are the two neighbors of $v$ not on $P$.

\begin{subclaim}
\label{claim.16.1}
At least one of $ux$ and $uy$ is not an edge.
\end{subclaim}
\proof Suppose, to the contrary, that $u$ is adjacent to both $x$ and $y$. By Claim~\ref{claim.12}, $x,y \in \cL$. If $xy$ is an edge, then the graph $G$ is determined and $\gamma_r(G) = 2$ and $\w(G) = 26$, a contradiction. Hence, $xy$ is not an edge. Let $x_1$ and $y_1$ be the neighbors of $x$ and $y$, and let $G' = G - \{u,v_1,v_2\}$. Suppose $G' \in \cB_{\rdom}$. In this case, $xvy$ is a path in $G'$ where $x,v$ and $y$ all have degree~$2$ in $G'$, implying that $G' \in \{R_1,R_3,R_9\}$. If $G' = R_1$, then $G = R_5$, a contradiction. If $G' = R_3$, then $G$ could contain a $3$-linkage, a contradiction. If $G' = R_9$, then $G$ is determined and $\gamma_r(G) \le 6$ and $\w(G) = 60$, a contradiction.
Hence, $G' \notin \cB_{\rdom}$, implying that $\w(G) = 14 + (\w(G') - 3) = \w(G') + 11$. Let $S'$ be a $\gamma_r$-set of $G'$. If $v \notin S$, then either $x \in S'$ and $y \notin S'$ or $x \notin S'$ and $y \in S'$. In this case, we let $S = S' \cup \{v_2\}$. If $v \in S'$ and neither $x$ not $y$ belongs to $S'$, then we let $S = S' \cup \{u\}$. If $v \in S'$ and at least one of $x$ and $y$ belongs to $S'$, then we let $S = S' \cup \{v_2\}$. In all cases, $S$ is a RD-set of $G$, and so $\gamma_r(G) \le |S| = |S'| + 1 = \gamma_r(G') + 1$, and so $\w(G) < \w(G') + 10$. This contradicts our earlier observation that $\w(G) = \w(G') + 11$.~\smallqed

\begin{subclaim}
\label{claim.16.2}
Neither $ux$ nor $uy$ is an edge.
\end{subclaim}
\proof Suppose, to the contrary, that $u$ is adjacent to exactly one of $x$ and $y$. We may assume that $uy$ is an edge. By Claim~\ref{claim.12}, $y \in \cL$. By Claim~\ref{claim.16.1}, $ux$ is not an edge. Let $G'$ be obtained from $G - \{v,v_1,v_2\}$ by adding the edge $ux$.  The graph $G'$ is a connected special subcubic graph of order less than~$n$. Let $S'$ be a $\gamma_r$-set of $G'$. If $u \in S$, then let $S = S' \cup \{v\}$. If $u \notin S'$ and $x \in S'$, then we let $S = S' \cup \{v_1\}$. If $u \notin S'$ and $x \notin S'$, then we let $S = S' \cup \{v_2\}$. In all cases, $S$ is a RD-set of $G$, and so $\gamma_r(G) \le |S| = |S'| + 1 = \gamma_r(G') + 1$, and so $\w(G) < \w(G') + 10$. If $G' \notin \{R_1,R_4,R_5\}$, then $\w(G) \ge 14 + (\w(G') - 1 - 3) = \w(G') + 10$, a contradiction. Hence, $G' \in \{R_1,R_4,R_5\}$. Since $u$ is a vertex of degree~$3$ in $G'$, we note that $G' \ne R_1$. If $G' = R_4$, then $G = R_8$, a contradiction. If $G' = R_5$, then $G = R_7$, a contradiction.~\smallqed

\medskip
By Claim~\ref{claim.16.2}, the vertex $u$ is adjacent to neither $x$ nor $y$. Renaming vertices if necessary, we may assume that $\deg_G(x) \le \deg_G(y)$.

\begin{subclaim}
\label{claim.16.4}
$x, y \in \cS$.
\end{subclaim}
\proof Suppose that $y \in \cL$. Let $G'$ be obtained from $G - \{v,v_1,v_2\}$ by adding the edge $ux$.  As in the proof of Claim~\ref{claim.16.2}, $\gamma_r(G) \le |S| = |S'| + 1 = \gamma_r(G') + 1$, implying that $\w(G) < \w(G') + 10$. If no component of $G'$ belongs to $\cB_{\rdom}$, then $\w(G) = \w(G') + 13$, a contradiction. Hence, at least one component in $G'$ belongs to $\cB_{\rdom}$. Suppose that $G'$ is connected. Since $u$ is a vertex of degree~$3$ in $G'$, we note that $G' \ne R_1$. If $G' \in \{R_4,R_5\}$, then the graph $G$ is determined and in all cases, $\gamma_r(G) = 4$ and $\w(G) = 49$, a contradiction. Hence, $G \notin \{R_1,R_4,R_5\}$, implying that $\w(G) \ge 14 + (\w(G') - 1 - 3) = \w(G') + 10$, a contradiction.

Hence, $G'$ is disconnected with two components. Let $G_x$ be the component of $G'$ containing the vertices $u$ and $x$, and let $G_y$ be the component containing the vertex~$y$. Both $G_x$ and $G_y$ are connected special subcubic graphs. Further we note that the edge $vy$ is a bridge in $G$, implying by Claim~\ref{claim.bridge2} that $G_y \notin \cB_{\rdom}$ and therefore $G_x \in \cB_{\rdom}$. Let $G_v$ be the component of $G - vy$ that contains the vertex~$v$. Thus, $G_v$ is obtained from the graph $G'$ by subdividing the edge $ux$ three times resulting in the path $u v_1 v_2 v x$.

Suppose that $G' = R_6$. In this case, $\gamma_r(G_x) \le 5$ and $\w(G_x) = 54$, and so $10\gamma_r(G_x) \le \w(G_x) - 4$. By the minimality of $G$, we have $10\gamma_r(G_y) \le \w(G_y)$. Thus, $10\gamma_r(G) \le 10(\gamma_r(G_x) + \gamma_r(G_y)) \le \w(G_x) + \w(G_y) - 4 = \w(G) - 6$, a contradiction. Hence, $G' \ne R_6$. Let $S_v^1$ be a minimum type-$1$ NeRD-set of $G_v$ with respect to the vertex~$v$, and let $S_v^2$ be a minimum type-$2$ NeRD-set of $G_v$ with respect to the vertex $v$. By Observation~\ref{obser-4}, $|S_v^1| = \gamma_{r,\ndom}(G^*;v_1) \le \gamma_{r}(G')$ and $|S_v^2| = \gamma_{r,\dom}(G^*;v_1) \le \gamma_{r}(G')$. Let $S_y$ be a $\gamma_r$-set of $G_y$. If $y \in S_y$, let $S = S_y \cup S_v^1$. If $y \notin S_y$, let $S = S_y \cup S_v^2$. In both cases, $S$ is a RD-set of $G$, implying that $\gamma_r(G) \le \gamma_r(G_y) + \gamma_{r}(G')$, and so $\w(G) < \w(G_y) + \w(G')$. However, $\w(G) \ge 14 + (\w(G') - 4) + (\w(G_y) - 1) = \w(G_y) + \w(G') + 9$, a contradiction. Hence, $y \in \cS$. By our choice of the vertex~$x$, this implies that $x \in \cS$.~\smallqed

\medskip
By Claim~\ref{claim.16.4}, $x \in \cS$ and $y \in \cS$. Thus, all three neighbors of $v$ are small vertices. Interchanging the roles of $u$ and $v$, analogous arguments show that all three neighbors of $u$ are small vertices. Recall that $ux \notin E(G)$ and $uy \notin E(G)$, and so $u$ and $v$ do not have a common neighbor. By Claim~\ref{claim.11a}, no small vertex belongs to a triangle, implying that $xy \notin E(G)$. Let $x_1$ and $y_1$ be the neighbors of $x$ and $y$, respectively, different from~$v$. Possibly, $x_1 = y_1$.

\begin{subclaim}
\label{claim.16.5}
$x_1 \ne y_1$.
\end{subclaim}
\proof Suppose, to the contrary, that $x_1 = y_1$. In this case, we let $z = x_1$. Since $G$ has no handle, $z \in \cL$. Thus, $C \colon vxzyv$ is a $4$-cycle in $G$, where $v,z \in \cL$ and $x,y \in \cS$. Let $z_1$ be the neighbor of $z$ different from $x$ and $y$. Since all three neighbors of $u$ belong to $\cS$, we note that $uz \notin E(G)$. Thus, $u \ne z_1$. Let $G'$ be obtained from $G - \{v,v_2,x,y,z\}$ by adding the edge $v_1z_1$. The resulting graph $G'$ is a connected subcubic graph.

Suppose that $G' \in \cB_{\rdom}$. Let $G^* = G - y$, that is, $G^*$ is obtained from $G'$ by subdividing the added edge $v_1z_1$ four times resulting in the path $v_1v_2vxzz_1$. By Observation~\ref{obser-5}, there exists a RD-set $S^*$ of $G^*$ such that $|S^*| \le \gamma_{r}(G') + 1$ and $S^* \cap \{v_2,v,x,z\} = \{v_2,z\}$. The set $S^*$ is a RD-set of $G$, and so $\gamma_r(G) \le |S^*| = \gamma_{r}(G') + 1$, implying that $\w(G) < \w(G') + 10$. Noting that $G' \ne R_1$ and the degrees of the vertices in $G'$ are the same as their degrees in $G$, we have $\w(G) \ge 23 + (\w(G') - 4) = \w(G') + 19$, a contradiction.
Hence, $G' \notin \cB_{\rdom}$, implying that $\w(G) = 22 + \w(G')$. Let $S'$ be a $\gamma_r$-set of $G'$. If at least one of $v_1$ and $z_1$ belongs to $S'$, let $S = S' \cup \{v_2,z\}$. If $u_1 \notin S'$ and $z_1 \notin S'$, let $S = S' \cup \{v,x\}$. In both cases, $S$ is a RD-set of $G$ and $|S| \le |S'| + 2 = \gamma_r(G') + 2$. Thus, $\w(G) < \w(G') + 20$, a contradiction.~\smallqed

\medskip
By Claim~\ref{claim.16.5}, $x_1 \ne y_1$. By Claim~\ref{claim.11} if $x_1 \in \cS$, then $ux_1 \notin E(G)$ and if $y_1 \in \cS$, then $uy_1 \notin E(G)$.

\begin{subclaim}
\label{claim.16.6}
$x_1y_1 \in E(G)$.
\end{subclaim}
\proof Suppose that $x_1y_1 \notin E(G)$. Let $G'$ be obtained from $G - \{v,v_1,v_2,x,y\}$ by adding the edge $x_1y_1$. The resulting graph $G'$ is a subcubic graph. We note that either $G'$ is connected or has two components. Let $G_{xy}$ be the component of $G'$ containing the added edge $x_1y_1$. If $G'$ is disconnected, then let $G_u$ be the second component of $G'$ which necessarily contains the vertex~$u$. In this case the edge $uv_1$ is a bridge in $G$, implying by Claim~\ref{claim.bridge2} that $G_u \notin \cB_{\rdom}$. Therefore, the component $G_{xy}$ is the only possible component of $G'$ that belongs to~$\cB_{\rdom}$.

Let $S'$ be a $\gamma_r$-set of $G'$. If $x_1 \in S'$, let $S = S' \cup \{v_1,y\}$. If $x_1 \notin S'$ and $y_1 \in S'$, let $S = S' \cup \{v_1,x\}$. If $x_1 \notin S'$, $y_1 \notin S'$ and $u \in S'$, let $S = S' \cup \{v\}$. If $x_1 \notin S'$, $y_1 \notin S'$ and $u \notin S'$, let $S = S' \cup \{v,v_2\}$. In all cases, $S$ is a RD-set of $G$ and $|S| \le |S'| + 2 = \gamma_r(G') + 2$. Thus, $\w(G) < 10\gamma_r(G) \le \w(G') + 20$. If $G'$ has no component in $\cB_{\rdom}$, then $\w(G) \ge 24 + (\w(G') - 1) = \w(G') + 23$, a contradiction. Hence by our earlier observations, $G'$ has exactly one component in $\cB_{\rdom}$, namely the component $G_{xy}$. By our earlier properties of the graph $G$, we note that $G_{xy} \ne R_1$. If $G_{xy} \notin \{R_4,R_5\}$, then $\w(G) \ge 24 + (\w(G') - 1 - 3) = \w(G') + 20$, a contradiction. Hence, $G_{xy} \in \{R_4,R_5\}$, implying that $\w(G) = 24 + (\w(G') - 1 - 4) = \w(G') - 19$.

If $G'$ is connected, then $G' = G_{xy}$ and we let $G^* = G - \{v_1,v_2\}$. If $G'$ is disconnected, then $G'$ consists of the two components $G_u$ and $G_{xy}$ and we let $G^* = G - (V(G_u) \cup \{v_1,v_2\})$. In both cases, $G^*$ is the graph obtained from $G_{xy}$ by subdividing the added edge $x_1y_1$ three times resulting in the path $x_1 x v y y_1$. Recall that $G_{xy} \in \{R_4,R_5\}$. Applying Observation~\ref{obser-4} we have $\gamma_{r,\dom}(G^*;v) \le \gamma_{r}(G_{xy})$. Thus, there exists a type-$2$ NeRD-set $S^*$ in $G^*$ with respect to the vertex $v$ such that $|S^*| \le \gamma_{r}(G_{xy})$. If $G'$ is connected, then let $S = S^*$, and note that in this case, $|S| \le \gamma_{r,\dom}(G^*;v) \le \gamma_{r}(G_{xy}) = \gamma_{r}(G')$. If $G'$ is disconnected, let $S = S^* \cup S_u$ where $S_u$ is a $\gamma_r$-set of $G_u$, and note that in this case, $|S| \le \gamma_{r,\dom}(G^*;v) + \gamma_r(G_u) \le \gamma_{r}(G_{xy} + \gamma_r(G_u) = \gamma_{r}(G')$. In both cases, $|S| \le \gamma_{r}(G')$. Further in both cases $S \cup \{v_1\}$ is a RD-set of $G$, implying that $\gamma_r(G) \le |S| + 1 \le \gamma_{r}(G') + 1$. Hence, $\w(G) < 10\gamma_r(G) \le \w(G') + 10$, a contradiction.~\smallqed

\medskip
By Claim~\ref{claim.16.6}, $x_1y_1 \in E(G)$. Since $G$ has no handle, at least one of $x_1$ and $y_1$ is large. If exactly one of $x_1$ and $y_1$ is large, then we would contradict Claim~\ref{claim.12}. Hence, $x_1 \in \cL$ and $y_1 \in \cL$. Let $G' = G - \{v,v_1,v_2,x,y\}$. The resulting graph $G'$ is a special subcubic graph. Let $G_{xy}$ be the component of $G'$ containing the edge $x_1y_1$, and let $G_u$ be the component of $G'$ containing the vertex~$u$. If $G'$ is connected, then $G_u = G_{xy}$. If $G'$ is disconnected, then $G_{xy}$ and $G_u$ are the two components of $G'$. Further in this case, $uv_1$ is a bridge in $G$, implying by Claim~\ref{claim.bridge2} that $G_u \notin \cB_{\rdom}$. Therefore, the component $G_{xy}$ is the only possible component of $G'$ that belongs to~$\cB_{\rdom}$.

Let $S'$ be a $\gamma_r$-set of $G'$. If $\{u,x_1,y_1\} \subseteq S'$, let $S = S' \cup \{v_1,v_2\}$. If $S' \cap \{u,x_1,y_1\} = \{u,y_1\}$, let $S = S' \cup \{v_1,x\}$.  If $S' \cap \{u,x_1,y_1\} = \{u,y_1\}$, let $S = S' \cup \{v_1,y\}$. If $S' \cap \{u,x_1,y_1\} = \{u\}$, let $S = S' \cup \{v\}$. If $S' \cap \{u,x_1,y_1\} = \{x_1,y_1\}$, let $S = S' \cup \{v_2\}$. If $S' \cap \{u,x_1,y_1\} = \{x_1\}$, let $S = S' \cup \{v_2,y\}$. If $S' \cap \{u,x_1,y_1\} = \{y_1\}$, let $S = S' \cup \{v_2,x\}$. If $S' \cap \{u,x_1,y_1\} = \emptyset$, let $S = S' \cup \{v,v_2\}$. In all cases, $S$ is a RD-set of $G$ and $|S| \le |S'| + 2 = \gamma_r(G') + 2$. Thus, $\w(G) < 10\gamma_r(G) \le \w(G') + 20$.

If $G_{xy} \notin \cB_{\rdom}$ or if $G_{xy} \in \cB_{\rdom,1}$, then $\w(G) \ge 24 + (\w(G') - 3 - 1) = \w(G') + 20$, a contradiction. Hence, $G_{xy} \in \{R_1,R_2,R_3,R_4,R_5,R_{10}\}$. We note that $x_1$ and $y_1$ are adjacent vertices of degree~$2$ in $G_{xy}$, implying that $G_{xy} \notin \{R_4, R_{10}\}$.
If $G_{xy} = R_1$, then necessarily $G' = G_{xy}$ and the graph $G$ is determined. In this case, $\gamma_r(G) = 4$ and $\w(G) = 46$, a contradiction. Hence, $G_{xy} \ne R_1$.
If $G_{xy} = R_3$, then since $G$ has no $3$-linkage, $G' = G_{xy}$ and the graph $G$ is determined. In this case, $\gamma_r(G) = 5$ and $\w(G) = 59$, a contradiction. Hence, $G_{xy} \ne R_3$. Therefore, $G_{xy} \in \{R_2,R_5\}$. By our earlier observations, the vertex $u$ is adjacent to neither $x_1$ nor $y_1$. Further, we note that the vertex $u$ and its two neighbors in $G'$, as well as $x_1$ and $y_1$, all have degree~$2$ in $G'$. Moreover, $x_1y_1$ in an edge. These properties implies that $G'$ is disconnected. Thus, $G'$ has two components, namely $G_{xy}$ and $G_u$. As observed earlier, $G_u \notin \cB_{\rdom}$. Let $x_2$ and $y_2$ be the neighbors of $x_1$ and $y_1$, respectively, in $G_{xy}$.

Suppose that $G_{xy} = R_2$. We note that $x_2$ and $y_2$ are the two large vertices in $R_2$. Let $z_1$ and $z_2$ be the two common neighbors of $x_2$ and $y_2$ in $G_{xy}$. Thus, the graph in Figure~\ref{rdom:fig-20} is a subgraph of $G$. Let $S_u$ be a $\gamma_r$-set of $G_u$. If $u \in S_u$, let $S = \{v,x,y_2,z_1\}$. If $u \notin S_u$, let $S = \{v_2,x_1,x_2,y_1\}$. In both cases, $S$ is a RD-set of $G$, and so $\gamma_r(G) \le |S_u| + 4 = \gamma_r(G_u) + 4$. Hence, $\w(G) < \w(G_u) + 40$. However, $\w(G) = 50 + (\w(G_u) - 1) = \w(G_u) + 49$, a contradiction.

\begin{figure}[htb]
\begin{center}
\begin{tikzpicture}[scale=.8,style=thick,x=0.8cm,y=0.8cm]
\def\vr{2.5pt} 
\path (0,1) coordinate (u1);
\path (-0.5,0.5) coordinate (a);
\path (-0.5,1.5) coordinate (b);
\path (1,1) coordinate (u2);
\path (2,1) coordinate (u3);
\path (3,1) coordinate (u4);
\path (4,0) coordinate (u5);
\path (4,2) coordinate (u6);
\path (5,0) coordinate (u7);
\path (5,2) coordinate (u8);
\path (7,0) coordinate (u9);
\path (7,2) coordinate (u10);
\path (6,1) coordinate (u11);
\path (8,1) coordinate (u12);
\draw (u1)--(u2)--(u3)--(u4)--(u5)--(u7)--(u8)--(u6)--(u4);
\draw (u7)--(u9)--(u11)--(u10)--(u12)--(u9);
\draw (u8)--(u10);
\draw (a)--(u1)--(b);
\draw (u1) [fill=white] circle (\vr);
\draw (u2) [fill=white] circle (\vr);
\draw (u3) [fill=white] circle (\vr);
\draw (u4) [fill=white] circle (\vr);
\draw (u5) [fill=white] circle (\vr);
\draw (u6) [fill=white] circle (\vr);
\draw (u7) [fill=white] circle (\vr);
\draw (u8) [fill=white] circle (\vr);
\draw (u9) [fill=white] circle (\vr);
\draw (u10) [fill=white] circle (\vr);
\draw (u11) [fill=white] circle (\vr);
\draw (u12) [fill=white] circle (\vr);
%
%
\draw[anchor = south] (u1) node {{\small $u$}};
\draw[anchor = south] (u2) node {{\small $v_1$}};
\draw[anchor = south] (u3) node {{\small $v_2$}};
\draw[anchor = south] (u4) node {{\small $v$}};
\draw[anchor = north] (u5) node {{\small $y$}};
\draw[anchor = south] (u6) node {{\small $x$}};
\draw[anchor = north] (u7) node {{\small $y_1$}};
\draw[anchor = south] (u8) node {{\small $x_1$}};
\draw[anchor = north] (u9) node {{\small $y_2$}};
\draw[anchor = south] (u10) node {{\small $x_2$}};
\draw[anchor = east] (u11) node {{\small $z_1$}};
\draw[anchor = west] (u12) node {{\small $z_2$}};

\end{tikzpicture}
\end{center}
\begin{center}
\vskip -0.5 cm
\caption{A subgraph in the proof of Claim~\ref{claim.16} when $G_{xy} = R_2$}
\label{rdom:fig-20}
\end{center}
\end{figure}
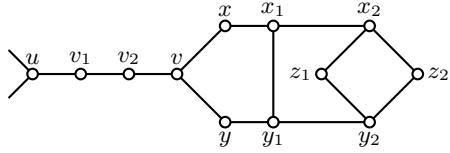

Hence, $G_{xy} \ne R_2$, and so $G_{xy} = R_5$. Let $x_3$ and $y_3$ be the two common neighbors of $x_2$ and $y_2$ in $G_{xy}$, and let $x_4$ and $y_4$ be the remaining vertices in $G_{xy}$, where $x_3x_4y_4y_3$ is a path. Thus, the graph in Figure~\ref{rdom:fig-21} is a subgraph of $G$. Let $S_u$ be a $\gamma_r$-set of $G_u$, and let $S = S_u \cup \{v_1,x,y,x_3,y_3\}$. The set $S$ is a RD-set of $G$, and so $\gamma_r(G) \le |S_u| + 5 = \gamma_r(G_u) + 5$. Hence, $\w(G) < \w(G_u) + 50$. However, $\w(G) = 58 + (\w(G_u) - 1) = \w(G_u) + 57$, a contradiction. This completes the proof of Claim~\ref{claim.16}.~\smallqed

\begin{figure}[htb]
\begin{center}
\begin{tikzpicture}[scale=.8,style=thick,x=0.8cm,y=0.8cm]
\def\vr{2.5pt} 
\path (0,1) coordinate (u1);
\path (-0.5,0.5) coordinate (a);
\path (-0.5,1.5) coordinate (b);
\path (1,1) coordinate (u2);
\path (2,1) coordinate (u3);
\path (3,1) coordinate (u4);
\path (4,0) coordinate (u5);
\path (4,2) coordinate (u6);
\path (7.25,0.00) coordinate (x1);
\path (8.13,0.37) coordinate (x2);
\path (8.50,1.25) coordinate (x3);
\path (8.13,2.13) coordinate (x4);
\path (7.25,2.50) coordinate (x5);
\path (6.37,2.13) coordinate (x6);
\path (6.00,1.25) coordinate (x7);
\path (6.00,1.15) coordinate (x7p);
\path (6.37,0.37) coordinate (x8);
\draw (u1)--(u2)--(u3)--(u4)--(u5);
\draw (u4)--(u6);
\draw (a)--(u1)--(b);
\draw (x1)--(x2)--(x3)--(x4)--(x5)--(x6)--(x7)--(x8)--(x1);
\draw (x1)--(x5);
\draw (x2)--(x6);
\draw (u5)--(x8);
\draw (u6)--(x7);
\draw (u1) [fill=white] circle (\vr);
\draw (u2) [fill=white] circle (\vr);
\draw (u3) [fill=white] circle (\vr);
\draw (u4) [fill=white] circle (\vr);
\draw (u5) [fill=white] circle (\vr);
\draw (u6) [fill=white] circle (\vr);
%
\draw (x1) [fill=white] circle (\vr);
\draw (x2) [fill=white] circle (\vr);
\draw (x3) [fill=white] circle (\vr);
\draw (x4) [fill=white] circle (\vr);
\draw (x5) [fill=white] circle (\vr);
\draw (x6) [fill=white] circle (\vr);
\draw (x7) [fill=white] circle (\vr);
\draw (x8) [fill=white] circle (\vr);
\draw[anchor = south] (u1) node {{\small $u$}};
\draw[anchor = south] (u2) node {{\small $v_1$}};
\draw[anchor = south] (u3) node {{\small $v_2$}};
\draw[anchor = south] (u4) node {{\small $v$}};
\draw[anchor = north] (u5) node {{\small $y$}};
\draw[anchor = south] (u6) node {{\small $x$}};
\draw[anchor = east] (x7p) node {{\small $x_1$}};
\draw[anchor = north] (x8) node {{\small $y_1$}};
\draw[anchor = north] (x1) node {{\small $y_2$}};
\draw[anchor = east] (x6) node {{\small $x_2$}};
\draw[anchor = north] (x2) node {{\small $y_3$}};
\draw[anchor = south] (x5) node {{\small $x_3$}};
\draw[anchor = west] (x4) node {{\small $x_4$}};
\draw[anchor = west] (x3) node {{\small $y_4$}};
\end{tikzpicture}

\end{center}
\begin{center}
\vskip -0.45 cm
\caption{A subgraph in the proof of Claim~\ref{claim.16} when $G_{xy} = R_5$}
\label{rdom:fig-21}
\end{center}
\end{figure}
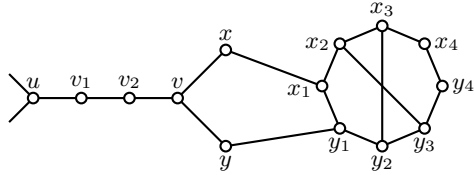

By Claim~\ref{claim.16}, there is no $2$-linkage in $G$. By our earlier observations, every vertex of degree~$2$ in $G$, if any, therefore belongs to a $1$-linkage.

\begin{claim}
\label{claim.18}
Two large vertices cannot be the ends of two common $1$-linkages.
\end{claim}
\proof Suppose, to the contrary, that there are two large vertices $u$ and $v$ such that $u v_1 v$ and $v v_2 u$ are $1$-linkages in $G$. Thus, $C \colon u v_1 v v_2 u$ is a $4$-cycle in $G$, where $u,v \in \cL$ and $v_1,v_2 \in \cS$.

\begin{subclaim}
\label{claim.18.1}
The vertices $v_1$ and $v_2$ are the only two common neighbors of $u$ and $v$.
\end{subclaim}
\proof Suppose that $u$ and $v$ have a third common neighbor $v_3$.
If $v_3 \in \cS$, then $\gamma_r(G) = 2$ and $\w(G) = 23$, a contradiction. Hence, $v_3 \in \cL$. Let $v_4$ be the neighbor of $v_3$ different from $u$ and $v$. If $v_4 \in \cL$, then let $G' = G - \{u,v,v_1,v_2,v_3\}$. By Claim~\ref{claim.bridge2}, $G' \notin \cB_{\rdom}$. Every $\gamma_r$-set of $G'$ can be extended to a RD-set of $G$ by adding to it $v$ and $v_1$, and so $\gamma_r(G) \le \gamma_r(G') + 2$, implying that $\w(G) < \w(G') + 20$. However, $\w(G) = 22 + (\w(G') - 1) = \w(G') + 21$, a contradiction. Hence, $v_4 \in \cS$. Let $v_5$ be the neighbor of $v_4$ different from $v_3$. If $v_5 \in \cL$, then let $G' = G - \{u,v,v_1,v_2,v_3,v_4\}$. By Claim~\ref{claim.bridge2}, $G' \notin \cB_{\rdom}$. Let $S'$ be a $\gamma_r$-set of $G'$. If $v_5 \in S'$, let $S = S' \cup \{v,v_1\}$. If $v_5 \notin S'$, let $S = S' \cup \{v,v_3\}$. In both cases, $S$ is a RD-set of $G$, and so $\gamma_r(G) \le |S| = |S'| + 2 = \gamma_r(G') + 2$, implying that $\w(G) < \w(G') + 20$. However, $\w(G) = 27 + (\w(G') - 1) = \w(G') + 26$, a contradiction.
Hence, $v_5 \in \cS$. Let $v_6$ be the neighbor of $v_5$ different from $v_4$. By Claim~\ref{claim.3}, $v_6 \in \cL$. In this case, let $G' = G - \{u,v,v_1,v_2,v_3,v_4,v_5\}$. By Claim~\ref{claim.bridge2}, $G' \notin \cB_{\rdom}$. Every $\gamma_r$-set of $G'$ can be extended to a RD-set of $G$ by adding to it $\{v,v_1,v_5\}$, and so $\gamma_r(G) \le \gamma_r(G') + 3$, implying that $\w(G) < \w(G') + 30$. However, $\w(G) = 32 + (\w(G') - 1) = \w(G') + 31$, a contradiction.~\smallqed

\medskip
By Claim~\ref{claim.18.1}, the vertices $v_1$ and $v_2$ are the only two common neighbors of $u$ and $v$. Let $x$ be the neighbor of $u$ different from $v_1$ and $v_2$, and let $y$ be the neighbor of $v$ different from $v_1$ and $v_2$.

\begin{subclaim}
\label{claim.18.2}
The vertices $x$ and $y$ are not adjacent.
\end{subclaim}
\proof Suppose that $x$ and $y$ are adjacent. If $x \in \cS$ and $y \in \cS$, then $G = R_2$, a contradiction. Hence at least one of $x$ and $y$ are large. Renaming vertices if necessary, assume that $y \in \cL$. Let $y_1$ be neighbor of $y$ different from $x$ and $v$. If $x \in \cS$, then the edge $yy_1$ is a bridge whose removal creates an $R_2$-component, contradicting Claim~\ref{claim.bridge2}. Hence, $x \in \cL$. Let $x_1$ be neighbor of $x$ different from $u$ and $y$.

Suppose that $x_1 \ne y_1$. In this case, let $G'$ be obtained from $G' - \{u,v,v_1,v_2,y\}$ by adding the edge $xy_1$. Let $S'$ be a $\gamma_r$-set of $G'$. If $x \in S'$, let $S = S' \cup \{v,y\}$. If $x \notin S'$ and $y_1 \in S$, let $S = S' \cup \{u,v_1\}$. If $x \notin S'$ and $y_1 \notin S$, let $S = S' \cup \{v,v_1\}$. In all cases, $S$ is a RD-set of $G$, and so $\gamma_r(G) \le |S| = |S'| + 2 = \gamma_r(G') + 2$, implying that $\w(G) < \w(G') + 20$. If $G' \notin \cB_{\rdom}$, then $\w(G) = \w(G') + 21$, a contradiction. Hence, $G' \in \cB_{\rdom}$. Let $G^* = G - v_2$, that is, $G^*$ is obtained from $G'$ by subdividing the edge $xy_1$ four times resulting in the path $xuv_1vyy_1$. By Observation~\ref{obser-5}, there exists a RD-set $S^*$ of $G^*$ such that $S^* \cap \{u,v_1,v,y\} = \{u,y\}$ and $|S^*| \le \gamma_{r}(G') + 1$. The set $S^*$ is a RD-set of $G$, and so $\gamma_r(G) \le |S^*| \le \gamma_{r}(G') + 1$, implying that $\w(G) < \w(G') + 10$. However, $\w(G) \le 22 + (\w(G') - 1 - 4) = \w(G') + 17$, a contradiction.

Hence, $x_1 = y_1$. In this case, we let $z = x_1$. Since no small vertex belongs to a triangle, we note that $z \in \cL$. Let $z_1$ be the neighbor of $z$ different from $x$ and $y$. If $z_1 \in \cL$, then we let $G' = G - \{u,v,v_1,v_2,x,y,z\}$. By Claim~\ref{claim.bridge2}, we note that $G' \notin \cB_{\rdom}$. Since $\gamma_r(G) \le \gamma_r(G') + 2$, we have $\w(G) < \w(G') + 20$. However, $\w(G) = 30 + (\w(G') - 1) = \w(G') + 29$, a contradiction. Hence, $z_1 \in \cS$. Let $z_2$ be the neighbor of $z_1$ different from $z$. Since every vertex of degree~$2$ belongs to a $1$-linkage, we note that $z_2 \in \cL$. We now let $G' = G - \{u,v,v_1,v_2,x,y,z,z_1\}$. By Claim~\ref{claim.bridge2}, we note that $G' \notin \cB_{\rdom}$. Since $\gamma_r(G) \le \gamma_r(G') + 3$, we have $\w(G) < \w(G') + 30$. However, $\w(G) = 35 + (\w(G') - 1) = \w(G') + 34$, a contradiction.~\smallqed

\medskip
By Claim~\ref{claim.18.2}, the vertices $x$ and $y$ are not adjacent.

\begin{subclaim}
\label{claim.18.3}
$x \in \cL$ and $y \in \cL$.
\end{subclaim}
\proof Suppose that at least one of $x$ and $y$ is small. Renaming vertices if necessary, we may assume that $x \in \cS$. Let $z$ be the neighbor of $x$ different from $u$. Necessarily, $z \in \cL$. Suppose that $yz \notin E(G)$. In this case, let $G'$ be the connected subcubic graph obtained from $G - \{u,v,v_1,v_2,x\}$ by adding the edge $yz$. Let $S'$ be a $\gamma_r$-set of $G'$. If at least one of $z$ and $y$ belongs to $S'$, let $S = S' \cup \{v,x\}$. If $z \notin S'$ and $y \notin S'$, let $S = S' \cup \{u,v_1\}$. In both cases, $S$ is a RD-set of $G$, and so $\gamma_r(G) \le |S| \le |S'| + 2 = \gamma_r(G') + 2$. Thus, $\w(G) < \w(G') + 20$.  If $G' \notin \{R_1,R_4,R_5\}$, then $\w(G) \ge \w(G') + 20$, a contradiction. Hence, $G' \in \{R_1,R_4,R_5\}$. Since every vertex of degree~$2$ in $G$ belongs to a $1$-linkage, $G' \notin \{R_1,R_5\}$, and so $G' = R_4$. Let $G^* = G - v_2$, and so $G^*$ is obtained from $G'$ by subdividing the added edge $yz$ four times resulting in the path $z x u v_1 v y$. By Observation~\ref{obser-5}, there exists a RD-set $S^*$ of $G^*$ such that $S^* \cap \{x,u,v_1,v\} = \{x,v\}$ and $|S^*| \le \gamma_{r}(G)$. The set $S^*$ is a RD-set of $G$, and so $\gamma_r(G) \le |S^*| \le \gamma_r(G')$, implying that $\w(G) < 10\gamma_r(G') \le \w(G')$. However, $w(G) = \w(G') + 19$, a contradiction. Hence, $yz \in E(G)$. Recall that $x \in \cS$ and $z \in \cL$.

Suppose that $y \in \cL$. In this case, let $G' = G - \{u,v,v_1,v_2,x\}$. We note that $\gamma_r(G) \le \gamma_r(G') + 2$, implying that $\w(G) < \w(G') + 20$. If $G' \notin \cB_{\rdom}$, then $\w(G) \ge \w(G') + 21$, a contradiction. Hence, $G' \in \cB_{\rdom}$. A type-$1$ NeRD-set of $G'$ with respect to the vertex~$y$ can be extended to a RD-set of $G$ by adding to it the set $\{v,x\}$. Therefore by Observation~\ref{obser-1}, $\gamma_r(G) \le \gamma_{r,\ndom}(G';y) + 2 \le (\gamma_{r}(G)-1) + 2 = \gamma_{r}(G) + 1$, implying that $\w(G) < \w(G') + 10$. However, $\w(G) \ge \w(G') + 17$, a contradiction. Hence, $y \in \cS$. Let $z_1$ be the neighbor of $z$ different from $x$ and $y$.

If $z_1 \in \cL$, then let $G' = G - \{u,v,v_1,v_2,x,y,z\}$. By Claim~\ref{claim.bridge2}, $G' \notin \cB_{\rdom}$. We note that $\gamma_r(G) \le \gamma_r(G') + 3$, implying that $\w(G) < \w(G') + 30$. However, $\w(G) = \w(G') + 31$, a contradiction. Hence, $z_1 \in \cS$. Let $z_2$ be the neighbor of $z_1$ different from~$z$. Necessarily, $z_2 \in \cL$. We now let $G' = G - \{u,v,v_1,v_2,x,y,z,z_1\}$. By Claim~\ref{claim.bridge2}, $G' \notin \cB_{\rdom}$. Every $\gamma_r$-set of $G'$ can be extended to a RD-set of $G$ by adding to it the set $\{v,x,z_1\}$, and so $\gamma_r(G) \le \gamma_r(G') + 3$, implying that $\w(G) < \w(G') + 30$. However, $\w(G) = \w(G') + 36$, a contradiction.~\smallqed

\medskip
By Claim~\ref{claim.18.3}, $x \in \cL$ and $y \in \cL$. Recall that $x$ and $y$ are not adjacent. Let $x_1$ and $x_2$ be the two neighbors of $x$ different from~$u$.

\begin{subclaim}
\label{claim.18.4}
$x_1x_2 \in E(G)$.
\end{subclaim}
\proof Suppose that $x_1x_2 \notin E(G)$. Let $G'$ be the subcubic graph obtained from $G - \{u,v,v_1,v_2,x\}$ by adding the edge $x_1x_2$. Let $G_x$ be the component of $G'$ containing the added edge $x_1x_2$. If $G'$ is disconnected, then let $G_y$ be the second component of $G'$ which necessarily contains the vertex~$y$. In this case the edge $vy$ is a bridge in $G$, implying by Claim~\ref{claim.bridge2} that $G_y \notin \cB_{\rdom}$. Therefore, the component $G_x$ is the only possible component of $G'$ that belongs to~$\cB_{\rdom}$.
Let $S'$ be a $\gamma_r$-set of $G'$. If at least one of $x_1$ and $x_2$ belongs to $S'$, let $S = S' \cup \{v,x\}$. If $x_1 \notin S'$ and $x_2 \notin S'$, let $S = S' \cup \{u,v_1\}$. In both cases, $S$ is a RD-set of $G$, and so $\gamma_r(G) \le |S| \le |S'| + 2 = \gamma_r(G') + 2$, implying that $\w(G) < \w(G') + 20$.  If $G_x \notin \cB_{\rdom}$, then $\w(G) \ge \w(G') + 21$, a contradiction. Hence, $G_x \in \cB_{\rdom}$. Our earlier properties of the graph $G$ imply that $G_x \ne R_1$.
Let $G^*$ be obtained from $G_x$ by subdividing the added edge $x_1x_2$ resulting in the path $x_1 x x_2$. Applying Observation~\ref{obser-2}, there exists a $\gamma_r$-set $S^*$ of $G^*$ such that $x \in S^*$ and $|S^*| = \gamma_r(G^*) \le \gamma_r(G_x)$. If $G'$ is connected, then let $S = S^* \cup \{v\}$. In this case, $|S| \le |S^*| + 1 \le \gamma_r(G_x) + 1 = \gamma_r(G') + 1$. If $G'$ is disconnected, let $S = S^* \cup S_y \cup \{v\}$, where $S_y$ is a $\gamma_r$-set of $G_y$. In this case, $|S| \le |S^*| + |S_y| + 1 \le \gamma_r(G_x) + \gamma_r(G_y) + 1 = \gamma_r(G') + 1$. In both cases, $S$ is a RD-set of $G$ and $|S| \le \gamma_r(G') + 1$, implying that $\w(G) < \w(G') + 10$. However, $\w(G) \ge \w(G') - 17$, a contradiction.~\smallqed

\medskip
We now return to the proof of Claim~\ref{claim.18}. By Claim~\ref{claim.18.4}, $x_1x_2 \in E(G)$. Since no vertex of degree~$2$ belongs to a triangle in $G$, we note that $x_1, x_2 \in \cL$. Recall that $y \in \cL$. If $y$ is adjacent to both $x_1$ and $x_2$, then the graph $G$ is determined and $\gamma_r(G) = 2$ and $\w(G) = 34$, a contradiction. Hence renaming vertices if necessary, we may assume that $x_1y \notin E(G)$. Let $G'$ be the connected subcubic graph obtained from $G - \{u,v,v_1,v_2,x\}$ by adding the edge $x_1y$. Let $S'$ be a $\gamma_r$-set of $G'$. If at least one of $x_1$ and $y$ belongs to $S'$, let $S = S' \cup \{v,x\}$. If $x_1 \notin S'$ and $y \notin S'$, let $S = S' \cup \{u,v_1\}$. In both cases, $S$ is a RD-set of $G$, and so $\gamma_r(G) \le |S| \le |S'| + 2 = \gamma_r(G') + 2$, implying that $\w(G) < \w(G') + 20$. If $G' \notin \cB_{\rdom}$, then $\w(G) \ge \w(G') + 21$, a contradiction. Hence, $G' \in \cB_{\rdom}$. Let $G^* = G - v_2$, that is, $G^*$ is obtained from $G'$ by subdividing the edge $x_1y$ four times resulting in the path $x_1xuv_1vy$. By Observation~\ref{obser-5}, there exists a RD-set $S^*$ of $G^*$ such that $S^* \cap \{x,u,v_1,v\} = \{x,v\}$ and $|S^*| \le \gamma_{r}(G') + 1$. The set $S^*$ is a RD-set of $G$, and so $\gamma_r(G) \le |S^*| \le \gamma_{r}(G') + 1$, implying that $\w(G) < \w(G') + 10$. However, $\w(G) \ge \w(G') + 17$, a contradiction. This completes the proof of Claim~\ref{claim.18}.~\smallqed

\medskip
By Claim~\ref{claim.18}, there is no $4$-cycle in $G$ that contains two small (non-adjacent) vertices.

\begin{claim}
\label{claim.20}
No small vertex in $G$ belongs to a $4$-cycle.
\end{claim}
\proof Suppose, to the contrary, that there is a vertex $v \in \cS$ that belongs to a $4$-cycle $C \colon v v_1 v_2 v_3 v$. By our earlier observations, $v_i \in \cL$ for $i \in [3]$. Let $u_i$ be the neighbor of $v_i$ that does not belong to $C$ for $i \in [3]$.

\begin{subclaim}
\label{claim.20.1}
$u_1 \in \cL$ and $u_3 \in \cL$.
\end{subclaim}
\proof  Suppose that at least one of $u_1$ and $u_3$ is small. Renaming vertices if necessary, we may assume that $u_1 \in \cS$. Since no small vertex belongs to a triangle, $u_1v_2 \notin E(G)$. Since no $4$-cycle contains two small vertices, $u_1v_3 \notin E(G)$. Let $u$ be the neighbor of $u_1$ different from $v_1$. By our earlier observations, $u \in \cL$ and $u \notin \{v_2,v_3\}$. If $u$ is adjacent to both $v_2$ and $v_3$, then the graph $G$ is determined and $\gamma_r(G) = 2$ and $\w(G) = 26$, a contradiction. Hence, $u$ is not adjacent to at least one of $v_2$ and $v_3$.

Suppose that $uv_3 \notin E(G)$. In this case, let $G'$ be the connected subcubic graph obtained from $G - \{v,v_1,u_1\}$ by adding the edge $uv_3$. Let $S'$ be a $\gamma_r$-set of $G'$. If $u \in S'$, let $S = S' \cup \{v\}$. If $u \notin S'$ and $v_3 \in S'$, let $S = S' \cup \{u_1\}$. If $u \notin S'$ and $v_3 \notin S'$, let $S = S' \cup \{v_1\}$. In all cases, $S$ is a RD-set of $G$, and so $\gamma_r(G) \le |S| \le |S'| + 1 = \gamma_r(G') + 1$, implying that $\w(G) < \w(G') + 10$. If $G' \notin \{R_1,R_4,R_5\}$, then $\w(G) \ge \w(G') + 10$, a contradiction. Hence, $G' \in \{R_1,R_4,R_5\}$. We note that $u$ and $v_3$ are adjacent vertices of degree~$3$ in $G'$, and so $G' \ne R_1$. Since there is no $2$-linkage in $G$, we note that $G' \ne R_5$. Hence, $G' = R_4$. The graph $G$ is now determined and satisfies $\gamma_r(G) = 4$ and $\w(G) = 49$, a contradiction.

Hence, $uv_3 \in E(G)$, that is, $u = u_3$. We now let $G' = G - \{v,v_1,u_1\}$. The graph $G'$ is a connected subcubic graph. Let $S'$ be a $\gamma_r$-set of $G'$. If $u \in S'$, let $S = S' \cup \{v\}$. If $u \notin S'$ and $v_3 \in S'$, let $S = S' \cup \{u_1\}$. If $u \notin S'$ and $v_3 \notin S'$, let $S = S' \cup \{v_1\}$. In all cases, $S$ is a RD-set of $G$, and so $\gamma_r(G) \le |S| \le |S'| + 1 = \gamma_r(G') + 1$, implying that $\w(G) < \w(G') + 10$. If $G' \notin \cB_{\rdom}$ or if $G' \in \cB_{\rdom,1}$, then $\w(G) \ge \w(G') + 10$, a contradiction. Hence, $G' \in  \cB_{\rdom,i}$ for some $i \in \{2,3,4,5\}$. Since $uv_3v_2$ is a path in $G'$, and $u$, $v_3$ and $v_2$ all have degree~$2$ in $G'$, either $G' = R_1$ or $G' = R_3$. If $G' = R_1$, then $G$ would contain a $2$-linkage, and if $G' = R_3$, then $G$ would contain a $3$-linkage. Both cases produce a contradiction.~\smallqed

\medskip
By Claim~\ref{claim.20.1}, $u_1 \in \cL$ and $u_3 \in \cL$. If $u_1 = u_2 = u_3$, then the graph $G$ is determined and $\gamma_r(G) = 2$ and $\w(G) = 21$, a contradiction. Renaming vertices if necessary, we may assume that $u_2 \ne u_3$. In this case, let $G'$ be the subcubic graph obtained from $G - \{v,v_1,v_2\}$ by adding the edge $u_2v_3$. Let $G_1$ be the component of $G'$ containing the vertex $u_1$ and let $G_2$ be the component of $G'$ containing the added edge $u_2v_3$. If $G'$ is connected, then $G' = G_1 = G_2$. If disconnected, then the edge $u_1v_1$ is a bridge in $G$, implying by Claim~\ref{claim.bridge2} that $G_1 \notin \cB_{\rdom}$. Therefore, the component $G_2$ is the only possible component of $G'$ that belongs to~$\cB_{\rdom}$.

Let $S'$ be a $\gamma_r$-set of $G'$. If $u_2 \in S'$, let $S = S' \cup \{v\}$. If $u_2 \notin S'$ and $v_3 \in S'$, let $S = S' \cup \{v_2\}$. If $u_2 \notin S'$ and $v_3 \notin S'$, let $S = S' \cup \{v_1\}$. In all cases, $S$ is a RD-set of $G$, and so $\gamma_r(G) \le |S| \le |S'| + 1 = \gamma_r(G') + 1$, implying that $\w(G) < \w(G') + 10$. By our earlier properties of the graph $G$, we note that $G_2 \ne R_1$. If $G_2 \notin \{R_4,R_5,R_{10}\}$, then $\w(G) \ge \w(G) + 10$, a contradiction. Hence, $G_2 \in \{R_4,R_5,R_{10}\}$. If $G'$ is connected, then the graph $G$ is determined and either $G_2 \in \{R_4,R_5\}$, in which case $\gamma_r(G) \le 4$ and $\w(G) = 47$, or $G_2 = R_{10}$, in which case $\gamma_r(G) \le 5$ and $\w(G) = 58$. In both cases, we have a contradiction. Hence, $G'$ is disconnected.

Since every small vertex in $G$ belongs to a $1$-linkage, the case $G_2 = R_5$ cannot occur, and so $G_2 \in \{R_4,R_{10}\}$. Let $G_v$ be the component of $G - v_1u_1$ that contains the vertex~$v$. Thus, $G_v$ is obtained from $G_2$ by subdividing the added edge $u_2v_3$ three times resulting in the path $u_2 v_2 v_1 v v_3$ and adding the edge $v_2v_3$. If $G_2 = R_4$, then $\gamma_r(G_v) \le 4$, and so $\gamma_r(G) \le \gamma_r(G_1) + \gamma_r(G_v) \le \gamma_r(G_1) + 4$, implying that $\w(G) < \w(G_1) + 40$. However in this case, $\w(G) = \w(G_1) + 47$, a contradiction. If $G_2 = R_{10}$, then $\gamma_r(G_v) \le 5$, and so $\gamma_r(G) \le \gamma_r(G_1) + \gamma_r(G_v) \le \gamma_r(G_1) + 5$, implying that $\w(G) < \w(G_1) + 50$. However in this case, $\w(G) = \w(G_1) + 58$, a contradiction. This completes the proof of Claim~\ref{claim.20}.~\smallqed

\medskip
Recall that no small vertex belongs to a triangle. By Claim~\ref{claim.20}, no small vertex in $G$ belongs to a $4$-cycle. Hence every cycle that contains a small vertex in $G$ has at least five vertices.

\begin{claim}
\label{claim.21}
No large vertex has two small neighbors and one large neighbor.
\end{claim}
\proof Suppose, to the contrary, that $v \in \cL$ and $N(v) = \{v_1,v_2,v_3\}$ where $v_1,v_2 \in \cS$ and $v_3 \in \cL$. Let $u_1$ and $u_2$ be the neighbors of $v_1$ and $v_2$, respectively, different from $v$. Since no vertex of degree~$2$ belongs to a triangle or a $4$-cycle, $\{u_1,u_2\} \cap N(v) = \emptyset$. Further, $u_1v_1$ and $u_2v_2$ are the only edges between $\{u_1,u_2\}$ and $N(v)$. Let $u_3$ and $w_3$ be the neighbors of $v_3$ different from $v_3$. The graph illustrated in Figure~\ref{rdom:fig-31} is a subgraph of~$G$.

\begin{figure}[htb]
\begin{center}
\begin{tikzpicture}[scale=.8,style=thick,x=0.8cm,y=0.8cm]
\def\vr{2.5pt} 
\path (0,0) coordinate (u1);
\path (-0.5,-0.5) coordinate (u11);
\path (0.5,-0.5) coordinate (u12);
\path (0,1) coordinate (u2);
\path (2,0) coordinate (u3);
\path (1.5,-0.5) coordinate (u31);
\path (2.5,-0.5) coordinate (u32);
\path (2,1) coordinate (u4);
\path (2,2) coordinate (u5);
\path (3.25,0) coordinate (u6);
\path (3.3,0) coordinate (u6p);
\path (2.75,-0.5) coordinate (u61);
\path (3.75,-0.5) coordinate (u62);
\path (4,1) coordinate (u7);
\path (4.75,0) coordinate (u8);
\path (4.8,0) coordinate (u8p);
\path (4.25,-0.5) coordinate (u81);
\path (5.25,-0.5) coordinate (u82);
\draw (u1)--(u2)--(u5)--(u4)--(u3);
\draw (u5)--(u7)--(u6);
\draw (u7)--(u8);
\draw (u11)--(u1)--(u12);
\draw (u31)--(u3)--(u32);
\draw (u61)--(u6)--(u62);
\draw (u81)--(u8)--(u82);
\draw (u1) [fill=white] circle (\vr);
\draw (u2) [fill=white] circle (\vr);
\draw (u3) [fill=white] circle (\vr);
\draw (u4) [fill=white] circle (\vr);
\draw (u5) [fill=white] circle (\vr);
\draw (u6) [fill=white] circle (\vr);
\draw (u7) [fill=white] circle (\vr);
\draw (u8) [fill=white] circle (\vr);
\draw[anchor = east] (u1) node {{\small $u_1$}};
\draw[anchor = east] (u2) node {{\small $v_1$}};
\draw[anchor = east] (u3) node {{\small $u_2$}};
\draw[anchor = east] (u4) node {{\small $v_2$}};
\draw[anchor = south] (u5) node {{\small $v$}};
\draw[anchor = west] (u6p) node {{\small $u_3$}};
\draw[anchor = west] (u8p) node {{\small $w_3$}};
\draw[anchor = west] (u7) node {{\small $v_3$}};
\end{tikzpicture}
\end{center}
\begin{center}
\vskip -0.5 cm
\caption{A subgraph in the proof of Claim~\ref{claim.21}}
\label{rdom:fig-31}
\end{center}
\end{figure}
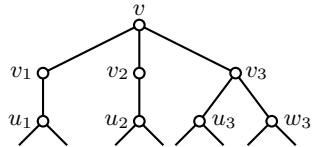

\vskip -0.5 cm
Let $G = G - \{v,v_1,v_2\}$. The graph $G'$ is a subcubic graph with at most three components.
Let $S'$ be a $\gamma_r$-set of $G'$. If $u_1 \in S'$, let $S = S' \cup \{v_2\}$. If $u_1 \notin S'$ and $u_2 \in S'$, let $S = S' \cup \{v_1\}$. If $u_1 \notin S'$ and $u_2 \notin S'$, let $S = S' \cup \{v\}$. In all cases, $S$ is a RD-set of $G$, and so $\gamma_r(G) \le |S| \le |S'| + 1 = \gamma_r(G') + 1$, implying that $\w(G) < \w(G') + 10$. If no component belongs to $\cB_{\rdom}$, then $\w(G) \ge \w(G') + 11$, a contradiction. Hence at least one component of $G'$ belongs to~$\cB_{\rdom}$. Let $H$ be such a component of $G'$. Possibly, $G' = H$. Since the removal of a bridge cannot create a component that belongs to $\cB_{\rdom}$, the component $H$ necessarily contains at least two vertices from the set $\{u_1,u_2,v_3\}$. We note that each of $u_1$, $u_2$ and $v_3$ has degree~$2$ in $G'$. Thus, at least one of $u_1$ and $u_2$ belong to the component $H$.

Suppose that exactly one of $u_1$ and $u_2$ belong to $H$. Let $G^*$ be obtained from $G'$ by adding the edge $u_1u_2$. The resulting graph $G^*$ is a connected subcubic graph that contains a bridge, namely the added edge $u_1u_2$. Since no graph in $\cB_{\rdom}$ contains a bridge, $G^* \notin \cB_{\rdom}$. Let $S^*$ be a $\gamma_r$-set of $G^*$. If $u_1 \in S^*$, let $S = S^* \cup \{v_2\}$. If $u_1 \notin S^*$ and $u_2 \in S^*$, let $S = S^* \cup \{v_1\}$. If $u_1 \notin S^*$ and $u_2 \notin S^*$, let $S = S^* \cup \{v\}$. In all cases, $S$ is a RD-set of $G$, and so $\gamma_r(G) \le |S| \le |S^*| + 1 = \gamma_r(G^*) + 1$, implying that $\w(G) < \w(G^*) + 10$. However since $G^* \notin \cB_{\rdom}$, we have $\w(G) = \w(G^*) + 13$, a contradiction. Hence, $\{u_1,u_2\} \subset V(H)$.

Let $X = \{u_1,u_2\}$, and so $X \subset V(H)$. As observed earlier, $u_1$ and $u_2$ have degree~$2$ in $G'$. Let $S_H$ be a minimum type-$2$ NeRD-set in $H$ with respect to the set $X$. By Observation~\ref{obser-1}(f), we have $|S_H| =  \gamma_{r,\dom}(H;X) \le \gamma_{r}(H) - 1$. Suppose that $G' = H$. In this case, the set $S_H^* = S_H \cup \{v\}$ is a RD-set of $G$, and so $\gamma_r(G) \le |S_H^*| + 1 \le (\gamma_{r}(H) - 1) + 1 = \gamma_r(G')$, implying that $\w(G) < \w(G')$. However, $\w(G) \ge \w(G') + 7$, a contradiction. Hence, $G' \ne H$. Let $G_3$ be the component of $G'$ containing the vertex~$v_3$, and so $G' = H \cup G_3$. We note that the removal of the bridge $vv_3$ creates the component $G_3$, implying that $G_3 \notin \cB_{\rdom}$. Let $S_3$ be a $\gamma_r$-set of $G_3$. In this case, the set $S_H^* = S_H \cup S_3 \cup \{v\}$ is a RD-set of $G$, and so $\gamma_r(G) \le |S_H^*| + |S_3| + 1 \le (\gamma_{r}(H) - 1) + \gamma_r(G_3) + 1 = \gamma_r(G')$. Therefore,  $\w(G) < \w(G')$. However, $\w(G) \ge 14 + (\w(H) - 2 - 4) + (\w(G_3) - 1) = (\w(H) + \w(G_3)) + 7 = \w(G') + 7$, a contradiction. This completes the proof of Claim~\ref{claim.21}.~\smallqed

\medskip
By Claim~\ref{claim.21}, no large vertex has two small neighbors and one large neighbor. Hence if a large vertex has a small neighbor, then it has either one small neighbor or three small neighbors.

\begin{claim}
\label{claim.22}
No large vertex has exactly one small neighbor.
\end{claim}
\proof Suppose, to the contrary, that $v \in \cL$ has exactly one small neighbor. Let $N(v) = \{v_1,v_2,v_3\}$, where $v_1 \in \cS$ and $v_2,v_3 \in \cL$. Let $u$ be the neighbor of $v_1$ different from~$v$. Necessarily, $u \in \cL$. Let $u_1$ and $u_2$ be the two neighbors of $u$ different from~$v_1$. Since no vertex of degree~$2$ belongs to a $3$-cycle or $4$-cycle, the vertices $u_1,u_2,v_2,v_3$ are pairwise distinct.

\begin{subclaim}
\label{claim.22.1}
$\{u_1,u_2\} \subset \cS$.
\end{subclaim}
\proof Suppose, to the contrary, that at least one of $u_1$ and $u_2$ is large. By Claim~\ref{claim.21} this implies that both $u_1$ and $u_2$ are large. Let $G' = G - \{u,v,v_1\}$. Suppose that $G'$ contains a component $F$ such that $F \in \cB_{\rdom}$. Since the removal of a bridge cannot create a component that belongs to $\cB_{\rdom}$, the component $F$ contains at least two vertices from the set $\{u_1,u_2,v_2,v_3\}$.

Suppose that $F$ contains a vertex from both $\{u_1,u_2\}$ and $\{v_2,v_3\}$. By symmetry, and renaming vertices if necessary, we may assume that $\{u_2,v_3\} \subset V(F)$. We note that both $u_2$ and $v_3$ have degree~$2$ in $F$. Applying Observation~\ref{obser-1}(f) to the graph $F$ with $X = \{u_2,v_3\}$, we have $\gamma_{r,\dom}(F;X) \le \gamma_{r}(F) - 1$. Let $S_F$ be a minimum type-$2$ NeRD-set of $F$ with respect to the set $X$, and so $|S_F| = \gamma_{r,\dom}(F;X) \le \gamma_{r}(F) - 1$.

Suppose that $F$ is the only component of $G'$ that belongs to $\cB_{\rdom}$. If $G'$ is connected, then the set $S_F \cup \{v_1\}$ is a RD-set of $G$. If $G'$ is disconnected, then the set $S_F \cup \{v_1\}$ can be extended to a RD-set by adding to it a $\gamma_r$-set from the component(s) of $G'$ different from~$F$. In both cases, we infer that $\gamma_r(G) \le 1 + (\gamma_r(G') - 1) = \gamma_r(G')$, implying that $\w(G) < \w(G')$. Since $G'$ has exactly one component that belongs to $\cB_{\rdom}$, we have $\w(G) \ge \w(G') + 4$, a contradiction.
Hence, the graph $G'$ contains a component $H$, different from $F$, that belongs to $\cB_{\rdom}$. In this case, $\{u_1,v_2\} \subset V(H)$ and, analogously as with the component $F$, there exists a minimum type-$2$ NeRD-set of $H$ with respect to the set $\{u_1,v_2\}$ satisfying $|S_H| \le \gamma_{r}(H) - 1$. The set $S_F \cup S_H \cup \{v_1\}$ is a RD-set of $G$, and so $\gamma_r(G) \le 1 + |S_F| + |S_H| \le 1 + (\gamma_r(F) - 1) + (\gamma_r(H) - 1) = \gamma_r(G') - 1$, noting that $G' = F \cup H$. This implies that $\w(G) \le \w(G') - 10$. However, $\w(G) \ge 13 + (\w(G') - 4 - 5 - 5) = \w(G') - 1$, a contradiction.
Hence, if the graph $G'$ contains a component $C$ in $\cB_{\rdom}$, then either $\{u_1,u_2\} \subseteq V(C)$ and $\{v_2,v_3\} \cap V(C) = \emptyset$ or $\{v_2,v_3\} \subseteq V(C)$ and $\{u_1,u_2\} \cap V(C) = \emptyset$. If all edges are present between $\{u_1,u_2\}$ and $\{v_2,v_3\}$, then $G = R_{11}$, a contradiction. Hence renaming vertices if necessary, we may assume that $u_1v_2 \notin E(G)$.

Let $G''$ be the graph obtained from $G'$ by adding the edge $u_1v_2$. The resulting graph $G'$ is a subcubic graph with at most three components. Let $G_1$ be the component of $G'$ containing the added edge $u_1v_2$, and let $G_2$ and $G_3$ be the components of $G'$ containing $v_3$ and $u_2$, respectively. If $G'$ is connected, then $G_1 = G_2 = G_3$. Let $S'$ be a $\gamma_r$-set of $G'$. If $u_1 \in S'$, let $S = S' \cup \{v\}$. If $u_1 \notin S'$ and $v_2 \in S'$, let $S = S' \cup \{u\}$. If $u_1 \notin S'$ and $v_2 \notin S'$, let $S = S' \cup \{v_1\}$. In all cases, $S$ is a RD-set of $G$, and so $\gamma_r(G) \le |S| \le |S'| + 1 = \gamma_r(G') + 1$, implying that $\w(G) < \w(G') + 10$. If no component of $G'$ belongs to $\cB_{\rdom}$, then $\w(G) \ge \w(G') + 11$, a contradiction. Hence at least one component of $G'$ belongs to~$\cB_{\rdom}$. Let $H$ be such a component of $G$. If $H \ne G_1$, then since the removal of a bridge in $G$ cannot create a component in $\cB_{\rdom}$, this implies that $\{u_2,v_3\} \subset V(H)$ and $\{u_1,v_2\} \cap V(H) = \emptyset$. However, such a component $H$ is a component in $G'$, contradicting our earlier properties of a component of $G'$ that belongs to~$\cB_{\rdom}$. Hence, $H = G_1$ and $H$ is the only component of $G'$ that belongs to~$\cB_{\rdom}$.

Suppose that $G'$ is connected, and so $G' = G_1 \in \cB_{\rdom}$ and the graph $G$ is determined. We note that the vertices $u_1$ and $v_2$ are adjacent vertices of degree~$3$ in $G'$, and so $G' \notin \{R_1,R_2\}$. Further, we note that $u_2$ and $v_3$ have degree~$2$ in $G'$. Reconstructing the graph $G$ from $G' \in \cB_{\rdom}$ it can be readily checked that $10\gamma_r(G) \le \w(G)$, a contradiction. 
Hence, $G'$ is disconnected, and so $G_1 \ne G_2$ or $G_1 \ne G_3$. By symmetry and renaming vertices if necessary, we may assume that $G_1 \ne G_3$. Let $G_u$ be obtained from $G_1$ by subdividing the added edge $u_1v_2$ of $G_1$ three times resulting path in the path $u_1 u v_1 v v_2$. Let $S_u^1$ be a minimum type-$1$ NeRD-set of $G_u$ with respect to the vertex~$u$, and let $S_v^2$ be a minimum type-$2$ NeRD-set of $G_u$ with respect to the vertex~$u$. By Observation~\ref{obser-4}, $|S_u^1| = \gamma_{r,\ndom}(G_u;u) \le \gamma_r(G_1)$ and $|S_u^1| = \gamma_{r,\dom}(G_u;u) \le \gamma_r(G_1)$. Recall that $G_1 \ne G_3$. Let $S_3$ be a $\gamma_r$-set of $G_3$. If $u_2 \in S_3$, then let $S = S_u^1 \cup S_3$, while if $u_2 \notin S_3$, then let $S = S_u^2 \cup S_3$. In both cases, $|S| \le \gamma_r(G_1) + \gamma_r(G_3)$.

If $G_2 = G_1$ or $G_2 = G_3$, then $\gamma(G') = \gamma_r(G_1) + \gamma_r(G_3)$ and $S$ is a RD-set of $G$. If $G_2 \ne G_1$ or $G_2 \ne G_3$, then $\gamma(G') = \gamma_r(G_1) + \gamma_r(G_2) + \gamma_r(G_3)$ and $S$ can be extended to a RD-set of $G$ by adding to it a $\gamma_r$-set of $G_2$. In both cases, we have that $\gamma_r(G) \le |S| \le \gamma_r(G')$, and we infer that $\w(G) < \w(G')$. As observed earlier, $G_1$ is the only component of $G'$ that belongs to~$\cB_{\rdom}$. Hence, $\w(G) \ge 13 + (\w(G') - 2 - 4) = \w(G') + 7$, a contradiction.~\smallqed

\medskip
By Claim~\ref{claim.22.1}, $u_1 \in \cS$ and $u_2 \in \cS$.

\begin{subclaim}
\label{claim.22.2}
There is no edge between $\{u_1,u_2\}$ and $\{v_2,v_3\}$.
\end{subclaim}
\proof Suppose that there is an edge between $\{u_1,u_2\}$ and $\{v_2,v_3\}$. Renaming vertices if necessary, we may assume that $u_1v_2 \in E(G)$. Since no small vertex belongs to a $4$-cycle, we note that $u_2v_2 \notin E(G)$.
Suppose that $u_2v_3 \in E(G)$. If $v_2v_3 \in E(G)$, then the graph $G$ is determined and $\gamma_r(G) = 3$ and $\w(G) = 31$, a contradiction. Hence, $v_2v_3 \notin E(G)$. In this case, let $G'$ be the connected subcubic graph obtained from $G - \{u,v,v_1,u_1,u_2\}$ by adding the edge $v_2v_3$. Let $S'$ be a $\gamma_r$-set of $G'$. If $v_2 \in S'$, let $S = S' \cup \{u_2,v\}$. If $v_2 \notin S'$ and $v_3 \in S'$, let $S = S' \cup \{v,u_1\}$. If $v_2 \notin S'$ and $v_3 \notin S'$, let $S = S' \cup \{u,v_1\}$. In all cases, $S$ is a RD-set of $G$, and so $\gamma_r(G) \le |S| \le |S'| + 2 = \gamma_r(G') + 2$, implying that $\w(G) < \w(G') + 20$. We note that $v_2$ and $v_3$ are adjacent vertices of degree~$2$ in $G'$. By our earlier properties of the graph $G$, we infer that $G' \notin \{R_1,R_4,R_5\}$. Let $G^*$ be obtained from $G'$ by subdividing the edge $v_2v_3$ four times resulting in the path $v_2 u_1 u v_1 v v_3$. By Observation~\ref{obser-6}(a), there exists a RD-set $S^*$ in $G^*$ such that $v_1 \in S^*$ and $|S^*| \le \gamma_r(G')$. The set $S^* \cup \{u_2\}$ is a RD-set of $G$, and so $\gamma_r(G) \le |S^*| + 1 \le \gamma_r(G') + 1$, implying that $\w(G) < \w(G') + 10$. However, $\w(G) \ge \w(G) + 18$, a contradiction. Hence, $u_2v_3 \notin E(G)$. Let $x$ be the neighbor of $u_2$ different from~$u$. By our earlier observations, $x \in \cL$ and $x \notin \{v_2,v_3\}$.

We show next that $xv_2 \in E(G)$. Suppose, to the contrary, that $xv_2 \notin E(G)$. In this case, let $G'$ be the subcubic graph obtained from $G - \{u,v,v_1,u_1,u_2\}$ by adding the edge $xv_2$. Let $G_x$ be the component containing~the added edge $xv_2$, and let $G_3$ be the component containing the vertex~$v_3$. If $G'$ is connected, then $G_x = G_3$. If $G'$ is disconnected, then it has two components, $G_x$ and $G_3$. In this case, since the removal of a bridge cannot create a component in $\cB_{\rdom}$, we note that $G_3  \notin \cB_{\rdom}$.

Let $S'$ be a $\gamma_r$-set of $G'$. If $v_2 \in S'$, let $S = S' \cup \{u_2,v\}$. If $v_2 \notin S'$ and $x \in S'$, let $S = S' \cup \{v,u_1\}$. If $v_2 \notin S'$ and $x \notin S'$, let $S = S' \cup \{u,v_1\}$. In all cases, $S$ is a RD-set of $G$, and so $\gamma_r(G) \le |S| \le |S'| + 2 = \gamma_r(G') + 2$, implying that $\w(G) < \w(G') + 20$. We note that $v_2$ and $x$ are adjacent vertices of degree~$2$ and degree~$3$, respectively, in $G_x$. In particular, $G_x \ne R_1$. If $G_x = R_5$, then $G' = G_x$, and the graph $G$ is determined and $\gamma_r(G) \le 5$ and $\w(G) = 57$, a contradiction. Hence, $G_x \ne R_5$.

Suppose that $G_x = R_4$. If $G' = G_x$, then the graph $G$ is determined and $\gamma_r(G) \le 5$ and $\w(G) = 57$, a contradiction. Hence, $G'$ is disconnected. In this case, let $G_v$ be the component of $G - vv_3$ that contains the vertex~$v$. We infer from the structure of the graph $G_v$ (using the structure of $G_x$) that a $\gamma_r$-set of $G_3$ can be extended to a RD-set of $G$ by adding to it five vertices from $G_v$, and so $\gamma_r(G) \le 5 + \gamma_r(G_3)$. This implies that $\w(G) < \w(G_3) + 50$. However, $\w(G) = \w(G_3) + 56$, a contradiction. Hence, $G_x \ne R_4$.

Suppose that $G_x = R_{10}$. If $G' = G_x$, then the graph $G$ is determined and $\gamma_r(G) \le 6$ and $\w(G) = 60$, a contradiction. Hence, $G'$ is disconnected. In this case, let $G_v$ be the component of $G - vv_3$ that contains the vertex~$v$. We infer from the structure of the graph $G_v$ (using the structure of $G_x$) that a $\gamma_r$-set of $G_3$ can be extended to a RD-set of $G$ by adding to it six vertices from $G_v$, and so $\gamma_r(G) \le 6 + \gamma_r(G_3)$. This implies that $\w(G) < \w(G_3) + 60$. However, $\w(G) = \w(G_3) + 60$, a contradiction. Hence, $G_x \ne R_{10}$.

Hence, $G_x \notin \{R_1,R_4,R_5,R_{10}\}$. Recall that if $G_3 \ne G_x$, then $G_3 \notin \cB_{\rdom}$. Thus there is at most one bad component in $G'$, and such a component does not belong to~$\{R_1,R_4,R_5,R_{10}\}$. Hence, $\w(G) \ge \w(G') + 20$, a contradiction. Hence, $xv_2 \in E(G)$. The graph $G$ therefore contains the subgraph shown in Figure~\ref{rdom:fig-22.2}(a). Let $C$ be the cycle $vv_1uu_1v_2v$, and let $G'$ be the connected special subcubic graph obtained from $G - V(C)$ by adding the edge $u_2v_3$. Let $S'$ be a $\gamma_r$-set of $G'$. If $u_2 \in S'$, let $S = S' \cup \{v,v_2\}$. If $u_2 \notin S'$ and $v_3 \in S'$, let $S = S' \cup \{u,u_1\}$. If $u_2 \notin S'$ and $v_3 \notin S'$, let $S = S' \cup \{v_1,v_2\}$. In all cases, $S$ is a RD-set of $G$, and so $\gamma_r(G) \le |S| \le |S'| + 2 = \gamma_r(G') + 2$, implying that $\w(G) < \w(G') + 20$. If $G' \notin \cB_{\rdom}$, then $\w(G) = \w(G') + 21$, a contradiction. Hence, $G' \in \cB_{\rdom}$.  We note that $P \colon xu_2v_3$ is a path in $G'$ where the vertices $x$, $u_2$ and $v_3$ have degrees~$2$,~$2$, and~$3$, respectively in $G'$. Our earlier properties of the graph $G$, together with the existence of the path $P$ in $G'$, imply that $G' = R_8$. Reconstructing the graph $G$ from $G'$ now yields the graph shown in Figure~\ref{rdom:fig-22.2}(b) that satisfies $\gamma_r(G) = 6$ and $\w(G) = 70$, a contradiction. (The six shaded vertices, for example, shown in Figure~\ref{rdom:fig-22.2}(b) form a $\gamma_r$-set in $G$.) This completes the proof of Claim~\ref{claim.22.2}.~\smallqed

\begin{figure}[htb]
\begin{center}
\begin{tikzpicture}[scale=.8,style=thick,x=0.8cm,y=0.8cm]
\def\vr{2.5pt} 
\path (1.5,1) coordinate (u1);
\path (1.25,0) coordinate (u2);
\path (1,2) coordinate (u3);
\path (1,2.1) coordinate (u3p);
\path (2.5,1) coordinate (u4);
\path (2.1,1) coordinate (u4p);
\path (2,3) coordinate (u5);
\path (3.5,0) coordinate (u6);
\path (3,2) coordinate (u7);
\path (4.5,0) coordinate (u8);
\path (5.35,0) coordinate (u10);
\path (4.75,3) coordinate (u9);
\path (5.5,3) coordinate (u11);
\path (4.75,2.35) coordinate (u12);
\draw (u1)--(u4)--(u7)--(u5)--(u3)--(u1)--(u2)--(u6)--(u7);
\draw (u5)--(u9);
\draw (u6)--(u8);
\draw (u12)--(u9)--(u11);
%
\draw (u1) [fill=white] circle (\vr);
\draw (u2) [fill=white] circle (\vr);
\draw (u3) [fill=white] circle (\vr);
\draw (u4) [fill=white] circle (\vr);
\draw (u5) [fill=white] circle (\vr);
\draw (u6) [fill=white] circle (\vr);
\draw (u7) [fill=white] circle (\vr);
\draw (u9) [fill=white] circle (\vr);
\draw[anchor = east] (u1) node {{\small $u$}};
\draw[anchor = north] (u2) node {{\small $u_2$}};
\draw[anchor = east] (u3p) node {{\small $v_1$}};
\draw[anchor = north] (u4) node {{\small $u_1$}};
\draw[anchor = south] (u5) node {{\small $v$}};
\draw[anchor = north] (u6) node {{\small $x$}};
\draw[anchor = west] (u7) node {{\small $v_2$}};
\draw[anchor = south] (u9) node {{\small $v_3$}};
%
\draw (2.5,-1) node {{\small (a)}};
%
\path (9.5,1) coordinate (u1);
\path (9.25,0) coordinate (u2);
\path (9,2) coordinate (u3);
\path (9,2.1) coordinate (u3p);
\path (10.5,1) coordinate (u4);
\path (10.1,1) coordinate (u4p);
\path (10,3) coordinate (u5);
\path (11.5,0) coordinate (u6);
\path (11,2) coordinate (u7);
\path (12.75,0) coordinate (u8);
\path (13.35,0) coordinate (u10);
\path (12.75,3) coordinate (u9);
\path (13.5,3) coordinate (u11);
\path (12.75,2.35) coordinate (u12);
\path (12.75,1.50) coordinate (p1);
\path (13.5,3) coordinate (p2);
\path (14.25,0) coordinate (p3);
\path (14.25,0.750) coordinate (p4);
\path (14.25,1.50) coordinate (p5);
\path (14.25,3) coordinate (p6);
\path (15.5,1.50) coordinate (p7);
\draw (u1)--(u4)--(u7)--(u5)--(u3)--(u1)--(u2)--(u6)--(u7);
\draw (u5)--(u9);
\draw (u6)--(u8);
\draw (u8)--(p1)--(u9)--(u11);
\draw (p2)--(p6)--(p5)--(p4)--(p3)--(u8);
\draw (p1)--(p5);
\draw (p3)--(p7)--(p6);
\draw (u1) [fill=black] circle (\vr);
\draw (u2) [fill=white] circle (\vr);
\draw (u3) [fill=black] circle (\vr);
\draw (u4) [fill=white] circle (\vr);
\draw (u5) [fill=black] circle (\vr);
\draw (u6) [fill=white] circle (\vr);
\draw (u7) [fill=white] circle (\vr);
\draw (u8) [fill=black] circle (\vr);
\draw (u9) [fill=white] circle (\vr);
\draw (p1) [fill=white] circle (\vr);
\draw (p2) [fill=white] circle (\vr);
\draw (p3) [fill=white] circle (\vr);
\draw (p4) [fill=white] circle (\vr);
\draw (p5) [fill=black] circle (\vr);
\draw (p6) [fill=black] circle (\vr);
\draw (p7) [fill=white] circle (\vr);
%
%
\draw (12,-1) node {{\small (b)}};
%
\end{tikzpicture}
\end{center}
\begin{center}
\vskip -0.5 cm
\caption{Subgraphs in the proof of Claim~\ref{claim.22.2}}
\label{rdom:fig-22.2}
\end{center}
\end{figure}
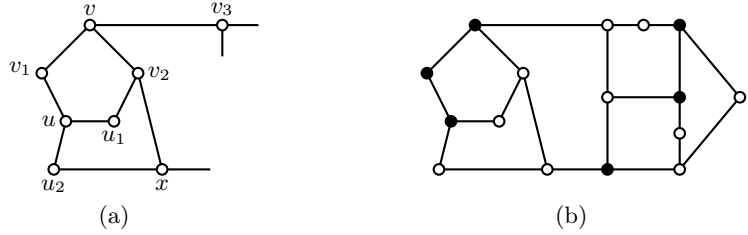

\vskip -0.5 cm
Let $x$ and $y$ be the neighbors of $u_1$ and $u_2$, respectively, different from~$u$. By Claim~\ref{claim.22.2}, there is no edge between $\{u_1,u_2\}$ and $\{v_2,v_3\}$, implying that $\{x,y\} \cap \{v_2,v_3\} = \emptyset$. Hence, the graph illustrated in Figure~\ref{rdom:fig-35} is a subgraph of $G$, where possibly edges between $\{x,y\}$ and $\{v_2,v_3\}$ may exist. By our earlier observations, $\{v_1,u_1,u_2\} \subseteq \cS$ and $\{u,v,v_2,v_3,x,y\} \subseteq \cL$.

\begin{figure}[htb]
\begin{center}
\begin{tikzpicture}[scale=.8,style=thick,x=0.8cm,y=0.8cm]
\def\vr{2.5pt} 
\path (0,0) coordinate (u1);
\path (-0.5,-0.5) coordinate (u11);
\path (-0.5,0.5) coordinate (u12);
\path (0,2) coordinate (u2);
\path (0,2.05) coordinate (u2p);
\path (-0.5,1.5) coordinate (u21);
\path (-0.5,2.5) coordinate (u22);
\path (1,0) coordinate (u3);
\path (1,2) coordinate (u4);
\path (2,1) coordinate (u5);
\path (3,1) coordinate (u6);
\path (4,1) coordinate (u7);
\path (5,0) coordinate (u8);
\path (5.5,-0.5) coordinate (u81);
\path (5.5,0.5) coordinate (u82);
\path (5,2) coordinate (u9);
\path (5.5,1.5) coordinate (u91);
\path (5.5,2.5) coordinate (u92);
\draw (u11)--(u1)--(u12);
\draw (u21)--(u2)--(u22);
\draw (u81)--(u8)--(u82);
\draw (u91)--(u9)--(u92);
\draw (u1)--(u3)--(u5)--(u4)--(u2);
\draw (u8)--(u7)--(u9);
\draw (u5)--(u6)--(u7);
\draw (u1) [fill=white] circle (\vr);
\draw (u2) [fill=white] circle (\vr);
\draw (u3) [fill=white] circle (\vr);
\draw (u4) [fill=white] circle (\vr);
\draw (u5) [fill=white] circle (\vr);
\draw (u6) [fill=white] circle (\vr);
\draw (u7) [fill=white] circle (\vr);
\draw (u8) [fill=white] circle (\vr);
\draw (u9) [fill=white] circle (\vr);
\draw[anchor = north] (u1) node {{\small $y$}};
\draw[anchor = south] (u2p) node {{\small $x$}};
\draw[anchor = north] (u3) node {{\small $u_2$}};
\draw[anchor = south] (u4) node {{\small $u_1$}};
\draw[anchor = south] (u5) node {{\small $u$}};
\draw[anchor = south] (u6) node {{\small $v_1$}};
\draw[anchor = south] (u7) node {{\small $v$}};
\draw[anchor = west] (u8) node {{\small $v_3$}};
\draw[anchor = west] (u9) node {{\small $v_2$}};
\end{tikzpicture}
\end{center}
\begin{center}
\vskip -0.65 cm
\caption{A subgraph in the proof of Claim~\ref{claim.22}}
\label{rdom:fig-35}
\end{center}
\end{figure}
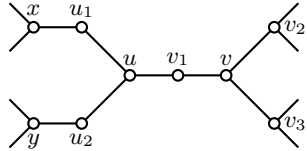

\begin{subclaim}
\label{claim.22.3}
$xy \notin E(G)$.
\end{subclaim}
\proof  Suppose that $xy \in E(G)$. Thus, $C \colon x u_1 u u_2 y x$ is a $5$-cycle in $G$. Let $x_1$ and $y_1$ be the neighbors of $x$ and $y$, respectively, that do not belong to the $5$-cycle $C$. (Possibly, $x_1 = y_1$.) By Claim~\ref{claim.21}, $x_1,y_1 \in \cL$. Let $G'$ be the special subcubic graph obtained from $G - V(C)$ by adding the edge $v_1x_1$. Let $G_x$ be the component of $G'$ containing the added edge $v_1x_1$, and let $G_y$ be the component of $G'$ containing $y_1$. If $G'$ is connected, then $G_x = G_y$. If $G'$ is disconnected, then $G'$ has two components, $G_x$ and $G_y$. In this case, since the removal of a bridge cannot create a component in $\cB_{\rdom}$, we note that $G_y \notin \cB_{\rdom}$.

Let $S'$ be a $\gamma_r$-set of $G'$. If $x_1 \in S'$, let $S = S' \cup \{u,u_2\}$. If $x_1 \notin S'$ and $v_1 \in S'$, let $S = S' \cup \{x,y\}$. If $x_1 \notin S'$ and $v_1 \notin S'$, let $S = S' \cup \{u_1,y\}$. In all cases, $S$ is a RD-set of $G$, and so $\gamma_r(G) \le |S| \le |S'| + 2 = \gamma_r(G') + 2$, implying that $\w(G) < \w(G') + 20$. If $G_x \notin \cB_{\rdom}$, then $\w(G) = w(G') + 21$, a contradiction. Hence, $G_x \in \cB_{\rdom}$. 

Let $G_x^*$ be the component of $G - \{u_2,y\}$ that contains the vertex~$x$. Thus, $G_x^*$ is obtained from the graph $G_x$ by subdividing the added edge $v_1x_1$ three times resulting in the path $v_1 u u_1 x x_1$. Recall that $G_x \in \cB_{\rdom}$.

Suppose that $G_x \ne R_6$ or $G_x = R_6$ and the added edge $v_1x_1$ does not belong to the $4$-cycle in $G_x$. Let $S_x^*$ be a minimum type-$2$ NeRD-set of $G_x^*$ with respect to the vertex $x$. Thus the set $S_x^*$ is a dominating set in $G_x^*$. Further, $x \notin S_x^*$ and the vertex $x$ is the only possible vertex in $G_x^*$ with all its neighbors in~$S_x^*$. By Observation~\ref{obser-4}, we have $|S_x^*| = \gamma_{r,\dom}(G_x^*;x) \le \gamma_{r}(G_x)$.  Let $S^* = S_x^* \cup \{u_2\}$.  If $G'$ is connected, then $\gamma_{r}(G_x) = \gamma_r(G')$ and $S^*$ is a RD-set of $G$. In this case, $\gamma_r(G) \le |S^*| = |S_x^*| + 1 \le \gamma_r(G') + 1$. If $G'$ is disconnected, then $G_x \ne G_y$ and $S^* \cup S_y$ is a RD-set of $G$, where $S_y$ is a $\gamma_r$-set of $G_y$. In this case, $\gamma_r(G) \le |S^*| + |S_y| = |S_x^*| + 1 + |S_y| \le \gamma_{r}(G_x) + 1 + \gamma_r(G_y) = \gamma_r(G') + 1$. In both cases, $\gamma_r(G) \le \gamma_r(G') + 1$, implying that $\w(G) < \w(G') + 10$. However, $\w(G) \ge \w(G') + 17$, a contradiction.

Hence, $G_x = R_6$ and the added edge $v_1x_1$ belongs to the $4$-cycle in $G_x$. We note that $v_1$ is necessarily the vertex of degree~$2$ in the $4$-cycle in $R_6$. By the structure of $R_6$, the vertex~$x_1$ is adjacent to exactly one of~$v_2$ and~$v_3$. Renaming vertices if necessary, we may assume that $x_1v_2$ is an edge. If we had added the edge $v_1y_1$ instead of adding the edge $v_1x_1$, then we would infer analogously that vertex~$y_1$ is adjacent to exactly one of~$v_2$ and~$v_3$. Since~$y_1$ has degree~$2$ in the graph $G_x$, the vertex~$y_1$ is necessarily the vertex of degree~$2$ in $G_x$ adjacent to~$v_3$. Thus, $G_x$ is the graph illustrated in Figure~\ref{fig:special-case}(a). Reconstructing the graph $G$ from $G_x$ now yields the graph illustrated in Figure~\ref{fig:special-case}(b). Thus, $\gamma_r(G) = 5$, where a $\gamma_r$-set of $G$ is indicated by the shaded vertices in Figure~\ref{fig:special-case}(b). However, $\w(G) = 60$, and so $10\gamma_r(G) \le \w(G)$, a contradiction. This completes the proof of Claim~\ref{claim.22.3}.~\smallqed

\vskip -1 cm 
\begin{figure}[htb]
\begin{center}
\begin{tikzpicture}[scale=.75,style=thick,x=0.75cm,y=0.75cm]
\def\vr{2.5pt} 

\path (8,0) coordinate (u1);
\path (8,2.4) coordinate (u2);
\path (9.25,0.75) coordinate (u3);
\path (10.5,1) coordinate (u4);
\path (10.5,1.7) coordinate (u5);
\path (10.5,3.2) coordinate (u6);
\path (11.75,0.75) coordinate (u7);
\path (13,0) coordinate (u8);
\path (13,2.4) coordinate (u9);
\draw (u1)--(u2)--(u6)--(u9)--(u8)--(u1);
\draw (u2)--(u5)--(u9);
\draw (u1)--(u3)--(u4)--(u7)--(u8);
\draw (u4)--(u5);
\draw (u1) [fill=white] circle (\vr);
\draw (u2) [fill=white] circle (\vr);
\draw (u3) [fill=white] circle (\vr);
\draw (u4) [fill=white] circle (\vr);
\draw (u5) [fill=white] circle (\vr);
\draw (u6) [fill=white] circle (\vr);
\draw (u7) [fill=white] circle (\vr);
\draw (u8) [fill=white] circle (\vr);
\draw (u9) [fill=white] circle (\vr);
\draw[anchor = west] (u8) node {{\small $v_3$}};
\draw[anchor = east] (u2) node {{\small $x_1$}};
\draw[anchor = west] (u9) node {{\small $v$}};
\draw[anchor = south] (u7) node {{\small $y_1$}};
\draw[anchor = south] (u5) node {{\small $v_2$}};
\draw[anchor = south] (u6) node {{\small $v_1$}};
\draw (10.5,-1) node {{\small (a) $G_x = R_6$}};
\path (18,0) coordinate (u1);
\path (18,2.4) coordinate (u2);
\path (19.25,0.75) coordinate (u3);
\path (20.5,1) coordinate (u4);
\path (20.5,1.7) coordinate (u5);
\path (19,3.2) coordinate (a);
\path (19,4.2) coordinate (a1);
\path (20,3.2) coordinate (b);
\path (21,3.2) coordinate (c);
\path (21,4.2) coordinate (c1);
\path (22,3.2) coordinate (d);
\path (21.75,0.75) coordinate (u7);
\path (23,0) coordinate (u8);
\path (23,2.4) coordinate (u9);
\draw (u1)--(u2)--(a)--(b)--(c)--(d)--(u9)--(u8)--(u1);
\draw (u2)--(u5)--(u9);
\draw (u1)--(u3)--(u4)--(u7)--(u8);
\draw (u4)--(u5);
\draw (a)--(a1)--(c1)--(c);
\draw (a1) to[out=35,in=0, distance=3.5cm] (u7);
\draw (u1) [fill=white] circle (\vr);
\draw (u2) [fill=black] circle (\vr);
\draw (u3) [fill=white] circle (\vr);
\draw (u4) [fill=black] circle (\vr);
\draw (u5) [fill=black] circle (\vr);
\draw (u7) [fill=black] circle (\vr);
\draw (u8) [fill=white] circle (\vr);
\draw (u9) [fill=white] circle (\vr);
\draw (a) [fill=white] circle (\vr);
\draw (a1) [fill=white] circle (\vr);
\draw (b) [fill=white] circle (\vr);
\draw (c) [fill=black] circle (\vr);
\draw (c1) [fill=white] circle (\vr);
\draw (d) [fill=white] circle (\vr);
\draw[anchor = west] (u8) node {{\small $v_3$}};
\draw[anchor = east] (u2) node {{\small $x_1$}};
\draw[anchor = north] (a) node {{\small $x$}};
\draw[anchor = south] (a1) node {{\small $y$}};
\draw[anchor = north] (b) node {{\small $u_1$}};
\draw[anchor = north] (c) node {{\small $u$}};
\draw[anchor = south] (c1) node {{\small $u_2$}};
\draw[anchor = north] (d) node {{\small $v_1$}};
\draw[anchor = west] (u9) node {{\small $v$}};
\draw[anchor = south] (u7) node {{\small $y_1$}};
\draw[anchor = south] (u5) node {{\small $v_2$}};
\draw (20.5,-1) node {{\small (b) $G$}};

\end{tikzpicture}
\end{center}
\begin{center}
\vskip -0.70 cm
\caption{The graph $G$ in the proof of Claim~\ref{claim.22.3} when $G_x = R_6$}
\label{fig:special-case}
\end{center}
\end{figure}

\medskip
We now return to the proof of Claim~\ref{claim.22}. By Claim~\ref{claim.22.3}, the vertices $x$ and $y$ are not adjacent in $G$. Let $G'$ be the special subcubic graph obtained from $G - \{u,u_1,u_2,v,v_1\}$ by adding the edge $xy$. Let $G_x$ be the component of $G'$ containing the added edge $xy$, and let $G_2$ and $G_3$ be the components of $G'$ containing $v_2$ and $v_3$, respectively. If $G'$ is connected, then $G_x = G_2 = G_3$. Let $S'$ be a $\gamma_r$-set of $G'$. If $x \in S'$, let $S = S' \cup \{u_2,v\}$. If $x \notin S'$ and $y \in S'$, let $S = S' \cup \{u_1,v\}$. If $x \notin S'$, $y \notin S'$ and $v_2 \in S$, let $S = S' \cup \{u\}$. If $x \notin S'$, $y \notin S'$ and $v_2 \notin S$, let $S = S' \cup \{u,v_1\}$. In all cases, $S$ is a RD-set of $G$, and so $\gamma_r(G) \le |S| \le |S'| + 2 = \gamma_r(G') + 2$, implying that $\w(G) < \w(G') + 20$. If no component of $G'$ belongs to $\cB_{\rdom}$, then $\w(G) = \w(G') + 21$, a contradiction. Hence, there is a component in $G'$ that belongs to $\cB_{\rdom}$.

Suppose that $G_2$ or $G_3$ is different from $G_x$ and belongs to $\cB_{\rdom}$. Renaming vertices if necessary, by symmetry we may assume that $G_2 \ne G_x$ and $G_2 \in \cB_{\rdom}$. Since the removal of a bridge cannot create a component in $\cB_{\rdom}$, we infer that $G_2 = G_3$. Further, both $v_2$ and~$v_3$ have degree~$2$ in $G_2$. Applying Observation~\ref{obser-1}(f) to the graph $G_2$ with $X = \{v_2,v_3\}$, we have $\gamma_{r,\dom}(G_2;X) \le \gamma_{r}(G_2) - 1$. Let $S^*$ be a minimum type-$2$ NeRD-set of $G_2$ with respect to the set $X$. Let $S_x$ be a $\gamma_r$-set of $G_x$. If $x \in S_x$, let $S = S_x \cup S^* \cup \{u_2,v_1\}$. If $x \notin S_x$ and $y \in S_x$, let $S = S_x \cup S^* \cup \{u_1,v_1\}$. If $x \notin S_x$ and $y \notin S_x$, let $S = S_x \cup S^* \cup \{u,v_1\}$. In all cases, $S$ is a RD-set of $G$, and so $\gamma_r(G) \le |S_x| + |S^*| + 2 \le \gamma_r(G_x) + (\gamma_{r}(G_2) - 1) + 2 = \gamma_{r}(G') + 1$, implying that $\w(G) < \w(G') + 10$. However, $\w(G) \ge \w(G') + 13$, a contradiction. Hence if $G_2 \ne G_x$, then $G_2 \notin \cB_{\rdom}$, and if $G_3 \ne G_x$, then $G_3 \notin \cB_{\rdom}$. 

Since there is a component in $G'$ that belongs to $\cB_{\rdom}$, we infer that $G_x$ is the only such component of $G'$. If $G_x \in \cB_{\rdom,1}$, then $\w(G) = \w(G') + 20$, a contradiction. Hence, $G_x \notin \{R_6,R_7,R_8,R_9\}$. We note that $x$ and $y$ are adjacent vertices of degree~$3$ in $G_x$, implying that $G_x \ne \{R_1,R_2\}$. If $G_x = R_3$, then our properties of the graph $G$ imply that $G'$ is connected and $v_2$ and $v_3$ are the vertices of degree~$2$ in $R_3$ that have no degree~$3$ neighbor. In this case, the graph $G$ is determined and $\gamma_r(G) = 4$ and $\w(G) = 59$, a contradiction. Hence, $G_x \ne R_3$, implying that $G_x \in \{R_4,R_5,R_{10}\}$. 

Let $G^*$ be obtained from $G_x$ by subdividing the added edge $xy$ three times resulting in the path $x u_1 u u_2 y$. Applying Observation~\ref{obser-4}(b) we have $\gamma_{r,\dom}(G^*;u) \le \gamma_{r}(G_x)$. Thus, there exists a type-$2$ NeRD-set $S^*$ in $G^*$ with respect to the vertex $u$ such that $|S^*| \le \gamma_{r}(G_x)$. The set $S^*$ is a dominating set in $G^*$. Further, $u \notin S^*$ and the vertex $u$ is the only possible vertex in $G^*$ with all its neighbors in~$S^*$. Let $S = S^* \cup \{v\}$. If $G'$ is connected, then $S^* \cup \{v\}$ is a RD-set of $G$, and so $\gamma_r(G) \le |S^*| + 1 \le \gamma_{r}(G_x) + 1 = \gamma_{r}(G') + 1$. If $G'$ is disconnected, then every $\gamma_r$-set of $G' - V(G_x)$ can be extended to a RD-set of $G$ by adding to it the set $S^* \cup \{v\}$, implying once again that $\gamma_r(G) \le \gamma_{r}(G') + 1$. Hence in both cases we infer that $\w(G) < \w(G') + 10$. However, $\w(G) = w(G') + 17$, a contradiction. This completes the proof of Claim~\ref{claim.22}.~\smallqed

\begin{claim}
\label{claim.cubic}
The graph $G$ is a cubic graph.
\end{claim}
\proof Suppose, to the contrary, that $G$ contains a small vertex. By Claim~\ref{claim.21}, no large vertex has exactly two small neighbors. By Claim~\ref{claim.22}, no large vertex has exactly one small neighbor. Hence if a large vertex has a small neighbor, then all three of its neighbors are small. Thus the three neighbors of every large vertex are either all small or all large. Since $G$ is connected and contains at least one small vertex, this implies that $G$ is a bipartite subcubic graph with partite sets $\cS$ and $\cL$. Thus, by Lemma~\ref{lem:1}, $\gamma_r(G) \le |\cL|$, and so $\w(G) < 10\gamma_r(G) \le 10|\cL|$. However in this case, $3|\cL| = 2|\cS|$, and so $\w(G) = 5|\cS| + 4|\cL| = 5 \times \frac{3}{2}|\cL| + 4|\cL| > 10|\cL$, a contradiction.~\smallqed

\medskip
By Claim~\ref{claim.cubic}, $G$ is a (connected) cubic graph. Recall by Claim~\ref{claim.no.R10} that $R_{11}$ is not a subgraph of~$G$. We note that $R_{10}$ contains three small vertices, and every graph in $\cB_{\rdom} \setminus \{R_{10},R_{11}\}$ contains at least four small vertices. Our earlier observations therefore yield the following properties of the graph $G$.

\begin{claim}
\label{claim:edge-cut}
If $E'$ is a $k$-edge-cut in $G$ and $G'$ is a component of $G - E'$ that belongs to~$\cB_{\rdom}$, then $k \ge 3$ and the following properties hold. \\[-22pt]
\begin{enumerate}
\item[{\rm (a)}] If $k = 3$, then $G' = R_{10}$.
\item[{\rm (b)}] If $k = 4$, then $G' \in \{R_2,R_4,R_5,R_{10}\}$.
\item[{\rm (c)}] If $k = 5$, then $G' \in \{R_1,R_7,R_8,R_9\}$.
\end{enumerate}
\end{claim}

\begin{claim}
\label{claim:bad-comp}
If $G' \in \cB_{\rdom}$ is a special subcubic component of $G - S$ where $S \subset V(G)$, then $G'$ contains at least three vertices of degree~$2$.
\end{claim}

\begin{claim}
\label{claim.no-diamond}
The graph $G$ contains no diamond.
\end{claim}
\proof Suppose, to the contrary, that $G$ contains a diamond $D$, where $V(D) = \{v_1,v_2,v_3,v_4\}$ and where $v_1v_2$ is the missing edge in $D$. Let $u_i$ be the neighbor of $v_i$ not in $D$ for $i \in [2]$. Suppose that $u_1 = u_2$. Let $u$ be the neighbor of $u_1$ different from $v_1$ and $v_2$, and let $G' = G - (V(D) \cup \{u,u_1\})$. The graph $G'$ is a special subcubic graph that contains exactly two small vertices, and so by Claim~\ref{claim:bad-comp} no component of $G'$ belongs to $\cB_{\rdom}$. Every $\gamma_r$-set of $G'$ can be extended to a RD-set of $G$ by adding to it the set $\{v_3,u\}$, and so $\gamma_r(G) \le \gamma_r(G) + 2$, implying that $\w(G) < \w(G') + 20$. However, $\w(G) = \w(G) + 22$, a contradiction. Hence, $u_1 \ne u_2$. In this case, let $G' = G - V(D)$. The graph $G'$ is a special subcubic graph that contains exactly two small vertices, and so no component of $G'$ belongs to $\cB_{\rdom}$. Every $\gamma_r$-set of $G'$ can be extended to a RD-set of $G$ by adding to it the vertex $v_3$, and so $\gamma_r(G) \le \gamma_r(G) + 1$, implying that $\w(G) < \w(G') + 10$. In this case, $\w(G) = \w(G) + 14$, a contradiction.~\smallqed

\begin{claim}
\label{claim.no-triangle}
The graph $G$ contains no triangle.
\end{claim}
\proof Suppose, to the contrary, that $T$ is a triangle in $G$ where $V(T) = \{v_1, v_2, v_3\}$. Let $x_i$ be the third neighbor of $v_i$ that does not belong to $T$ for $i \in [3]$. By Claim~\ref{claim.no-diamond}, the graph $G$ contains no diamond, and so the vertices $x_1$, $x_2$ and $x_3$ are pairwise distinct. Let $X = \{x_1,x_2,x_3\}$. Suppose that $G[X]$ contains a vertex of degree~$2$. Renaming vertices if necessary, we may assume that $\{x_1x_2,x_2x_3\} \subset E(G)$. If $x_1x_3 \in E(G)$, then $G$ is the $3$-prism $C_3 \cp K_2$, and so $\gamma_r(G) = 2$ and $\w(G) = 24$, a contradiction. Hence, $x_1x_3 \notin E(G)$. Let $y_i$ be the neighbor of $x_i$ different from $x_2$ and $v_i$ for $i \in \{1,3\}$. If $y_1 = y_3$, then we let $Q = V(T) \cup X \cup \{y_1\}$ and $G' = G - Q$. In this case $G'$ is a special connected subcubic graph that contains exactly one small vertex, and so, by Claim~\ref{claim:bad-comp}, $G' \notin \cB_{\rdom}$. Since $\gamma_r(G) \le \gamma_r(G) + 2$, we have $\w(G) < \w(G') + 20$. However, $\w(G) = \w(G) + 27$, a contradiction. Hence, $y_1 \ne y_3$. We now let $Q = V(T) \cup X$ and $G' = G - Q$. In this case, $G'$ is a special subcubic graph that contains exactly two small vertices, and so, by Claim~\ref{claim:bad-comp}, no component of $G'$ belongs to $\cB_{\rdom}$. Once again $\gamma_r(G) \le \gamma_r(G) + 2$, implying that $\w(G) < \w(G') + 20$. However, $\w(G) = 24 + (\w(G') - 2) = \w(G) + 22$, a contradiction.

Hence, $G[X]$ contains no vertex of degree~$2$, implying that $G[X]$ contains at least one isolated vertex. By symmetry, we may assume that $x_1$ is isolated in $G[X]$, that is, $x_1$ is adjacent to neither $x_2$ nor $x_3$. Let $y_1$ and $y_2$ be the two neighbors of $x_1$ different from $v_1$. We now let $Q = V(T) \cup \{x_1\}$ and let $G' = G - Q$. The graph $G'$ is a special subcubic graph. We note that $k' + r' \le 4$.

Let $S'$ be a $\gamma_r$-set of $G'$. If $y_1 \in S'$, let $S = S' \cup \{v_2\}$. If $y_1 \notin S'$, let $S = S' \cup \{v_1\}$. In both cases, $S$ is a RD-set of $G$, and so $\gamma_r(G) \le |S| = |S'| + 1 = \gamma_r(G') + 1$, implying that $\w(G) < \w(G') + 10$. If no component of $G'$ belongs to $\cB_{\rdom}$, then $\w(G) = \w(G') + 12$, a contradiction. Hence, $G'$ contains a component $G_1$ that belongs to $\cB_{\rdom}$. By Claim~\ref{claim:edge-cut}, there is only one such component and $G_1 \in \{R_2,R_4,R_5,R_{10}\}$. Thus, $G_1$ contains at least three small vertices. Let $X_1 \subset V(G_1) \cap \{y_1,y_2,x_2\}$ be chosen so that $|X_1| = 2$. Let $S_1$ be a minimum type-$2$ NeRD-set of $G_1$ with respect to the set $X_1$. By Observation~\ref{obser-1}(f), $|S_1| = \gamma_{r,\dom}(G_1;X_1) \le \gamma_{r}(G_1) - 1$.

Suppose that $G_1 \in \{R_2,R_4,R_5\}$. By Claim~\ref{claim:edge-cut}, $G' = G_1$, and so $k' = 1$ and $r' = 0$. In this case, the set $S_1 \cup \{v_1\}$ is a RD-set of $G$, and so $\gamma_r(G) \le |S_1| + 1 \le \gamma_{r}(G_1) = \gamma_r(G')$. Suppose that $G_1 = R_{10}$, implying that $k' = r' = 1$.  In this case, we let $G_2$ be the second component of $G'$, and so $G_2 \notin \cB_{\rdom}$. Let $S_2$ be a $\gamma_r$-set of $G_2$. The set $S_2$ can be extended to a RD-set of $G$ by adding to it the set $S_1 \cup \{v_1\}$, and so $\gamma_r(G) \le |S_1| + 1 + |S_2| \le \gamma_{r}(G_1) + \gamma_r(G_2) = \gamma_r(G')$. In both cases, $\gamma_r(G) \le \gamma_r(G')$, implying that $\w(G) < \w(G')$. However, $\w(G) \ge \w(G') + 8$, a contradiction.~\smallqed


\begin{claim}
\label{claim.no-K23}
The graph $G$ contains no $K_{2,3}$ as a subgraph.
\end{claim}
\proof
Suppose, to the contrary, that $H$ is a subgraph of $G$, where $H \cong K_{2,3}$. Let $X$ and $Y$ be the partite sets of $H$ where $X = \{x_1,x_2,x_3\}$ and $Y = \{y_1,y_2\}$. Since $G$ is triangle-free, the sets $X$ and $Y$ are independent. Let $v_i$ be the neighbor of $x_i$ not in $H$ for $i \in [3]$. If $v_1 = v_2 = v_3$, then $G = K_{3,3}$. In this case, $\gamma_r(G) = 2$ and $\w(G) = 24$, a contradiction. Hence renaming vertices if necessary, we may assume that $v_1 \ne v_2$.

\begin{subclaim}
\label{claim.no-K23.1}
The vertices $v_1, v_2, v_3$ are pairwise distinct.
\end{subclaim}
\proof
Suppose, to the contrary, that the vertices $v_1, v_2, v_3$ are not pairwise distinct, and so $v_1 = v_3$ or $v_2 = v_3$. Renaming vertices if necessary, we may assume that $v_2 = v_3$. Suppose that $v_1v_2 \in E(G)$. In this case, let $v$ denote the neighbor of $v_1$ different from $x_1$ and $v_2$. Thus, $vv_1$ is a bridge in $G$. Let $G'$ be the component of $G - vv_1$ that contains the vertex~$v$. By Claim~\ref{claim:edge-cut}, $G' \notin \cB_{\rdom}$. Let $S'$ be a $\gamma_r$-set of $G'$. If $v \in S'$, let $S = \{y_1,x_2\}$. If $v \notin S'$, let $S = \{x_1,v_2\}$. In both cases, $S$ is a RD-set of $G$, and so $\gamma_r(G) \le |S| = \gamma_r(G') + 2$, implying that $\w(G) < \w(G') + 20$. However, $\w(G) = w(G') + 27$, a contradiction. Hence, $v_1v_2 \notin E(G)$. We now let $G' = G - (V(H) \cup \{v_2\})$. The graph $G'$ is a special subcubic graph that contains exactly two small vertices. By Claim~\ref{claim:edge-cut}, no component of $G'$ belongs to $\cB_{\rdom}$. Every $\gamma_r$-set of $G'$ can be extended to a RD-set of $G$ by adding to it the set $\{y_1,x_2\}$, and so $\gamma_r(G) \le \gamma_r(G') + 2$, implying that $\w(G) < \w(G') + 20$. However, $\w(G) = \w(G') + 22$, a contradiction.~\smallqed

\medskip
By Claim~\ref{claim.no-K23.1}, the vertices $v_1, v_2, v_3$ are pairwise distinct.

\begin{subclaim}
\label{claim.no-K23.2}
The graph $G[\{v_1,v_2,v_3\}]$ is isolate-free.
\end{subclaim}
\proof
Suppose, to the contrary, that $G[\{v_1,v_2,v_3\}]$ contains an isolated vertex. Renaming vertices if necessary, we may assume that the vertex $v_1$ is adjacent to neither $v_2$ nor $v_3$. Let $G' = G - (V(H) \cup \{v_1\})$. Thus, $G'$ is a special subcubic graph that contains exactly four small vertices. Let $S'$ be a $\gamma_r$-set of $G'$. Let $u_1$ and $u_2$ be two neighbors of $v_1$ in $G$ different from $x_1$. If $u_1 \notin S'$, let $S = S' \cup \{x_1,y_1\}$. If $u_1 \in S'$, let $S = S' \cup \{y_1,x_2\}$. In both cases, $S$ is a RD-set of $G$, and so $\gamma_r(G) \le |S| = \gamma_r(G') + 2$, implying that $\w(G) < \w(G') + 20$. If no component of $G'$ belongs to $\cB_{\rdom}$, then $\w(G) = \w(G') + 20$, a contradiction. Hence, $G'$ contains a component $G_1$ that belongs to $\cB_{\rdom}$. By Claim~\ref{claim:edge-cut}, there is only one such component and $G_1 \in \{R_2,R_4,R_5,R_{10}\}$.

At least one of $u_1$ and $u_2$, and at least one of $v_2$ and $v_3$ belong to $G_1$. Renaming vertices if necessary, we may assume that $\{u_1,v_2\} \subset V(G_1)$. Let $S_1$ be a minimum type-$2$ NeRD-set of $G_1$ with respect to the set $X_1 = \{u_1,v_2\}$. By Observation~\ref{obser-1}(f), $|S_1| =  \gamma_{r,\dom}(G_1;X_1) \le \gamma_{r}(G_1) - 1$. Suppose that $G_1 \in \{R_2,R_4,R_5\}$. By Claim~\ref{claim:edge-cut}, $G' = G_1$, and so $k' = 1$ and $r' = 0$. In this case, the set $S_1 \cup \{x_1,y_1\}$ is a RD-set of $G$, and so $\gamma_r(G) \le |S_1| + 2 \le \gamma_{r}(G_1) + 1 = \gamma_r(G') + 1$. Suppose that $G_1 = R_{10}$, implying that $k' = r' = 1$.  In this case, we let $G_2$ be the second component of $G'$, and so $G_2 \notin \cB_{\rdom}$. Let $S_2$ be a $\gamma_r$-set of $G_2$. The set $S_2$ can be extended to a RD-set of $G$ by adding to it the set $S_1 \cup \{x_1,y_1\}$, and so $\gamma_r(G) \le |S_1| + 2 + |S_2| \le \gamma_{r}(G_1) + \gamma_r(G_2) + 1 = \gamma_r(G') + 1$. In both cases, $\gamma_r(G) \le \gamma_r(G') + 1$, implying that $\w(G) < \w(G') + 10$. However, $\w(G) \ge \w(G') + 16$, a contradiction.~\smallqed

\medskip
By Claim~\ref{claim.no-K23.2}, the graph $G[\{v_1,v_2,v_3\}]$ is isolate-free. Renaming vertices if necessary, we may assume that $v_1v_2$ and $v_2v_3$ are edges. Since $G$ is triangle-free, we note that $v_1v_3$ is not an edge. Let $u_1$ be the neighbor of $v_1$ different from $x_1$ and $v_2$, and let $u_3$ be the neighbor of $v_3$ different from $x_3$ and $v_2$. Suppose that $u_1 \ne u_3$. Hence, the graph illustrated in Figure~\ref{rdom:fig-no-K23}(a) is a subgraph of $G$. In this case, let $Q = V(H) \cup \{v_1,v_2,v_3\}$ and let $G' = G - Q$. We note that $G'$ is a special subcubic graph and is obtained by deleting the edges of a $2$-edge-cut in $G$. By Claim~\ref{claim:edge-cut}, no component of $G'$ belongs to $\cB_{\rdom}$. Every $\gamma_r$-set of $G'$ can be extended to a RD-set of $G$ by adding to it the set $\{v_2,x_2,y_2\}$, and so $\gamma_r(G) \le \gamma_r(G') + 3$, implying that $\w(G) < \w(G') + 30$. However, $\w(G) = \w(G') + 30$, a contradiction.

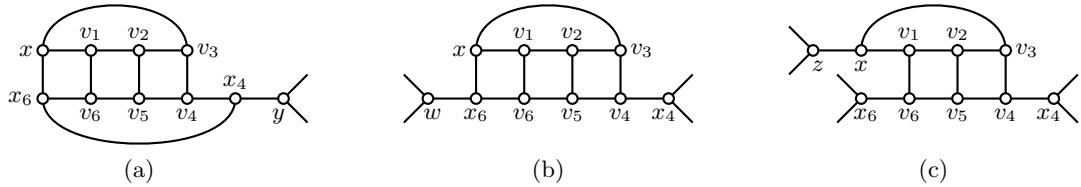
\begin{figure}[htb]
\begin{center}
\begin{tikzpicture}[scale=.8,style=thick,x=0.8cm,y=0.8cm]
\def\vr{2.5pt} 
\path (0,1) coordinate (u1);
\path (0,2) coordinate (u2);
\path (0,3) coordinate (u3);
\path (-0.5,3.5) coordinate (u31);
\path (0.5,3.5) coordinate (u32);
\path (1,0) coordinate (u4);
\path (2,1) coordinate (u5);
\path (2,2) coordinate (u6);
\path (3,0) coordinate (u7);
\path (4,1) coordinate (u8);
\path (4,2) coordinate (u9);
\path (4,3) coordinate (u10);
\path (3.5,3.5) coordinate (u101);
\path (4.5,3.5) coordinate (u102);

\draw (u4)--(u1)--(u2)--(u3);
\draw (u7)--(u8)--(u9)--(u10);
\draw (u4)--(u5)--(u7);
\draw (u2)--(u6)--(u9);
\draw (u31)--(u3)--(u32);
\draw (u5)--(u6);
\draw (u4)--(u8);
\draw (u1)--(u7);
\draw (u101)--(u10)--(u102);
\draw (u1) [fill=white] circle (\vr);
\draw (u2) [fill=white] circle (\vr);
\draw (u3) [fill=white] circle (\vr);
\draw (u4) [fill=white] circle (\vr);
\draw (u5) [fill=white] circle (\vr);
\draw (u6) [fill=white] circle (\vr);
\draw (u7) [fill=white] circle (\vr);
\draw (u8) [fill=white] circle (\vr);
\draw (u9) [fill=white] circle (\vr);
\draw (u10) [fill=white] circle (\vr);
\draw[anchor = east] (u1) node {{\small $x_1$}};
\draw[anchor = east] (u2) node {{\small $v_1$}};
\draw[anchor = east] (u3) node {{\small $u_1$}};
\draw[anchor = north] (u4) node {{\small $y_1$}};
\draw[anchor = east] (u5) node {{\small $x_2$}};
\draw[anchor = south] (u6) node {{\small $v_2$}};
\draw[anchor = north] (u7) node {{\small $y_2$}};
\draw[anchor = west] (u8) node {{\small $x_3$}};
\draw[anchor = west] (u9) node {{\small $v_3$}};
\draw[anchor = west] (u10) node {{\small $u_3$}};
%
\draw (2,-1) node {{\small (a)}};
%
\path (8,1) coordinate (u1);
\path (8,2) coordinate (u2);
\path (10,3) coordinate (u);
\path (9.9,3.1) coordinate (up);
\path (10,4) coordinate (u3);
\path (9.5,4.5) coordinate (w31);
\path (10.5,4.5) coordinate (w32);
\path (9,0) coordinate (u4);
\path (10,1) coordinate (u5);
\path (10,2) coordinate (u6);
\path (11,0) coordinate (u7);
\path (12,1) coordinate (u8);
\path (12,2) coordinate (u9);

\draw (u4)--(u1)--(u2)--(u);
\draw (u7)--(u8)--(u9)--(u);
\draw (u4)--(u5)--(u7);
\draw (u2)--(u6)--(u9);
\draw (w31)--(u3)--(w32);
\draw (u5)--(u6);
\draw (u4)--(u8);
\draw (u1)--(u7);
\draw (u2)--(u)--(u3);
%
\draw (u) [fill=white] circle (\vr);
\draw (u1) [fill=white] circle (\vr);
\draw (u2) [fill=white] circle (\vr);
\draw (u3) [fill=white] circle (\vr);
\draw (u4) [fill=white] circle (\vr);
\draw (u5) [fill=white] circle (\vr);
\draw (u6) [fill=white] circle (\vr);
\draw (u7) [fill=white] circle (\vr);
\draw (u8) [fill=white] circle (\vr);
\draw (u9) [fill=white] circle (\vr);
%
\draw[anchor = east] (up) node {{\small $u$}};
\draw[anchor = east] (u1) node {{\small $x_1$}};
\draw[anchor = east] (u2) node {{\small $v_1$}};
\draw[anchor = east] (u3) node {{\small $w$}};
\draw[anchor = north] (u4) node {{\small $y_1$}};
\draw[anchor = east] (u5) node {{\small $x_2$}};
\draw[anchor = south] (u6) node {{\small $v_2$}};
\draw[anchor = north] (u7) node {{\small $y_2$}};
\draw[anchor = west] (u8) node {{\small $x_3$}};
\draw[anchor = west] (u9) node {{\small $v_3$}};
%
\draw (10,-1) node {{\small (b)}};
%
\end{tikzpicture}
\end{center}
\begin{center}
\vskip -0.5 cm
\caption{Subgraphs in the proof of Claim~\ref{claim.no-K23}}
\label{rdom:fig-no-K23}
\end{center}
\end{figure}

Hence, $u_1 = u_3$ and let us rename this common neighbor of $v_1$ and $v_3$ by $u$. Let $w$ be the third neighbor of $u$ different from $v_1$ and $v_3$. Thus, the graph illustrated in Figure~\ref{rdom:fig-no-K23}(b) is a subgraph of $G$. In this case, let $Q = V(H) \cup \{v_1,v_2,v_3,u\}$ and let $G' = G - Q$. We note that $G'$ is a connected special subcubic graph and is obtained by deleting the cut-edge $uw$ in $G$. By Claim~\ref{claim:edge-cut}, $G' \notin \cB_{\rdom}$. Every $\gamma_r$-set of $G'$ can be extended to a RD-set of $G$ by adding to it the set $\{u,x_2,y_2\}$, and so $\gamma_r(G) \le \gamma_r(G') + 3$, implying that $\w(G) < \w(G') + 30$. However, $\w(G) = \w(G') + 35$, a contradiction. This completes the proof of Claim~\ref{claim.no-K23}.~\smallqed


\begin{claim}
\label{claim.no-domino}
The graph $G$ contains no domino as a subgraph.
\end{claim}
\proof
Suppose, to the contrary, that $G$ contains a domino $F$ as a subgraph. Let $V(F) = \{v_1,v_2,\ldots, v_6\}$ where $v_1v_2 \ldots v_6v_1$ is a $6$-cycle and $v_2v_5$ is an edge. Since $G$ is triangle-free and $K_{2,3}$-free, we note that $F$ is an induced subgraph of $G$. Let $x_i$ be the neighbor of $v_i$ that does not belong to $F$ for $i \in \{1,3,4,6\}$. Since $G$ is triangle-free, $x_1 \ne x_6$ and $x_3 \ne x_4$.

\begin{subclaim}
\label{claim.no-domino.1}
$x_1 \ne x_3$ and $x_4 \ne x_6$.
\end{subclaim}
\proof Suppose, to the contrary, that $x_1 = x_3$ or $x_4 = x_6$. Renaming vertices if necessary, we may assume by symmetry that $x_1 = x_3$. Thus, $x_1$ is a common neighbor of $v_1$ and $v_3$ different from $v_2$. Let us rename the vertex $x_1$ by~$x$ for notational simplicity.

Suppose firstly that $x_4 = x_6$, and so $x_4$ is a common neighbor of $v_4$ and $v_6$ different from $v_5$. Let us rename the vertex $x_4$ by~$y$ for notational simplicity. If $xy \in E(G)$, then the graph $G$ is determined and $\gamma_r(G) = 2$ and $\w(G) = 32$, a contradiction. Hence, $xy \notin E(G)$. Let $x_1$ and $y_1$ be the neighbors of $x$ and $y$, respectively, that do not belong to $F$.
Suppose that $x_1 = y_1$, and let us rename this common neighbor of $x$ and $y$ by~$w$. Let $z$ be the third neighbor of $w$ different from $x$ and $y$. In this case, let $G'$ be the component of $G - wz$ that contains the vertex~$z$. We note that $G'$ is a connected special subcubic graph and the vertex $z$ is the only vertex of degree~$2$ in $G'$, and so $G' \notin \cB_{\rdom}$. Every $\gamma_r$-set of $G'$ can be extended to a RD-set of $G$ by adding to it the set $\{v_1,v_4,w\}$, and so $\gamma_r(G) \le \gamma_r(G') + 3$, implying that $\w(G) < \w(G') + 30$. However, $\w(G) = \w(G') + 35$, a contradiction. Hence, $x_1 \ne y_1$.
In this case, we let $G' = G - (V(F) \cup \{x,y\})$. We note that $G'$ is a special subcubic graph that contains exactly two vertices of degree~$2$. By Claim~\ref{claim:edge-cut}, no component of $G'$ belongs to $\cB_{\rdom}$. Every $\gamma_r$-set of $G'$ can be extended to a RD-set of $G$ by adding to it the set $\{v_1,v_4\}$, and so $\gamma_r(G) \le \gamma_r(G') + 2$, implying that $\w(G) < \w(G') + 20$. However, $\w(G) = w(G') + 30$, a contradiction.

Hence, $x_4 \ne x_6$, that is, $v_5$ is the only common neighbor of $v_4$ and $v_6$. Since $G$ is triangle-free, $x \ne x_4$ and $x \ne x_6$, that is, the vertices $x, x_4, x_6$ are pairwise distinct. Suppose that $x$ is adjacent to $x_4$ or $x_6$. Renaming vertices if necessary, we may assume $xx_6 \in E(G)$. Suppose that $x_4x_6 \in E(G)$. In this case, let $y$ be the neighbor of $x_4$ different from $v_4$ and $x_6$. Hence, the graph illustrated in Figure~\ref{rdom:fig-no-dom}(a) is a subgraph of $G$. Let $G'$ be the component of $G - x_4y$ that contains the vertex~$y$. We note that $G' \notin \cB_{\rdom}$. Every $\gamma_r$-set of $G'$ can be extended to a RD-set of $G$ by adding to it the set $\{x,x_4,v_5\}$, and so $\gamma_r(G) \le \gamma_r(G') + 3$, implying that $\w(G) < \w(G') + 30$. However, $\w(G) = \w(G') + 35$, a contradiction. Hence, $x_4x_6 \notin E(G)$. In this case, let $w$ be the neighbor of $x_6$ different from $x$ and $v_6$. Hence, the graph illustrated in Figure~\ref{rdom:fig-no-dom}(b) is a subgraph of $G$. We now let $G' = G - (V(F) \cup \{x,x_6\})$. The special subcubic graph $G'$ contains exactly two vertices of degree~$2$, and so by Claim~\ref{claim:edge-cut} no component of $G'$ belongs to $\cB_{\rdom}$. Every $\gamma_r$-set of $G'$ can be extended to a RD-set of $G$ by adding to it the set $\{v_1,v_4,x_6\}$, and so $\gamma_r(G) \le \gamma_r(G') + 3$, implying that $\w(G) < \w(G') + 30$. However, $\w(G) = \w(G') + 30$, a contradiction.

\begin{figure}[htb]
\begin{center}
\begin{tikzpicture}[scale=.8,style=thick,x=0.8cm,y=0.8cm]
\def\vr{2.5pt} 
\path (0,0) coordinate (u1);
\path (0,1) coordinate (u2);
\path (1,0) coordinate (u3);
\path (1,1) coordinate (u4);
\path (2,0) coordinate (u5);
\path (2,1) coordinate (u6);
\path (3,0) coordinate (u7);
\path (3,1) coordinate (u8);
\path (4,0) coordinate (u9);
\path (5,0) coordinate (u10);
\path (4.9,0) coordinate (u10p);
\path (5.5,-0.5) coordinate (u101);
\path (5.5,0.5) coordinate (u102);
\draw (u2)--(u1)--(u3)--(u5)--(u7)--(u9)--(u10);
\draw (u2)--(u4)--(u6)--(u8);
\draw (u101)--(u10)--(u102);
\draw (u3)--(u4);
\draw (u5)--(u6);
\draw (u7)--(u8);
\draw (u1) to[out=270,in=270, distance=1cm] (u9);
\draw (u2) to[out=90,in=90, distance=1cm] (u8);
\draw (u1) [fill=white] circle (\vr);
\draw (u2) [fill=white] circle (\vr);
\draw (u3) [fill=white] circle (\vr);
\draw (u4) [fill=white] circle (\vr);
\draw (u5) [fill=white] circle (\vr);
\draw (u6) [fill=white] circle (\vr);
\draw (u7) [fill=white] circle (\vr);
\draw (u8) [fill=white] circle (\vr);
\draw (u9) [fill=white] circle (\vr);
\draw (u10) [fill=white] circle (\vr);
\draw[anchor = east] (u1) node {{\small $x_6$}};
\draw[anchor = east] (u2) node {{\small $x$}};
\draw[anchor = north] (u3) node {{\small $v_6$}};
\draw[anchor = south] (u4) node {{\small $v_1$}};
\draw[anchor = north] (u5) node {{\small $v_5$}};
\draw[anchor = south] (u6) node {{\small $v_2$}};
\draw[anchor = north] (u7) node {{\small $v_4$}};
\draw[anchor = west] (u8) node {{\small $v_3$}};
\draw[anchor = south] (u9) node {{\small $x_4$}};
\draw[anchor = north] (u10p) node {{\small $y$}};
%
\draw (2,-1.5) node {{\small (a)}};
%
\path (7.5,-0.5) coordinate (u01);
\path (7.5,0.5) coordinate (u02);
\path (8,0) coordinate (u0);
\path (8.1,0) coordinate (u0p);
\path (9,0) coordinate (u1);
\path (9,1) coordinate (u2);
\path (10,0) coordinate (u3);
\path (10,1) coordinate (u4);
\path (11,0) coordinate (u5);
\path (11,1) coordinate (u6);
\path (12,0) coordinate (u7);
\path (12,1) coordinate (u8);
\path (13,0) coordinate (u10);
\path (12.9,0) coordinate (u10p);
\path (13.5,-0.5) coordinate (u101);
\path (13.5,0.5) coordinate (u102);
\draw (u2)--(u1)--(u3)--(u5)--(u7)--(u10);
\draw (u2)--(u4)--(u6)--(u8);
\draw (u101)--(u10)--(u102);
\draw (u01)--(u0)--(u02);
\draw (u0)--(u1);
\draw (u3)--(u4);
\draw (u5)--(u6);
\draw (u7)--(u8);
%
\draw (u2) to[out=90,in=90, distance=1cm] (u8);
\draw (u0) [fill=white] circle (\vr);
\draw (u1) [fill=white] circle (\vr);
\draw (u2) [fill=white] circle (\vr);
\draw (u3) [fill=white] circle (\vr);
\draw (u4) [fill=white] circle (\vr);
\draw (u5) [fill=white] circle (\vr);
\draw (u6) [fill=white] circle (\vr);
\draw (u7) [fill=white] circle (\vr);
\draw (u8) [fill=white] circle (\vr);
\draw (u10) [fill=white] circle (\vr);
\draw[anchor = north] (u0p) node {{\small $w$}};
\draw[anchor = north] (u1) node {{\small $x_6$}};
\draw[anchor = east] (u2) node {{\small $x$}};
\draw[anchor = north] (u3) node {{\small $v_6$}};
\draw[anchor = south] (u4) node {{\small $v_1$}};
\draw[anchor = north] (u5) node {{\small $v_5$}};
\draw[anchor = south] (u6) node {{\small $v_2$}};
\draw[anchor = north] (u7) node {{\small $v_4$}};
\draw[anchor = west] (u8) node {{\small $v_3$}};
\draw[anchor = north] (u10p) node {{\small $x_4$}};
%
\draw (10.5,-1.5) node {{\small (b)}};
%
\path (15.5,0.5) coordinate (u01);
\path (15.5,1.5) coordinate (u02);
\path (16,1) coordinate (u0);
\path (16.1,1) coordinate (u0p);
\path (17,0) coordinate (u1);
\path (17.1,0) coordinate (u1p);
\path (16.5,-0.5) coordinate (w01);
\path (16.5,0.5) coordinate (w02);
\path (17,1) coordinate (u2);
\path (18,0) coordinate (u3);
\path (18,1) coordinate (u4);
\path (19,0) coordinate (u5);
\path (19,1) coordinate (u6);
\path (20,0) coordinate (u7);
\path (20,1) coordinate (u8);
\path (21,0) coordinate (u10);
\path (20.9,0) coordinate (u10p);
\path (21.5,-0.5) coordinate (u101);
\path (21.5,0.5) coordinate (u102);
\draw (u1)--(u3)--(u5)--(u7)--(u10);
\draw (u2)--(u4)--(u6)--(u8);
\draw (u101)--(u10)--(u102);
\draw (u01)--(u0)--(u02);
\draw (w01)--(u1)--(w02);
\draw (u0)--(u2);
\draw (u3)--(u4);
\draw (u5)--(u6);
\draw (u7)--(u8);
%
\draw (u2) to[out=90,in=90, distance=1cm] (u8);
\draw (u0) [fill=white] circle (\vr);
\draw (u1) [fill=white] circle (\vr);
\draw (u2) [fill=white] circle (\vr);
\draw (u3) [fill=white] circle (\vr);
\draw (u4) [fill=white] circle (\vr);
\draw (u5) [fill=white] circle (\vr);
\draw (u6) [fill=white] circle (\vr);
\draw (u7) [fill=white] circle (\vr);
\draw (u8) [fill=white] circle (\vr);
\draw (u10) [fill=white] circle (\vr);
\draw[anchor = north] (u0p) node {{\small $z$}};
\draw[anchor = north] (u1p) node {{\small $x_6$}};
\draw[anchor = north] (u2) node {{\small $x$}};
\draw[anchor = north] (u3) node {{\small $v_6$}};
\draw[anchor = south] (u4) node {{\small $v_1$}};
\draw[anchor = north] (u5) node {{\small $v_5$}};
\draw[anchor = south] (u6) node {{\small $v_2$}};
\draw[anchor = north] (u7) node {{\small $v_4$}};
\draw[anchor = west] (u8) node {{\small $v_3$}};
\draw[anchor = north] (u10p) node {{\small $x_4$}};
%
\draw (18.5,-1.5) node {{\small (c)}};
\end{tikzpicture}
\end{center}
\begin{center}
\vskip -0.75 cm
\caption{Subgraphs in the proof of Claim~\ref{claim.no-domino.1}}
\label{rdom:fig-no-dom}
\end{center}
\end{figure}

\vskip -0.5 cm
Hence, $x$ is adjacent to neither $x_4$ nor $x_6$. In this case, let $z$ be the neighbor of $x$ different from $v_1$ and $v_3$. Hence, the graph illustrated in Figure~\ref{rdom:fig-no-dom}(c) is a subgraph of $G$ (where the edge $x_4x_6$ may or may not exist). We now let $G' = G - (V(F) \cup \{x\})$. Every $\gamma_r$-set of $G'$ can be extended to a RD-set of $G$ by adding to it the set $\{v_1,v_4\}$, and so $\gamma_r(G) \le \gamma_r(G') + 2$, implying that $\w(G) < \w(G') + 20$. Since $G'$ is obtained from $G$ by deleting the edges in a $3$-edge-cut, by Claim~\ref{claim:edge-cut} either no component of $G'$ belongs to $\cB_{\rdom}$ or $G'$ is connected and $G' = R_{10}$. Hence, $\w(G) \ge \w(G') + 22$, a contradiction.~\smallqed

\begin{subclaim}
\label{claim.no-domino.2}
$x_1 \ne x_4$ and $x_3 \ne x_6$.
\end{subclaim}
\proof Suppose, to the contrary, that $x_1 = x_4$ or $x_3 = x_6$. Renaming vertices if necessary, we may assume by symmetry that $x_1 = x_4$. Thus, $x_1$ is a common neighbor of $v_1$ and $v_4$. Let us rename the vertex $x_1$ by~$x$ for notational simplicity. Suppose firstly that $x_3 = x_6$, and so $x_3$ is a common neighbor of $v_3$ and $v_6$. Let us rename the vertex $x_3$ by~$y$ for notational simplicity. If $xy \in E(G)$, then the graph $G$ is determined and $\gamma_r(G) = 3$ and $\w(G) = 32$, a contradiction. Hence, $xy \notin E(G)$. Let $x_1$ and $y_1$ be the neighbors of $x$ and $y$, respectively, that do not belong to $F$.
Suppose that $x_1 = y_1$, and let us rename this common neighbor of $x$ and $y$ by~$w$. Let $z$ be the third neighbor of $w$ different from $x$ and $y$. In this case, let $G'$ be the component of $G - wz$ that contains the vertex~$z$, and so $G'$ is a connected special subcubic graph. Further, the vertex $z$ is the only vertex of degree~$2$ in $G'$, and so $G' \notin \cB_{\rdom}$. Every $\gamma_r$-set of $G'$ can be extended to a RD-set of $G$ by adding to it the set $\{v_1,v_4,y\}$, and so $\gamma_r(G) \le \gamma_r(G') + 3$, implying that $\w(G) < \w(G') + 30$. However, $\w(G) = \w(G') + 35$, a contradiction. Hence, $x_1 \ne y_1$.
We now let $G' = G - (V(F) \cup \{x,y\})$, and so $G'$ is a special subcubic graph that contains exactly two vertices of degree~$2$. By Claim~\ref{claim:edge-cut}, no component of $G'$ belongs to $\cB_{\rdom}$. Every $\gamma_r$-set of $G'$ can be extended to a RD-set of $G$ by adding to it the set $\{x,y,v_2\}$, and so $\gamma_r(G) \le \gamma_r(G') + 3$, implying that $\w(G) < \w(G') + 30$. However, $\w(G) = \w(G') + 30$, a contradiction.

Hence, $x_3 \ne x_6$, that is, the vertices $x, x_3, x_6$ are pairwise distinct. Suppose that $x$ is adjacent to neither $x_3$ nor $x_6$. Let $x'$ be the neighbor of $x$ different from $v_1$ and $v_4$. In this case, we let $G' = G - (V(F) \cup \{x\})$. Let $S'$ be a $\gamma_r$-set of $G'$. If $x' \in S'$, let $S = S' \cup \{v_3,v_6\}$. If $x' \notin S'$, let $S = S' \cup \{v_1,v_4\}$. In both cases, the set $S$ is a RD-set of $G$, and so $\gamma_r(G) \le \gamma_r(G') + 2$, implying that $\w(G) < \w(G') + 20$. Since $G'$ is obtained from $G$ by deleting the edges in a $3$-edge-cut, by Claim~\ref{claim:edge-cut}  either no component of $G'$ belongs to $\cB_{\rdom}$ or $G'$ is connected and $G' = R_{10}$. Hence, $\w(G) \ge \w(G') + 22$, a contradiction. Thus, either $xx_3 \in E(G)$ or $xx_6 \in E(G)$.

Suppose that $xx_3 \in E(G)$. If $x_3x_6 \in E(G)$, then let $y$ be the neighbor of $x_6$ different from $x_3$ and $v_6$, and let $G'$ be the component of $G - x_6y$ that contains the vertex~$y$. Thus, $G' \notin \cB_{\rdom}$. Every $\gamma_r$-set of $G'$ can be extended to a RD-set of $G$ by adding to it the set $\{x_6,v_1,v_4\}$, and so $\gamma_r(G) \le \gamma_r(G') + 3$, implying that $\w(G) < \w(G') + 30$. However, $\w(G) = \w(G') + 35$, a contradiction. Hence, $x_3x_6 \notin E(G)$. In this case, we let $G' = G - (V(F) \cup \{x,x_3\})$. The special subcubic graph $G'$ contains exactly two vertices of degree~$2$, and so by Claim~\ref{claim:edge-cut} no component of $G'$ belongs to $\cB_{\rdom}$. Every $\gamma_r$-set of $G'$ can be extended to a RD-set of $G$ by adding to it the set $\{v_3,v_6,x_3\}$, and so $\gamma_r(G) \le \gamma_r(G') + 3$, implying that $\w(G) < \w(G') + 30$. However, $\w(G) = \w(G') + 30$, a contradiction.

Hence, $xx_3 \notin E(G)$, implying that $xx_6 \in E(G)$. If $x_3x_6 \in E(G)$, then let $y$ be the neighbor of $x_3$ different from $v_3$ and $x_6$, and let $G'$ be the component of $G - x_3y$ that contains the vertex~$y$. We note that $G' \notin \cB_{\rdom}$. Every $\gamma_r$-set of $G'$ can be extended to a RD-set of $G$ by adding to it the set $\{x_3,v_1,v_4\}$, and so $\gamma_r(G) \le \gamma_r(G') + 3$, implying that $\w(G) < \w(G') + 30$. However, $\w(G) = w(G') + 35$, a contradiction. Hence, $x_3x_6 \notin E(G)$. In this case, we let $G' = G - (V(F) \cup \{x,x_6\})$. The special subcubic graph $G'$ contains exactly two vertices of degree~$2$, and so by Claim~\ref{claim:edge-cut} no component of $G'$ belongs to $\cB_{\rdom}$. Every $\gamma_r$-set of $G'$ can be extended to a RD-set of $G$ by adding to it the set $\{v_3,v_6,x_6\}$, and so $\gamma_r(G) \le \gamma_r(G') + 3$, implying that $\w(G) < \w(G') + 30$. However, $\w(G) = \w(G') + 30$, a contradiction.~\smallqed

\medskip
By Claim~\ref{claim.no-domino.1}, $x_1 \ne x_3$ and $x_4 \ne x_6$. By Claim~\ref{claim.no-domino.2}, $x_1 \ne x_4$ and $x_3 \ne x_6$. Thus the vertices $x_1, x_3, x_4, x_6$ are pairwise distinct. Let $G' = G - V(F)$. The graph $G'$ is a special subcubic graph with exactly four vertices of degree~$2$. Every $\gamma_r$-set of $G'$ can be extended to a RD-set of $G$ by adding to it the set $\{v_1,v_4\}$, and so $\gamma_r(G) \le \gamma_r(G') + 2$, implying that $\w(G) < \w(G') + 20$. If no component of $G'$ belongs to $\cB_{\rdom}$, then $w(G) = \w(G') + 20$, a contradiction. Hence, $G'$ contains a component $G_1$ that belongs to $\cB_{\rdom}$. By Claim~\ref{claim:edge-cut}, there is only one such component and $G_1 \in \{R_2,R_4,R_5,R_{10}\}$. Necessarily, $G_1$ contains at least three vertices from the set $\{x_1,x_3,x_4,x_6\}$. In particular, $\{x_1,x_4\} \subset V(G_1)$ or $\{x_3,x_6\} \subset V(G_1)$. Renaming vertices if necessary, we may assume by symmetry that $\{x_3,x_6\} \subset V(G_1)$. Let $S_1$ be a minimum type-$2$ NeRD-set of $G_1$ with respect to the set $X_1 = \{x_3,x_6\}$. We note that $X_1 \cap S_1 = \emptyset$. By Observation~\ref{obser-1}(f), $|S_1| =  \gamma_{r,\dom}(G_1;X_1) \le \gamma_{r}(G_1) - 1$. If $G'$ is connected, then $G' = G_1$ and the set $S_1 \cup \{v_1,v_4\}$ is a RD-set of $G$, and so $\gamma_r(G) \le |S_1| + 2 \le (\gamma_{r}(G') - 1) + 2 = \gamma_{r}(G') + 1$. If $G'$ is disconnected, then let $G_2$ be the component of $G'$ different from $G_1$ which yields $G_1 = R_{10}$. In this case, let $S_2$ be a $\gamma_r$-set of $G_2$ and note that the set $S_1 \cup S_2 \cup \{v_1,v_4\}$ is a RD-set of $G$, implying that $\gamma_r(G) \le |S_1| + |S_2| + 2 \le (\gamma_{r}(G_1) - 1) + \gamma_{r}(G_2) + 2 = \gamma_{r}(G') + 1$. Thus in both cases, $\gamma_r(G) \le \gamma_{r}(G') + 1$, implying that $\w(G) < \w(G') + 10$. However, $\w(G) \ge \w(G') + 17$, a contradiction. This completes the proof of Claim~\ref{claim.no-domino}.~\smallqed

\medskip
By Claim~\ref{claim.no-domino}, the graph $G$ contains no domino as a subgraph.

\begin{claim}
\label{claim.4-cycle-property}
If the graph $G$ contains a $4$-cycle $C$, then the subgraph of $G$ induced by $V(C)$ and all neighbors in $G$ of vertices in $V(C)$ is isomorphic to the corona $C \circ K_1$ of the $4$-cycle $C$.
\end{claim}
\proof Suppose that $G$ contains a $4$-cycle $C \colon v_1v_2v_3v_4v_1v_4$. Since $G$ is triangle-free, the cycle $C$ is an induced cycle. Let $x_i$ be the neighbor of $v_i$ that does not belong to $C$ for $i \in [4]$. Since $G$ has no triangle and no $K_{2,3}$-subgraph, the vertices $x_1, x_2, x_3, x_4$ are pairwise distinct. Let $X = \{x_1,x_2,x_3,x_4\}$. To prove the claim, it suffices to show that the set $X$ is independent. Suppose, to the contrary, that $X$ is not an independent set. Since $G$ contains no domino as a subgraph, $x_ix_{i+1} \notin E(G)$ for all $i \in [4]$, where indices are taken modulo~$4$. Hence, $x_ix_{i+2} \in E(G)$ for some $i \in [4]$, where indices are taken modulo~$4$.
Renaming vertices if necessary, we may assume that $x_1x_3 \in E(G)$. Let $y_1$ be the third neighbor of $x_1$ different from $v_1$ and $x_3$, and let $y_3$ be the third neighbor of $x_3$ different from $v_3$ and $x_1$. Since $G$ is triangle-free, $y_1 \ne y_3$. Further since $G$ contains no domino as a subgraph, $\{y_1,y_3\} \cap \{x_2,x_4\} = \emptyset$. Thus, the vertices $x_2,x_4,y_1,y_3$ are pairwise distinct and the graph illustrated in Figure~\ref{rdom:fig-no-4cycle.1} is a subgraph of $G$. Let $G' = G - (V(C) \cup \{x_1,x_3\}$).

\vskip -0.75 cm
\begin{figure}[htb]
\begin{center}
\begin{tikzpicture}[scale=.8,style=thick,x=0.8cm,y=0.8cm]
\def\vr{2.5pt} 
\path (0.5,0.6) coordinate (u01);
\path (0.5,1.4) coordinate (u02);
\path (1,1) coordinate (u0);
\path (1.1,1) coordinate (u0p);
\path (2,0) coordinate (u1);
\path (2.1,0) coordinate (u1p);
\path (1.5,-0.4) coordinate (w01);
\path (1.5,0.4) coordinate (w02);
\path (2,1) coordinate (u2);
\path (1.9,1) coordinate (u2p);
\path (3,0) coordinate (u3);
\path (3,1) coordinate (u4);
\path (4,0) coordinate (u5);
\path (4,1) coordinate (u6);
\path (5,0) coordinate (u7);
\path (5,1) coordinate (u8);
\path (4.9,1) coordinate (u8p);
\path (5.5,0.6) coordinate (u81);
\path (5.5,1.4) coordinate (u82);
\path (6,0) coordinate (u10);
\path (5.9,0) coordinate (u10p);
\path (6.5,-0.4) coordinate (u101);
\path (6.5,0.4) coordinate (u102);
\draw (u1)--(u3)--(u5)--(u7)--(u10);
\draw (u2)--(u4)--(u6)--(u8);
\draw (u101)--(u10)--(u102);
\draw (u81)--(u8)--(u82);
\draw (u01)--(u0)--(u02);
\draw (w01)--(u1)--(w02);
\draw (u0)--(u2);
\draw (u3)--(u4);
\draw (u5)--(u6);
%
\draw (u2) to[out=45,in=25, distance=2.65cm] (u7);
\draw (u0) [fill=white] circle (\vr);
\draw (u1) [fill=white] circle (\vr);
\draw (u2) [fill=white] circle (\vr);
\draw (u3) [fill=white] circle (\vr);
\draw (u4) [fill=white] circle (\vr);
\draw (u5) [fill=white] circle (\vr);
\draw (u6) [fill=white] circle (\vr);
\draw (u7) [fill=white] circle (\vr);
\draw (u8) [fill=white] circle (\vr);
\draw (u10) [fill=white] circle (\vr);
\draw[anchor = south] (u0p) node {{\small $y_1$}};
\draw[anchor = north] (u1p) node {{\small $x_4$}};
\draw[anchor = south] (u2p) node {{\small $x_1$}};
\draw[anchor = north] (u3) node {{\small $v_4$}};
\draw[anchor = south] (u4) node {{\small $v_1$}};
\draw[anchor = north] (u5) node {{\small $v_3$}};
\draw[anchor = south] (u6) node {{\small $v_2$}};
\draw[anchor = north] (u7) node {{\small $x_3$}};
\draw[anchor = south] (u8p) node {{\small $x_2$}};
\draw[anchor = north] (u10p) node {{\small $y_3$}};
%
%
\end{tikzpicture}
\end{center}
\begin{center}
\vskip -0.65 cm
\caption{A subgraph in the proof of Claim~\ref{claim.4-cycle-property}}
\label{rdom:fig-no-4cycle.1}
\end{center}
\end{figure}
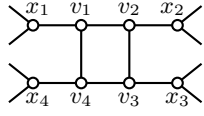

\vskip -0.5 cm
Let $S'$ be a $\gamma_r$-set of $G'$. If $y_1 \in S'$, let $S = S' \cup \{v_3,v_4\}$. If $y_1 \notin S'$, let $S = S' \cup \{v_1,x_3\}$. In both cases, the set $S$ is a RD-set of $G$, and so $\gamma_r(G) \le \gamma_r(G') + 2$, implying that $\w(G) < \w(G') + 20$. If no component of $G'$ belongs to $\cB_{\rdom}$, then $\w(G) = 24 + (\w(G') - 4) = \w(G') + 20$, a contradiction. Hence, $G'$ contains a component $G_1$ that belongs to $\cB_{\rdom}$. By Claim~\ref{claim:edge-cut}, there is only one such component and $G_1 \in \{R_2,R_4,R_5,R_{10}\}$. Necessarily, $G_1$ contains at least three vertices from the set $\{x_2,x_4,y_1,y_3\}$. At least one of $y_1$ and $y_3$ belong to $G_1$. By symmetry, we may assume that $y_1 \in V(G_1)$. Further at least one of $x_2$ and $x_4$ belongs to $G_1$. By symmetry, we may assume that $x_2 \in V(G_1)$.

Let $S_1$ be a minimum type-$2$ NeRD-set of $G_1$ with respect to the set $X_1 = \{y_1,x_2\} \subset V(G_1)$. We note that $X_1 \cap S_1 = \emptyset$. By Observation~\ref{obser-1}(f), $|S_1| =  \gamma_{r,\dom}(G_1;X_1) \le \gamma_{r}(G_1) - 1$. If $G'$ is connected, then $G' = G_1$ and the set $S_1 \cup \{v_1,x_3\}$ is a RD-set of $G$, and so $\gamma_r(G) \le |S_1| + 2 \le (\gamma_{r}(G') - 1) + 2 = \gamma_{r}(G') + 1$. If $G'$ is disconnected, then let $G_2$ be the component of $G'$ different from $G_1$. In this case, let $S_2$ be a $\gamma_r$-set of $G_2$. The set $S_1 \cup S_2 \cup \{v_1,x_3\}$ is a RD-set of $G$, implying that $\gamma_r(G) \le |S_1| + |S_2| + 2 \le (\gamma_{r}(G_1) - 1) + \gamma_{r}(G_2) + 2 = \gamma_{r}(G') + 1$. Thus in both cases, $\gamma_r(G) \le \gamma_{r}(G') + 1$, implying that $\w(G) < \w(G') + 10$. However, $\w(G) \ge \w(G') + 16$, a contradiction. Hence, the set $X$ is an independent set.~\smallqed

\begin{claim}
\label{claim.no-4-cycle-A}
If $G$ contains a $4$-cycle $C \colon v_1v_2v_3v_4v_1v_4$ where $x_i$ is the neighbor of $v_i$ that does not belong to $C$ for $i \in [4]$, then $|N(x_i) \cap N(x_{i+2})| \le 1$ for $i \in [2]$.
\end{claim}
\proof Let the cycle $C$ and the vertices $x_1,x_2,x_3,x_4$ be as in the statement of the claim. By Claim~\ref{claim.4-cycle-property} and our earlier observations, the graph illustrated in Figure~\ref{rdom:fig-no-4cycle.2} is a subgraph of $G$ where $\{x_1,x_2,x_3,x_4\}$ is an independent set. Suppose, to the contrary, that $|N(x_i) \cap N(x_{i+2})| = 2$ for some $i \in [2]$. By symmetry, we may assume that $|N(x_1)\cap N(x_3)|=2$. Let $u$ and $z$ be the two common neighbors of $x_1$ and $x_3$. Since $G$ is triangle-free, the vertices $u$ and $z$ are not adjacent. Let $u'$ and $z'$ be the third neighbors of $u$ and $z$, respectively, different from $x_1$ and $x_3$. Since $G$ has no $K_{2,3}$-subgraph, we note that $u' \ne z'$.

\begin{figure}[htb]
\begin{center}
\begin{tikzpicture}[scale=.8,style=thick,x=0.8cm,y=0.8cm]
\def\vr{2.5pt} 
\path (1.5,0.8) coordinate (u01);
\path (1.5,1.6) coordinate (u02);
\path (2,0) coordinate (u1);
\path (2.1,0) coordinate (u1p);
\path (1.5,-0.4) coordinate (w01);
\path (1.5,0.4) coordinate (w02);
\path (2,1.2) coordinate (u2);
\path (2.1,1.2) coordinate (u2p);
\path (3,0) coordinate (u3);
\path (3,1.2) coordinate (u4);
\path (4,0) coordinate (u5);
\path (4,1.2) coordinate (u6);
\path (5,0) coordinate (u7);
\path (5,1.2) coordinate (u8);
\path (4.9,1.2) coordinate (u8p);
\path (5.5,0.8) coordinate (u81);
\path (5.5,1.6) coordinate (u82);
\path (5.5,-0.4) coordinate (u101);
\path (5.5,0.4) coordinate (u102);
\draw (u1)--(u3)--(u5)--(u7);
\draw (u2)--(u4)--(u6)--(u8);
\draw (u101)--(u7)--(u102);
\draw (u81)--(u8)--(u82);
\draw (u01)--(u2)--(u02);
\draw (w01)--(u1)--(w02);
\draw (u3)--(u4);
\draw (u5)--(u6);
%
%
\draw (u1) [fill=white] circle (\vr);
\draw (u2) [fill=white] circle (\vr);
\draw (u3) [fill=white] circle (\vr);
\draw (u4) [fill=white] circle (\vr);
\draw (u5) [fill=white] circle (\vr);
\draw (u6) [fill=white] circle (\vr);
\draw (u7) [fill=white] circle (\vr);
\draw (u8) [fill=white] circle (\vr);
%
%
\draw[anchor = north] (u1p) node {{\small $x_4$}};
\draw[anchor = south] (u2p) node {{\small $x_1$}};
\draw[anchor = north] (u3) node {{\small $v_4$}};
\draw[anchor = south] (u4) node {{\small $v_1$}};
\draw[anchor = north] (u5) node {{\small $v_3$}};
\draw[anchor = south] (u6) node {{\small $v_2$}};
\draw[anchor = north] (u7) node {{\small $x_3$}};
\draw[anchor = south] (u8p) node {{\small $x_2$}};
%
%
\end{tikzpicture}
\end{center}
\begin{center}
\vskip -0.7 cm
\caption{A subgraph in the proof of Claim~\ref{claim.no-4-cycle-A}}
\label{rdom:fig-no-4cycle.2}
\end{center}
\end{figure}

\begin{subclaim}
\label{claim.no-4-cycle-A.1}
$\{u',z'\} \ne \{x_2,x_4\}$.
\end{subclaim}
\proof Suppose that $\{u',z'\} = \{x_2,x_4\}$. Renaming vertices if necessary, we may assume by symmetry that $u' = x_2$ and $z' = x_4$. Suppose that $x_2$ and $x_4$ have a common neighbor $x$. Let $y$ be the third neighbor of $x$. Let $G'$ be the component of $G - xy$ that contains the vertex~$y$, and so $G'$ is a connected special subcubic graph that contains exactly one vertex of degree~$2$. Every $\gamma_r$-set of $G'$ can be extended to a RD-set of $G$ by adding to it the set $\{x_1,x_4,v_2,v_3\}$, and so $\gamma_r(G) \le \gamma_r(G') + 4$, implying that $\w(G) < \w(G') + 40$. Since $G' \notin \cB_{\rdom}$, we have $\w(G) = w(G') + 43$, a contradiction. Hence, $x_2$ and $x_4$ have no common neighbor. Let $y_2$ be the neighbor of $x_2$ different from $u$ and $v_2$, and let $y_4$ be the neighbor of $x_4$ different from $z$ and $v_4$. By our earlier observations, $y_2 \ne y_4$. If $y_2y_4 \in E(G)$, then let $Q = V(C) \cup \{x_1,x_2,x_3,x_4\} \cup \{u,z,y_2,y_4\}$ and $G' = G - Q$. Thus, $G'$ is a special subcubic graph with exactly two vertices of degree~$2$, implying that no component of $G'$ belongs to $\cB_{\rdom}$. Every $\gamma_r$-set of $G'$ can be extended to a RD-set of $G$ by adding to it the set $\{v_2,v_3,x_1,y_4\}$, and so $\gamma_r(G) \le \gamma_r(G') + 4$, implying that $\w(G) < \w(G') + 40$. However, $\w(G) = \w(G') + 46$, a contradiction.

Hence, $y_2y_4 \notin E(G)$. We now let $Q = V(C) \cup \{x_1,x_2,x_3,x_4\} \cup \{u,z,y_2\}$ and $G' = G - Q$. Thus, $G'$ is a special subcubic graph with exactly three vertices of degree~$2$, implying that either no component of $G'$ belongs to $\cB_{\rdom}$ or $G'$ is connected and $G' = R_{10}$. Every $\gamma_r$-set of $G'$ can be extended to a RD-set of $G$ by adding to it the set $\{v_3,v_4,x_1,y_2\}$, and so $\gamma_r(G) \le \gamma_r(G') + 4$, implying that $\w(G) < \w(G') + 40$. If no component of $G'$ belongs to $\cB_{\rdom}$, then $\w(G) = \w(G') + 41$, a contradiction. Hence, $G' = G_9$. In this case the set $\{v_3,v_4,x_1,y_2\}$ can be extended to a RD-set of $G$ by adding to it $\gamma_r(G') - 1$ vertices from $G'$ applying Observation~\ref{obser-1}(d) with respect to the vertex~$x_4$, and so $\gamma_r(G) \le \gamma_r(G') + 3$, implying that $\w(G) < \w(G') + 30$. However since $G' = G_9$, we have $\w(G) = \w(G') + 38$, a contradiction.~\smallqed

\begin{subclaim}
\label{claim.no-4-cycle-A.2}
$\{u',z'\} \cap \{x_2,x_4\} = \emptyset$.
\end{subclaim}
\proof Suppose that $\{u',z'\} \cap \{x_2,x_4\} \ne \emptyset$, implying by  Claim~\ref{claim.no-4-cycle-A.1} that $|\{u',z'\} \cap \{x_2,x_4\}| = 1$.  Renaming vertices if necessary, we may assume by symmetry that $u' = x_2$. Thus, $z' \ne x_4$. Let $v$ be the neighbor of $u'$ different from $u$ and $v_2$. Since the set $\{x_1,x_2,x_3,x_4\}$ is independent, $v \ne x_4$. Since $G$ contains no $K_{2,3}$ as a subgraph, $u' \ne z'$, that is, $v \ne z$. If $v \ne z'$, then we let $Q = V(C) \cup \{x_1,x_2,x_3,u,z\}$ and let $G' = G - Q$. Thus, $G'$ is a special subcubic graph with exactly three vertices of degree~$2$. Every $\gamma_r$-set of $G'$ can be extended to a RD-set of $G$ by adding to it the set $\{x_1,x_2,v_3\}$, and so $\gamma_r(G) \le \gamma_r(G') + 3$, implying that $\w(G) < \w(G') + 30$. Since either no component of $G'$ belongs to $\cB_{\rdom}$ or $G'$ is connected and $G' = G_9$, we have $\w(G) = \w(G') + 30$, a contradiction. Hence, $v = z'$. If $vx_4 \in E(G)$, then we let $w$ be the neighbor of $x_4$ different from $v$ and $v_4$. Further we let $Q = V(C) \cup \{x_1,x_2,x_3,x_4,u,v,z\}$ and $G' = G - Q$. Thus, $G'$ is a connected special subcubic graph with exactly one vertex of degree~$2$. Every $\gamma_r$-set of $G'$ can be extended to a RD-set of $G$ by adding to it the set $\{x_1,v,v_3\}$, and so $\gamma_r(G) \le \gamma_r(G') + 3$, implying that $\w(G) < \w(G') + 30$. Since $G' \notin \cB_{\rdom}$, we have $\w(G) = \w(G') + 43$, a contradiction. Hence, $vx_4 \notin E(G)$. We now let $Q = V(C) \cup \{x_1,x_2,x_3,u,v,z\}$ and $G' = G - Q$. Thus, $G'$ is a special subcubic graph with exactly two vertices of degree~$2$. Once again every $\gamma_r$-set of $G'$ can be extended to a RD-set of $G$ by adding to it the set $\{x_1,v,v_3\}$, and so $\gamma_r(G) \le \gamma_r(G') + 3$, implying that $\w(G) < \w(G') + 30$. Since no component of $G'$ belongs to $\cB_{\rdom}$, we have $\w(G) = \w(G') + 38$, a contradiction.~\smallqed

\medskip
By Claim~\ref{claim.no-4-cycle-A.2}, $\{u',z'\} \cap \{x_2,x_4\} = \emptyset$. Let $Q = V(C) \cup \{x_1,x_3,u,z\}$ and let $G' = G - Q$. Thus, $G'$ is a special subcubic graph with exactly four vertices of degree~$2$. Every $\gamma_r$-set of $G'$ can be extended to a RD-set of $G$ by adding to it the set $\{x_1,v_3\}$, and so $\gamma_r(G) \le \gamma_r(G') + 2$, implying that $\w(G) < \w(G') + 20$. Since at most one component of $G'$ belongs to $\cB_{\rdom}$, we have $\w(G) = \w(G') + 24$, a contradiction. This completes the proof of Claim~\ref{claim.no-4-cycle-A}.~\smallqed

\begin{claim}
\label{claim.no-4-cycle-B}
If $G$ contains a $4$-cycle $C \colon v_1v_2v_3v_4v_1v_4$ where $x_i$ is the neighbor of $v_i$ that does not belong to $C$ for $i \in [4]$, then $N(x_i) \cap N(x_{i+2}) = \emptyset$ for $i \in [2]$.
\end{claim}
\proof Let the cycle $C$ and the vertices $x_1,x_2,x_3,x_4$ be as in the statement of the claim. Thus the graph illustrated in Figure~\ref{rdom:fig-no-4cycle.2} is a subgraph of $G$ where $\{x_1,x_2,x_3,x_4\}$ is an independent set. Suppose, to the contrary, that $|N(x_i) \cap N(x_{i+2})| \ge 1$ for some $i \in [2]$. By symmetry, we may assume that $|N(x_1) \cap N(x_{3})| \ge 1$. By Claim~\ref{claim.no-4-cycle-A}, $|N(x_1) \cap N(x_{3})| = 1$. Let $z$ be the common neighbor of $x_1$ and $x_3$, and let $z'$ be the third neighbor of $z$.

\begin{subclaim}
\label{claim.no-4-cycle-B.1}
The vertex $z$ is adjacent to neither $x_2$ nor $x_4$.
\end{subclaim}
\proof Suppose, to the contrary, that the vertex $z$ is adjacent to $x_2$ or $x_4$, that is, $z' = x_2$ or $z' = x_4$. By symmetry, we may assume that $z' = x_4$. Let $Q = V(C) \cup \{x_1,x_2,x_3,x_4,z\}$. Let $y_i$ be the neighbor of $x_i$ not in $Q$ for $i \in \{1,3,4\}$. Since the vertex $z$ is the only common neighbor of $x_1$ and $x_3$, we note that $y_1 \ne y_3$.

\begin{subsubclaim}
\label{claim.no-4-cycle-B.1}
The vertices $y_1, y_3, y_4$ are pairwise distinct.
\end{subsubclaim}
\proof Suppose that the vertices $y_1, y_3, y_4$ are not pairwise distinct. By symmetry, we may assume that $y_1 = y_4$. Suppose firstly that $y_1 = y_3 = y_4$. In this case, let $Q' = V(C) \cup \{x_1,x_3,x_4,y_1,z\}$ and let $G' = G - Q'$. Thus, $G'$ is a connected special subcubic graph with exactly small vertex, and so $G' \notin \cB_{\rdom}$ and $\w(G) \ge \w(G') + 35$. However, $\gamma_r(G) \le \gamma_r(G') + 3$, implying that $\w(G) < \w(G') + 30$, a contradiction. Hence, the vertices $y_1$ and $y_3$ are distinct.

Let $z_1$ be the neighbor of $y_1$ different from $x_1$ and $x_4$. Suppose that $x_2,y_3,z_1$ are pairwise distinct. In this case, we let $Q' = V(C) \cup \{x_1,x_3,x_4,y_1,z\}$  and $G' = G - Q'$. Thus, $G'$ is a special subcubic graph with three small vertices, and so by Claim~\ref{claim:edge-cut} either no component of $G'$ belongs to $\cB_{\rdom}$ or $G'$ is connected and $G' = R_{10}$. Hence, $\w(G) \ge \w(G') + 30$. Moreover, $\gamma_r(G) \le \gamma_r(G') + 3$, implying that $\w(G) < \w(G') + 30$, a contradiction. Hence, $x_2,y_3,z_1$ are not pairwise distinct vertices. Since $x_2$ and $x_3$ are not adjacent, $x_2 \ne y_3$. Hence either $z_1 = x_2$ or $z_1 = y_3$.

Suppose that $z_1 = x_2$. In this case, let $y_2$ be the neighbor of $x_2$ different from $v_2$ and $y_1$. If $y_2 \ne y_3$, then let $Q' = V(C) \cup \{x_1,x_2,x_3,x_4,y_1,z\}$  and $G' = G - Q'$. Thus, $G'$ is a special subcubic graph with two small vertices, and so no component of $G'$ belongs to $\cB_{\rdom}$, whence $\w(G) = \w(G') + 38$. However, $\gamma_r(G) \le \gamma_r(G') + 3$, implying that $\w(G) < \w(G') + 30$, a contradiction. If $y_2 = y_3$, then in this case let $Q' = V(C) \cup \{x_1,x_2,x_3,x_4,y_1,y_2,z\}$  and $G' = G - Q'$. Thus, $G'$ is a connected special subcubic graph with one small vertex, and so $G' \notin\cB_{\rdom}$ and $\w(G) = \w(G') + 43$. However, $\gamma_r(G) \le \gamma_r(G') + 4$, implying that $\w(G) < \w(G') + 40$, a contradiction.

Hence, $z_1 = y_3$. Suppose that $x_2y_3 \in E(G)$. In this case, we let $Q' = V(C) \cup \{x_1,x_2,x_3,x_4,y_1,y_3,z\}$  and let $G' = G - Q'$. Thus, $G'$ is a connected special subcubic graph with one small vertex, and so $G' \notin\cB_{\rdom}$ and $\w(G) = \w(G') + 43$. However, $\gamma_r(G) \le \gamma_r(G') + 4$, implying that $\w(G) < \w(G') + 40$, a contradiction. Hence, $x_2y_3 \notin E(G)$. In this case, we let $Q' = V(C) \cup \{x_1,x_3,x_4,y_1,y_3,z\}$  and let $G' = G - Q'$. Thus, $G'$ is a special subcubic graph with two small vertices, and so no component of $G'$ belongs to $\cB_{\rdom}$, yielding $\w(G) = \w(G') + 38$. However, $\gamma_r(G) \le \gamma_r(G') + 3$, a contradiction.~\smallqed

\begin{subsubclaim}
\label{claim.no-4-cycle-B.2}
The vertex $x_2$ is adjacent to at most one of $y_1$ and $y_4$.
\end{subsubclaim}
\proof
Suppose that $x_2$ is adjacent to both $y_1$ and $y_4$. In this case, we let $Q' = V(C) \cup \{x_1,x_2,x_3,x_4,y_1,y_4,z\}$  and let $G' = G - Q'$. Suppose that $G'$ is a special subcubic graph, and so $G'$ contains exactly three small vertices. By Claim~\ref{claim:edge-cut} either no component of $G'$ belongs to $\cB_{\rdom}$ or $G'$ is connected and $G' = R_{10}$. If no component of $G'$ belongs to $\cB_{\rdom}$, then $\w(G) = \w(G') + 41$, while if $G' = R_{10}$, then $\w(G) = \w(G') + 38$. However, $\gamma_r(G) \le \gamma_r(G') + 3$, implying that $\w(G) < \w(G') + 30$, a contradiction. Hence, $G'$ is not a special subcubic graph.

Let $z_1$ and $z_4$ be the neighbors of $y_1$ and $y_4$, respectively, in $G$ that do not belong to $Q$. Since $G'$ is not a special subcubic graph, the vertices $y_3,z_1,z_4$ are not pairwise distinct. If $y_3$ is adjacent to both $y_1$ and $y_4$, then the graph $G$ is determined and $\gamma_r(G) = 3$ and $\w(G) = 48$, a contradiction. If $y_3$ is adjacent to exactly one of $y_1$ and $y_4$, then by symmetry we may assume that $y_3y_4$ is an edge. In this case, we let $Q' = V(C) \cup \{x_1,x_2,x_3,x_4,y_1,y_3,y_4,z\}$  and let $G' = G - Q'$. Thus, $G'$ is a special subcubic graph with two small vertices, and so no component of $G'$ belongs to $\cB_{\rdom}$ and  $\w(G) = \w(G') + 46$. However, $\gamma_r(G) \le \gamma_r(G') + 3$, implying that $\w(G) < \w(G') + 30$, a contradiction. Hence, $y_3$ is adjacent to neither $y_1$ nor $y_4$, that is, $y_3 \ne z_1$ and $y_3 \ne z_4$, implying that $z_1 = z_4$.

If $z_1y_3 \in E(G)$, then we let $Q' = Q \cup \{y_1,y_3,y_4,z_1\}$ and let $G' = G - Q'$. Thus, $G'$ is a special subcubic graph with one small vertex, and so $G' \notin\cB_{\rdom}$ and $\w(G) = \w(G') + 51$. However, $\gamma_r(G) \le \gamma_r(G') + 4$, implying that $\w(G) < \w(G') + 40$, a contradiction. If $z_1y_3 \notin E(G)$, then let $z'$ be the neighbor of $z_1$ different from $y_1$ and $y_4$, and in this case let $Q' = Q \cup \{y_1,y_4,z_1\}$ and $G' = G - Q'$. Thus, $G'$ is a special subcubic graph with two small vertices, and so no component of $G'$ belongs to $\cB_{\rdom}$ and $\w(G) = \w(G') + 46$. However, $\gamma_r(G) \le \gamma_r(G') + 3$, implying that $\w(G) < \w(G') + 30$, a contradiction.~\smallqed

\medskip
By Claim~\ref{claim.no-4-cycle-B.2}, the vertex $x_2$ is adjacent to at most one of $y_1$ and $y_4$. By symmetry, we may assume that $x_2y_1 \notin E(G)$. Let $Q' = Q \setminus \{x_2\}$ and let $G'$ be obtained from $G - Q'$ by adding the edge $x_2y_1$. Thus, $G'$ is a subcubic graph with exactly two small vertices, namely $y_3$ and $y_4$, and so no component of $G'$ belongs to $\cB_{\rdom}$ and $\w(G) = \w(G') + 30$. Let $S'$ be a $\gamma_r$-set of $G'$. If $x_2 \in S'$, let $S = S' \cup \{v_3,v_4,x_1\}$.  If $x_2 \notin S'$ and $y_1 \in S'$, let $S = S' \cup \{v_2,v_3,x_4\}$. If $x_2 \notin S'$ and $y_1 \notin S'$, let $S = S' \cup \{v_1,x_3,x_4\}$. In all cases, $S$ is a RD-set of $G$, and so $\gamma_r(G) \le \gamma_r(G') + 3$, implying that $\w(G) < \w(G') + 30$, a contradiction.  This completes the proof of Claim~\ref{claim.no-4-cycle-B}.~\smallqed

\begin{claim}
\label{claim.no-4-cycle}
The graph $G$ has no $4$-cycle.
\end{claim}
\proof
Suppose, to the contrary, that $G$ contains a $4$-cycle $C \colon v_1v_2v_3v_4v_1v_4$. Let $x_i$ be the neighbor of $v_i$ that does not belong to $C$ for $i \in [4]$. Thus the graph illustrated in Figure~\ref{rdom:fig-no-4cycle.2} is a subgraph of $G$ where $\{x_1,x_2,x_3,x_4\}$ is an independent set.
By Claim~\ref{claim.no-4-cycle-B}, $N(x_i) \cap N(x_{i+2}) = \emptyset$ for $i \in [2]$. Thus, $x_1$ and $x_3$ have no common neighbor, and $x_2$ and $x_4$ have no common neighbor. Let $y_i$ and $z_i$ be the two neighbors of $x_i$ different from $v_i$ for $i \in \{1,3\}$. By our earlier observations, the vertices $x_2,x_4,y_1,y_3,z_1,z_3$ are pairwise distinct.

If $x_2$ is adjacent to both $y_3$ and $z_3$, then $C' \colon x_2 y_3 x_3 z_3 x_2$ is a $4$-cycle. However, in this case the neighbors $v_2$ and $v_3$ of the vertices $x_2$ and $x_3$, respectively, that do not belong to the cycle $C'$ are adjacent, contradicting Claim~\ref{claim.4-cycle-property}. Hence, the vertex $x_2$ is not adjacent to at least one of $y_3$ and $z_3$. Renaming vertices, if necessary, we may assume that $x_2$ is not adjacent to $y_3$. Let $Q = V(C) \cup \{x_1,x_3\}$ and let $G'$ be obtained from $G - Q$ by adding the edge $x_2y_3$. The resulting graph $G'$ is a special subcubic that contains exactly four small vertices, namely $x_4,y_1,z_1,z_3$. Thus at most one component of $G'$ belongs to $\cB_{\rdom}$.

Let $S'$ be a $\gamma_r$-set of $G'$, and let $S$ be the set defined as follows. If $y_3 \in S'$ and $y_1 \notin S'$, let $S = S' \cup \{v_1,v_2\}$. If $y_3 \in S'$ and $y_1 \in S'$, let $S = S' \cup \{v_2,v_4\}$. If $y_3 \notin S'$, $x_2 \in S'$ and $y_1 \notin S'$, let $S = S' \cup \{v_1,x_3\}$. If $y_3 \notin S'$, $x_2 \in S'$ and $y_1 \in S'$, let $S = S' \cup \{v_4,x_3\}$. If $x_2 \notin S'$ and $y_3 \notin S'$, let $S = S' \cup \{x_1,v_3\}$.
The resulting set $S$ is a RD-set of $G$, and so $\gamma_r(G) \le \gamma_r(G') + 2$, implying that $\w(G) < \w(G') + 20$. If $G'$ has no component in $\cB_{\rdom}$, then $\w(G) = 24 + (w(G') - 4) = \w(G') + 20$, a contradiction. Hence, $G'$ contains a component $G_1$ that belongs to $\cB_{\rdom}$. By Claim~\ref{claim:edge-cut}, there is only one such component and $G_1 \in \{R_2,R_4,R_5,R_{10}\}$.

Suppose that the added edge $x_2y_3$ belongs to $G_1$. In this case, $G_1$ contains two adjacent vertices of degree~$3$, and so $G_1 \in \{R_4,R_5,R_{10}\}$. Let $G_1^*$ be the graph obtained from $G_1$ by subdividing the edge $x_2y_3$ three times resulting in the path $x_2 v_2 v_3 x_3 y_3$. Let $S_1^*$ be a minimum type-$2$ NeRD-set of $G_1^*$ with respect to the vertex $v_3$. By Observation~\ref{obser-4}(b), $|S_1^*| = \gamma_{r,\dom}(G_1^*;v_3) \le \gamma_{r}(G_1)$. By our earlier observations, at least one of $y_1$ and $z_1$ belong to the graph $G_1$. Renaming vertices if necessary, we may assume that $y_1 \in V(G_1)$. If $y_1 \in S_1^*$, then let $S = S_1^* \cup \{v_3\}$.  If $y_1 \notin S_1^*$, then let $S = S_1^* \cup \{v_1\}$. If $G' = G_1$, then the set $S$ is a RD-set of $G$, and so $\gamma_r(G) \le |S_1^*| + 1 \le \gamma_{r}(G_1) + 1 = \gamma_{r}(G') + 1$. If $G' \ne G_1$, then $G_1 = R_{10}$. In this case, $G_1$ contains three vertices of degree~$2$ in $G'$, and so $G'$ is disconnected and contains a second component $G_2$. Since $G_2$ contains exactly one small vertex, the component $G_2$ does not belong to $\cB_{\rdom}$. Every $\gamma_r$-set of $G_2$ can be extended to a RD-set of $G$ by adding to it the set $S$, and so in this case $\gamma_r(G) \le |S| + |S_2| = |S_1^*| + 1 + |S_2| \le (\gamma_r(G_1) + 1) + \gamma_r(G_2) = \gamma_r(G') + 1$. In both cases, $\gamma_r(G) \le \gamma_r(G') + 1$, implying that $\w(G) < \w(G') + 10$. However, $\w(G) \ge \w(G') + 16$, a contradiction.

Hence, the added edge $x_2y_3$ does not belong to $G_1$, implying that $G'$ is disconnected. Let $G_2 = G' - V(G_1)$. We note that $G_2$ contains the added edge $x_2y_3$ and contains at most two components and contains at most one small vertex. Thus, no component of $G_2$ belongs to $\cB_{\rdom}$. Let $S_2$ be a $\gamma_r$-set of $G_2$. We note that $\gamma_r(G') = \gamma_r(G_1) + \gamma_r(G_2)$. Recall that $G_1$ contains at least three vertices from the set $\{x_4,y_1,z_1,z_3\}$.

Suppose that $y_3 \in S_2$ or $x_2 \in S_2$. In this case, we let $X_1 \subset V(G_1)$ such that $|X_1| = 2$ and $X_1 \subset \{x_4,y_1,z_1\}$ noting that at least two vertices in $\{x_4,y_1,z_1\}$ belong to the component $G_1$. Let $S_1$ be a minimum type-$2$ NeRD-set of $G_1$ with respect to the vertex $X_1$. By Observation~\ref{obser-1}(f), $|S_1| = \gamma_{r,\dom}(G_1;X_1) \le \gamma_r(G_1) - 1$. If $y_3 \in S_2$, let $S^* = S_1 \cup S_2 \cup \{v_1,v_2\}$. If $y_3 \notin S_2$ and $x_2 \in S_2$, let $S^* = S_1 \cup S_2 \cup \{v_1,x_3\}$.
Suppose that $y_3 \notin S_2$ and $x_2 \notin S_2$. At least one of $y_1$ and $z_1$ belong to the graph $G_1$. Renaming vertices if necessary, we may assume that $y_1 \in V(G_1)$. In this case, we let $S_1$ be a minimum type-$1$ NeRD-set of $G_1$ with respect to the vertex $y_1$. By Observation~\ref{obser-1}(d), $|S_1| = \gamma_{r,\ndom}(G_1;y_1) \le \gamma_r(G_1) - 1$. Let $S^* = S_1 \cup S_2 \cup \{x_1,v_3\}$.
In all cases, $|S_1| \le \gamma_r(G_1) - 1$ and the set $S^*$ is a RD-set of $G$, and so $\gamma_r(G) \le |S_1| + |S_2| + 2 \le (\gamma_r(G_1) - 1) + \gamma_r(G_2) + 2 = \gamma_r(G') + 1$, implying that $\w(G) < \w(G') + 10$. However, $\w(G) \ge \w(G') + 16$, a contradiction. This completes the proof of Claim~\ref{claim.no-4-cycle}.~\smallqed

\begin{claim}
\label{claim.no-5-cycle}
The graph $G$ has no $5$-cycle.
\end{claim}
\proof
Suppose to the contrary that $G$ contains a $5$-cycle $C \colon v_1 v_2 v_3 v_4 v_5 v_1$. Let $x_i$ be the neighbor of $v_i$ that does not belong to $C$ for $i\in [5]$. Let $X = \{x_1,\ldots,x_5\}$ and let $Q = V(C) \cup X$.

\begin{subclaim}
\label{claim.no-5-cycle-xi-neighbor-xj}
Each vertex $x \in X$ has at least one neighbor in $X$.
\end{subclaim}
\proof Suppose, to the contrary, that there is a vertex in $X$ with no neighbor in $X$. Renaming vertices if necessary that $x_1$ has no neighbor in $X$. Let $y_1$ and $z_1$ be the two neighbors of $x_1$ different from $v_1$. If $x_2$ is adjacent to both $y_1$ and $z_1$, then $x_1 y_1 x_2 z_1 x_1$ is a $4$-cycle in $G$, a contradiction. Hence we may assume that $x_2z_1 \notin E(G)$. Let $Q' = V(C) \cup \{x_1\}$ and let $G'$ be obtained from $G - Q'$ by adding the edge $e = x_2z_1$. Thus, $G'$ is a special subcubic graph that contains exactly four small vertices. By Claim~\ref{claim:edge-cut}, at most one component of $G'$ belongs to $\cB_{\rdom}$. Let $S'$ be a $\gamma_r$-set of $G'$. If $x_2 \in S'$, let $S = S' \cup \{x_1,v_4\}$. If $x_2 \notin S'$ and $z_1 \in S'$, let $S = S' \cup \{v_2,v_5\}$. If $x_2 \notin S'$ and $z_1 \notin S'$, let $S = S' \cup \{v_1,v_3\}$. In all cases, $S$ is a RD-set of $G$, and so $\gamma_r(G) \le |S| + 2 = \gamma_r(G') + 2$, implying that $\w(G) < \w(G') + 20$.
If no component of $G'$ belongs to~$\cB_{\rdom}$, then $\w(G) = \w(G') + 20$, a contradiction. Hence, $G'$ contains a component $G_1$ that belongs to $\cB_{\rdom}$. By Claim~\ref{claim:edge-cut}, there is only one such component and $G_1 \in \{R_2,R_4,R_5,R_{10}\}$. If $G_1 = R_2$, then since $G$ contains no $4$-cycle the added edge $e$ belongs to $G_1$, implying that $G_1$ contains two adjacent vertices of degree~$3$, a contradiction. Hence, $G_1 \in \{R_4,R_5,R_{10}\}$.

Suppose that $e \in E(G_1)$. If $G_1 \in \{R_4,R_5\}$, then the graph $G$ is determined (in the sense that $V(G) = V(C) \cup V(G_1)$) and $\gamma_r(G) \le 4$ and $\w(G) = 56$, a contradiction. Hence, $G_1 = R_{10}$. In this case, $G'$ is disconnected and contains two components. Let $G_2$ be the second component of $G'$, and so $G_2$ contains exactly one small vertex and $G_2 \notin \cB_{\rdom}$. Let $G^*$ be the subgraph of $G$ or order~$17$ induced by $V(C) \cup V(G_1)$.  Every $\gamma_r$-set of $G^*$ can be extended to a RD-set of $G$ by adding to it a $\gamma_r$-set of $G_2$, and so $\gamma_r(G) \le \gamma_r(G_2) + \gamma_r(G^*) \le \gamma_r(G_2) + 6$, implying that $\w(G) < \w(G_2) + 60$. However, $\w(G) = 17 \times 4 + (\w(G_2) - 1) = \w(G_2) + 67$, a contradiction. Hence, $e \notin E(G_1)$. Let $G_2 = G' - V(G_1)$, and so $e \in E(G_2)$ and $\gamma_r(G') = \gamma_r(G_1) + \gamma_r(G_2)$. Since $G_1$ contains at least three small vertices, the graph $G_2$ contains at most one small vertex. Further $G_2$ has at most two components, and so no component of $G_2$ belongs to $\cB_{\rdom}$. Let $S_2$ be a $\gamma_r$-set of $G_2$. We now define a RD-set $S$ in $G$ as follows.

Suppose that $x_2 \in S_2$. We note that at least one of $x_4$ and $y_1$ belongs to $G_1$. Let $v \in \{x_4,y_1\} \cap V(G_1)$. Let $S_1$ be a minimum type-$1$ NeRD-set of $G_1$ with respect to the vertex $v$. By Observation~\ref{obser-1}(d), $|S_1| = \gamma_{r,\ndom}(G_1;v) \le \gamma_{r}(G_1) - 1$. Let $S = S_1 \cup S_2 \cup \{x_1,v_4\}$.

Suppose that $x_2 \notin S_2$ and $z_1 \in S_2$. We note that at least one of $x_3$ and $x_4$ belongs to $G_1$. Let $v \in \{x_3,x_4\} \cap V(G_1)$. Let $S_1$ be a minimum type-$2$ NeRD-set of $G_1$ with respect to the vertex $v$. By Observation~\ref{obser-1}(e), $|S_1| = \gamma_{r,\dom}(G_1;v) \le \gamma_{r}(G_1) - 1$. Let $S = S_1 \cup S_2 \cup \{v_2,v_5\}$.

Suppose that $x_2 \notin S_2$ and $z_1 \notin S_2$. We note that at least one of $x_4$ and $x_5$ belongs to $G_1$. Let $v \in \{x_4,x_5\} \cap V(G_1)$. Let $S_1$ be a minimum type-$2$ NeRD-set of $G_1$ with respect to the vertex $v$. By Observation~\ref{obser-1}(e), $|S_1| = \gamma_{r,\dom}(G_1;v) \le \gamma_{r}(G_1) - 1$. Let $S = S_1 \cup S_2 \cup \{v_1,v_3\}$.

In all cases, $|S_1| \le \gamma_{r}(G_1) - 1$ and the set $S$ is a RD-set of $G$. Therefore, $\gamma_r(G) \le |S_1| + |S_2| + 2 \le (\gamma_{r}(G_1) - 1) + \gamma_r(G_2) + 2 = \gamma_r(G') + 1$, implying that $\w(G) < \w(G') + 10$. However, $\w(G) \ge \w(G') + 16$, a contradiction.~\smallqed

\medskip
By Claim~\ref{claim.no-5-cycle-xi-neighbor-xj}, each vertex $x \in X$ has at least one neighbor in $X$. Hence, $G[X]$ contains at least three edges. Since $G$ has no $4$-cycles, we infer that $G[X]$ contains at most five edges. Using symmetry, we may assume without loss of generality that $\{x_1x_4,x_2x_4,x_3x_5\} \subset E(G)$ noting that $G$ contains no $4$-cycles.

\begin{subclaim}
\label{claim.no-5-cycle-x1x3-x2x5-edge}
$G[X]$ contains at least four edges.
\end{subclaim}
\proof
Suppose, to the contrary, that $G[X]$ contains exactly three edges. Let $y_i$ be neighbor of $x_i$ not in $Q$ for $i \in \{1,2,3,5\}$. Since $G$ has no $4$-cycles, we note that $y_1 \ne y_2$, and since $G$ has no triangles, we note that $y_3 \ne y_5$. We show firstly that $y_2 \ne y_3$. Suppose, to the contrary, that $y_2=y_3$, and let $z$ be the third neighbor of $y_2$.

Suppose that $y_1=y_5$. Let $z'$ be the neighbor of $y_1$ different from $x_1$ and $x_5$. Suppose that $z = z'$. In this case, let $Q' = Q \cup \{y_1,y_2,z\}$ and let $G' = G - Q'$. Thus, $G'$ is a connected special subcubic graph that contains exactly one small vertex, and so $G' \notin \cB_{\rdom}$ and $\w(G) = \w(G') + 51$. However, $\gamma_r(G) \le \gamma_r(G') + 5$, implying that $\w(G) < \w(G') + 50$, a contradiction. Hence, $z \ne z'$. In this case, let $Q' = Q \cup \{y_1,y_2\}$ and let $G' = G - Q'$. Thus, $G'$ is a special subcubic graph that contains exactly two small vertices, and so no component of $G'$ belongs to $\cB_{\rdom}$ and $\w(G) = \w(G') + 46$. However, $\gamma_r(G) \le \gamma_r(G') + 4$, a contradiction.

Hence, $y_1 \ne y_5$. Since $G$ contain no $4$-cycle, $y_2y_5\notin E(G)$. Suppose that $y_1y_2 \notin E(G)$, and so $z$ is distinct from both $y_1$ and $y_5$. In this case, let $Q' = Q \cup \{y_2\}$ and let $G' = G - Q'$. Thus, $G'$ is a special subcubic graph that contains exactly three small vertices, and so either no component of $G'$ belongs to  $\cB_{\rdom}$ or $G'$ is connected and $G' = R_{10}$. Therefore, $\w(G) \ge \w(G') + 38$. However, $\gamma_r(G) \le \gamma_r(G') + 3$, implying that $\w(G) < \w(G') + 30$, a contradiction.

Hence, $y_1y_2\in E(G)$. If, in addition, $y_1y_5\in E(G)$, then let $Q' = Q \cup \{y_1,y_2,y_5\}$ and let $G' = G - Q'$. Thus, $G'$ is a connected special subcubic graph that contains exactly one small vertex and so $G'\notin \cB_{\rdom}$ and $\w(G) \ge \w(G') + 51$. However, $\gamma_r(G) \le \gamma_r(G') + 4$, implying that $\w(G) < \w(G') + 40$, a contradiction. On the other hand, if $y_1y_2\in E(G)$ and $y_1y_5\notin E(G)$, then let $Q' = Q \cup \{y_1,y_2\}$ and let $G' = G - Q'$.   Thus, $G'$ is a special subcubic graph that contains two small vertices and so no component of $G'$ is in $\cB_{\rdom}$ and $\w(G) \ge \w(G') + 46$. However, $\gamma_r(G) \le \gamma_r(G') + 4$, implying that $\w(G) < \w(G') + 40$, a contradiction. Hence, $y_2 \ne y_3$.

We show next that $y_2 \ne y_5$. Suppose, to the contrary, that $y_2=y_5$, and let $z$ be the third neighbor of $y_2$. Suppose that $y_1=y_3$. Let $z'$ be the neighbor of $y_1$ different from $x_1$ and $x_3$. Suppose that $z = z'$. In this case, let $Q' = Q \cup \{y_1,y_2,z\}$ and let $G' = G - Q'$. Thus, $G'$ is a connected special subcubic graph that contains exactly one small vertex, and so $G' \notin \cB_{\rdom}$ and $\w(G) = \w(G') + 51$. However, $\gamma_r(G) \le \gamma_r(G') + 5$, implying that $\w(G) < \w(G') + 50$, a contradiction. Hence, $z \ne z'$. In this case, let $Q' = Q \cup \{y_1,y_2\}$ and let $G' = G - Q'$. Thus, $G'$ is a special subcubic graph that contains exactly two small vertices, and so no component of $G'$ belongs to $\cB_{\rdom}$ and $\w(G) = \w(G') + 46$. However, $\gamma_r(G) \le \gamma_r(G') + 4$, a contradiction.

Hence, $y_1 \ne y_3$. Since $G$ contains no $4$-cycle, vertex $y_2$ is not adjacent to $y_3$, that is, $z \ne y_3$. Suppose that $z \ne y_1$. In this case, let $Q' = Q \cup \{y_2\}$ and let $G' = G - Q'$. Thus, $G'$ is a special subcubic graph that contains exactly three small vertices, and either no component of $G'$ belongs to $\cB_{\rdom}$ or $G'$ is connected and $G' = R_{10}$. Thus, $\w(G) = \w(G') + 38$. However, $\gamma_r(G) \le \gamma_r(G') + 3$, a contradiction.

Hence, $z=y_1$, that is, $y_1y_2\in E(G)$. If, in addition, $y_1y_3\in E(G)$, then let $Q' = Q \cup \{y_1,y_2,y_3\}$ and let $G' = G - Q'$. Thus, $G'$ is a connected special subcubic graph that contains exactly one small vertex and so $G'\notin \cB_{\rdom}$ and $\w(G) \ge \w(G') + 51$. However, $\gamma_r(G) \le \gamma_r(G') + 4$, implying that $\w(G) < \w(G') + 40$, a contradiction. On the other hand, if $y_1y_2\in E(G)$ and $y_1y_3\notin E(G)$, then let $Q' = Q \cup \{y_1,y_2\}$ and let $G' = G - Q'$. Thus, $G'$ is a special subcubic graph that contains two small vertices and so no component of $G'$ is in $\cB_{\rdom}$ and $\w(G) \ge \w(G') + 46$. However, $\gamma_r(G) \le \gamma_r(G') + 4$, implying that $\w(G) < \w(G') + 40$, a contradiction. Hence, $y_2 \ne y_5$.

By our earlier observations, the vertices $y_1,y_2,y_3$ and $y_5$ are pairwise distinct. Let $G' = G - Q$. Thus, $G'$ is a  special subcubic graph that contains exactly four small vertices. At most one component of $G'$ belongs to $\cB_{\rdom}$ and $\w(G) \ge \w(G') + 32$. However, $\gamma_r(G) \le \gamma_r(G') + 3$, implying that $\w(G) < \w(G') + 30$, a contradiction. This completes the proof of Claim~\ref{claim.no-5-cycle-x1x3-x2x5-edge}.~\smallqed

\medskip
By Claim~\ref{claim.no-5-cycle-x1x3-x2x5-edge}, the graph $G[X]$ contains at least four edges. Hence at least one of $x_1x_3$ and $x_2x_5$ is an edge. By symmetry, we may assume that $x_1x_3\in E(G)$. If $x_2x_5 \in E(G)$, then the graph $G$ is determined and is isomorphic to the Petersen graph shown in Figure~\ref{Petersen}. In this case, $\gamma_r(G) = 4$ and $\w(G) = 40$, a contradiction. Hence, $x_2x_5\notin E(G)$. Let $y_i$ be the neighbor of $x_i$ not in $Q$ for $i \in \{2,5\}$. Suppose that $y_2 = y_5$. In this case, let $Q' = Q \cup \{y_2\}$ and let $G' = G - Q'$. Thus, $G'$ is a connected special subcubic graph that contains exactly one small vertex, and so $G' \notin \cB_{\rdom}$ and $\w(G) = \w(G') + 41$. However, $\gamma_r(G) \le \gamma_r(G') + 3$, implying that $\w(G) < \w(G') + 30$, a contradiction. Hence, $y_2 \ne y_5$. We now let $G' = G - Q$. Thus, $G'$ is a special subcubic graph that contains exactly two small vertices, and so no component of $G'$ belongs to $\cB_{\rdom}$ and $\w(G) \ge \w(G') + 38$. However, $\gamma_r(G) \le \gamma_r(G') + 3$, implying that $\w(G) < \w(G') + 30$, a contradiction. This completes the proof of Claim~\ref{claim.no-5-cycle}.~\smallqed

\begin{claim}
\label{claim.no-6-cycle}
The graph $G$ has no $6$-cycle.
\end{claim}
\proof
Suppose, to the contrary, that $G$ contains a $6$-cycle $C \colon v_1v_2v_3v_4v_5v_6v_1$. Thus, $G$ has girth equal to~$6$. In particular, $C$ is an induced cycle in $G$. Let $x_i$ be the neighbor of $v_i$ that does not belong to $C$ for $i\in [6]$. The girth condition implies that $x_i \ne x_j$ for $1 \le i < j \le 6$. Let $X = \{x_1,\ldots,x_6\}$. The girth condition implies that the only possible edges in $G[X]$ are the edges $x_1x_4$, $x_2x_5$ and $x_3x_6$. Let $G'$ be the special subcubic graph obtained from $G - V(C)$ by adding the edge $x_1x_2$. Thus, $G'$ contains exactly four small vertices, namely $x_3, x_4, x_5, x_6$. By Claim~\ref{claim:edge-cut}, at most one component of $G'$ belongs to $\cB_{\rdom}$. Let $S'$ be a $\gamma_r$-set of $G'$. If $x_1 \in S'$, let $S = S' \cup \{v_2,v_5\}$. If $x_1 \notin S'$ and $x_2 \in S'$, let $S = S' \cup \{v_1,v_4\}$. If $x_1 \notin S'$ and $x_2 \notin S'$, let $S = S' \cup \{v_3,v_6\}$. In all cases, $S$ is a RD-set of $G$, and so $\gamma_r(G) \le |S'| + 2 = \gamma_r(G') + 2$, implying that $\w(G) < \w(G') + 20$.

If no component of $G'$ belongs to~$\cB_{\rdom}$, then $\w(G) = \w(G') + 20$, a contradiction. Hence, $G'$ contains a component $G_1$ that belongs to $\cB_{\rdom}$. By Claim~\ref{claim:edge-cut}, there is only one such component and $G_1 \in \{R_2,R_4,R_5,R_{10}\}$. The set $X_1 = \{x_3,x_4,x_5,x_6\}$ of small vertices in $G_1$ is either independent or induces a graph that contains exactly one edge, namely the edge $x_3x_6$. Further, every cycle of length less than~$6$ in $G_1$ must contain the added edge $x_1x_2$ since the graph $G$ contains no cycles of length~$3$,~$4$ or~$5$. If $G_1$ contains the edge $x_1x_2$, then $G_1$ contains two adjacent vertices of degree~$3$. From these properties of the graph $G'$ we infer that $G_1 \notin \{R_2,R_{10}\}$. Since $R_5$ contains two pairs of small vertices that are adjacent while the set $X_1$ contains at most one pair of small vertices that are adjacent, $G_1 \ne R_5$, implying that $G_1 = R_4$ and $X_1 \subset V(G_1)$. The structure of $R_4$ implies that in this case every small vertex in $R_4$ is at distance~$2$ from two other small vertices. In particular, the vertex $x_4$ is at distance~$2$ from at least one of $x_3$ or $x_5$ in $G'$. If $x_3$ and $x_4$ are at distance~$2$ in $G'$ and $w$ denotes their common neighbor in $G'$, then $x_4v_4v_3x_3wx_4$ is a $5$-cycle in $G$. If $x_4$ and $x_5$ are at distance~$2$ in $G'$ and $z$ denotes their common neighbor in $G'$, then $x_4v_4v_5x_5zx_4$ is a $5$-cycle in $G$. In both cases, we contradict the girth at least~$6$ condition in $G$.~\smallqed

\medskip
By Claim~\ref{claim.no-6-cycle}, the graph $G$ has no $6$-cycle. Let $u$ and $v$ be adjacent vertices in $G$, and let $N(u) = \{u_1,u_2,v\}$ and $N(v) = \{u,v_1,v_2\}$. Further, let $N(u_i) = \{u,u_{i1},u_{i2}\}$ and let $N(v_i) = \{v,v_{i1},v_{i2}\}$ for $i \in [2]$. Thus, $G$ contains the subgraph shown in Figure~\ref{fig-end-proof}. Let $X = \{u_{11},u_{12},u_{21},u_{22},v_{11},v_{12},v_{21},v_{22}\}$. Since the graph $G$ has girth at least~$7$, the set $X$ is an independent set. The subgraph shown in Figure~\ref{fig-end-proof} is therefore an induced subgraph of $G$.

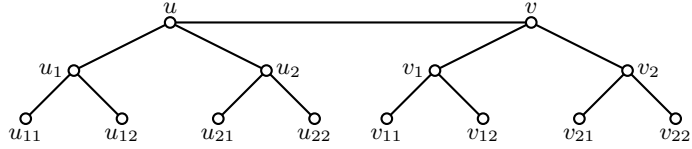
\begin{figure}[htb]
\begin{center}
\begin{tikzpicture}[scale=.8,style=thick,x=0.8cm,y=0.8cm]
\def\vr{2.5pt} 
\path (0,0) coordinate (u11);
\path (1,1) coordinate (u1);
\path (2,0) coordinate (u12);
\path (3,2) coordinate (u);
\path (4,0) coordinate (u21);
\path (5,1) coordinate (u2);
\path (6,0) coordinate (u22);
\path (7.5,0) coordinate (v11);
\path (8.5,1) coordinate (v1);
\path (9.5,0) coordinate (v12);
\path (10.5,2) coordinate (v);
\path (11.5,0) coordinate (v21);
\path (12.5,1) coordinate (v2);
\path (13.5,0) coordinate (v22);
\draw (u11)--(u1)--(u12);
\draw (u21)--(u2)--(u22);
\draw (u1)--(u)--(u2);
\draw (v11)--(v1)--(v12);
\draw (v21)--(v2)--(v22);
\draw (v1)--(v)--(v2);
\draw (u)--(v);
%
\draw (u) [fill=white] circle (\vr);
\draw (u1) [fill=white] circle (\vr);
\draw (u2) [fill=white] circle (\vr);
\draw (u11) [fill=white] circle (\vr);
\draw (u12) [fill=white] circle (\vr);
\draw (u21) [fill=white] circle (\vr);
\draw (u22) [fill=white] circle (\vr);
\draw (v) [fill=white] circle (\vr);
\draw (v1) [fill=white] circle (\vr);
\draw (v2) [fill=white] circle (\vr);
\draw (v11) [fill=white] circle (\vr);
\draw (v12) [fill=white] circle (\vr);
\draw (v21) [fill=white] circle (\vr);
\draw (v22) [fill=white] circle (\vr);
\draw[anchor = south] (u) node {{\small $u$}};
\draw[anchor = east] (u1) node {{\small $u_1$}};
\draw[anchor = west] (u2) node {{\small $u_2$}};
\draw[anchor = north] (u11) node {{\small $u_{11}$}};
\draw[anchor = north] (u12) node {{\small $u_{12}$}};
\draw[anchor = north] (u21) node {{\small $u_{21}$}};
\draw[anchor = north] (u22) node {{\small $u_{22}$}};
\draw[anchor = south] (v) node {{\small $v$}};
\draw[anchor = east] (v1) node {{\small $v_1$}};
\draw[anchor = west] (v2) node {{\small $v_2$}};
\draw[anchor = north] (v11) node {{\small $v_{11}$}};
\draw[anchor = north] (v12) node {{\small $v_{12}$}};
\draw[anchor = north] (v21) node {{\small $v_{21}$}};
\draw[anchor = north] (v22) node {{\small $v_{22}$}};
\end{tikzpicture}
\end{center}
\begin{center}
\vskip -0.5 cm
\caption{A subgraph in the graph $G$}
\label{fig-end-proof}
\end{center}
\end{figure}

Let $Q = \{u,u_1,u_2,v,v_1,v_2\}$ and let $G'$ be obtained from $G - Q$ be adding the edges $e = u_{12}u_{21}$ and $f = v_{12}v_{21}$. Thus, $G'$ is a special subcubic graph that contains exactly four small vertices, namely the vertices in the set $X' = \{u_{11},u_{22},v_{11},v_{22}\}$.  Let $S'$ be a $\gamma_r$-set of $G'$, and let $S = S' \cup \{u^*,v^*\}$ where the vertices $u^*$ and $v^*$ are defined as follows. If $u_{12} \in S'$, let $u^* = u_2$. If $u_{12} \notin S'$ and $u_{21} \in S'$, let $u^* = u_1$. If $u_{12} \notin S'$ and $u_{21} \notin S'$, let $u^* = u$. If $v_{12} \in S'$, let $v^* = v_2$. If $v_{12} \notin S'$ and $v_{21} \in S'$, let $v^* = v_1$. If $v_{12} \notin S'$ and $v_{21} \notin S'$, let $v^* = v$. The resulting set $S$ is a RD-set of $G$, and so $\gamma_r(G) \le \gamma_r(G') + 2$, implying that $\w(G) < \w(G') + 20$.

If $G'$ has no component in $\cB_{\rdom}$, then $\w(G) = \w(G') + 20$, a contradiction. Hence, $G'$ contains a component $G_1$ that belongs to $\cB_{\rdom}$. By Claim~\ref{claim:edge-cut}, there is only one such component and $G_1 \in \{R_2,R_4,R_5,R_{10}\}$. Necessarily, $G_1$ contains at least three vertices from the set $X'$. As observed earlier, the set $X$ is an independent set, and therefore so too is the subset $X'$ of $X$, implying that $G_1 \notin \{R_2,R_5\}$. Every cycle of length less than~$7$ in $G_1$ must contain at least one of the added edges $e$ and $f$ since the graph $G$ has girth at least~$7$. If $G_1$ contains the edge $e$ or $f$, then both ends of the added edge have degree~$3$ in $G_1$. From these properties of the graph $G'$, we deduce that if $G_1 = R_4$, then $G' = G_1$. But this would imply that $G[X] = C_8$, contradicting our earlier observation that $X$ is an independent set. Hence, $G_1 = R_{10}$. In this case, both  added edges $e$ and $f$ must belong to $G_1$. However removing any two edges from $R_{10}$ creates a graph which still contains a $5$-cycle. This implies that $G$ itself contains a $5$-cycle, which is a contradiction. This final contradiction concludes the proof of Theorem~\ref{thm:main-1}.~\QED

\section{Proof of main result}
\label{S:proof-main}

In this section, we prove our main result, namely Theorem~\ref{rdom:main-1}. As a consequence of key result, namely Theorem~\ref{thm:main-1}, we have the following upper bound on the restrained domination number of a cubic graph.

\begin{theorem}
\label{thm:main-2}
If $G$ is a cubic graph of order~$n$, then $\gamma_r(G) \le \frac{2}{5}n$.
\end{theorem}
\proof Let $G$ be a cubic graph of order~$n$. Thus, $n_2(G) = 0$ and $n_3(G) = n$. Since every graph in the family $\cB_{\rdom}$ contains a vertex of degree~$2$, no component of $G$ belongs to the family $\cB_{\rdom}$. The weight of $G$ is therefore $\w(G) = 4n$. Hence by Theorem~\ref{thm:main-1}, $10\gamma_r(G) \le \w(G) = 4n$, or, equivalently, $\gamma_r(G) \le \frac{2}{5}n$.~\QED

\medskip
By Theorem~\ref{thm:main-2}, $c_{\rdom} \le \frac{2}{5}$. As observed earlier, the Petersen graph shows that $c_{\rdom} \ge \frac{2}{5}$. Consequently, $c_{\rdom} = \frac{2}{5}$, yielding the result of Theorem~\ref{rdom:main-1}.
We remark that a classical result in domination theory due to Blank~\cite{Bl73} and McCuaig and Shepherd~\cite{McSh-89} states that if $G$ is a connected graph of order~$n \ge 8$ with $\delta(G) \ge 2$, then $\gamma(G) \le \frac{2}{5}n$. Hence by Theorem~\ref{thm:main-2} this $\frac{2}{5}$-bound for domination also holds for restrained domination if we replace the minimum degree at least~$2$ requirement with a $3$-regularity condition.

\section*{Acknowledgement}
B.B. acknowledges the financial support from the Slovenian Research Agency (research core funding No.\ P1-0297 and project grants J1-2452, J1-3002, and J1-4008).  
The authors thank Jesse Geneson for pointing out that the graph $R_6$ (see  Fig.~\ref{rdom:fig-2}) was missing from the original family $\cB_{\rdom}$.

\end{document}